\RequirePackage{fix-cm}
\PassOptionsToPackage{numbers,sort&compress}{natbib}
\documentclass[smallextended,runningheads,envcountsect,envcountsame,natbib]{svjour3}

\usepackage{amsmath,amssymb,mathtools}
\usepackage{enumitem,microtype}
\usepackage{mathrsfs,comment}
\usepackage{graphicx}
\usepackage{algorithm}
\usepackage{xcolor}
\usepackage{algpseudocode}
\usepackage{booktabs,multirow,tabularx,array}
\usepackage[hidelinks]{hyperref}
\setcitestyle{square,numbers,comma}

\numberwithin{equation}{section}
\spnewtheorem{assumption}[theorem]{Assumption}{\bfseries}{\itshape}
\smartqed
\let\endclassproof\endproof
\renewcommand{\endproof}{\leavevmode\qed\endclassproof}

\DeclareMathOperator{\prox}{Prox}
\DeclareMathOperator{\supp}{supp}

\DeclareMathOperator*{\dom}{dom}

\DeclareMathOperator{\ri}{ri}

\newcommand{\h}[1]{\mathbf{#1}}

\DeclareMathOperator{\sign}{sign}

\newcommand{\x}{\mathbf{x}}
\newcommand{\y}{\mathbf{y}}
\newcommand{\z}{\mathbf{z}}
\newcommand{\g}{\mathbf{g}}

\newcommand{\M}{\mathcal{M}}

\newcommand{\nn}{\nonumber}

\newcommand{\R}{\mathbb{R}}

\def\cB{{\mathscr{B}}}

\hypersetup{
  pdftitle={KL-Guided Support Selection and Subspace Cubic Newton Methods for Sparse Optimization},
  pdfauthor={Zihao Xia, Min Tao, and Guojian Yin},
  pdfkeywords={Nonconvex composite optimization, Subspace cubic regularized Newton, Support identification, Proximal gradient, Global convergence, Kurdyka--\L ojasiewicz}
}

\journalname{Journal of Scientific Computing}

\begin{document}

\title{KL-Guided Support Selection and Subspace Cubic Newton Methods for Sparse Optimization}
\titlerunning{KL-Guided Subspace Cubic Newton Methods}
\author{Zihao Xia \and Min Tao \and Guojian Yin}
\institute{Zihao Xia \at
  School of Mathematics, Nanjing University, Nanjing 210093, China
  \and
  Min Tao (corresponding author) \at
  School of Mathematics, Nanjing University, Nanjing 210093, China\\
  National Key Laboratory for Novel Software Technology, Nanjing 210093, China\\
  \email{\href{mailto:taom@nju.edu.cn}{taom@nju.edu.cn}}
  \and Guojian Yin \at
  Institute for Advanced Study, Shenzhen University, Shenzhen 518060, China
}
\date{}
\maketitle

\begin{abstract}
We study sparse composite optimization with a smooth loss and a possibly nonconvex, nonsmooth regularizer. The objective is locally smooth on each fixed support. Our method combines full-space proximal-gradient steps with inexact cubic-regularized Newton steps on a subspace selected by thresholding. If the iterates are bounded and the objective satisfies the Kurdyka--\L ojasiewicz (KL) property, the entire sequence converges to a first-order stationary point.
A threshold based on the KL exponent ensures that the selected coordinate set coincides with the support of the limit after finitely many iterations, without requiring a nondegeneracy assumption. We also provide a residual-based threshold when such an exponent exists and its value is unknown. Under an additional continuity assumption on the restricted Hessian, the line search eventually accepts full cubic steps. With bounded cubic regularization parameters and damping parameters chosen to meet the target tolerances, we bound the number of outer iterations needed to reach restricted approximate second-order stationarity after support identification and entry into the full-step regime. We conduct numerical experiments in two settings: when the KL exponent is known and when it is unknown. In the former case, we use a degenerate problem to demonstrate the effectiveness of the proposed method. In the latter case, we compare the proposed method and its adaptive variant with existing solvers on sparse recovery and logistic regression problems.
\end{abstract}

\section{Introduction}
	
Sparse optimization models often combine a smooth loss with a regularizer that is nonsmooth in the full space. However, many sparsity-promoting regularizers yield a twice continuously differentiable objective when restricted to a fixed support. We combine full-space proximal-gradient steps with cubic-regularized Newton steps on a selected coordinate subspace. A central question is
{\it whether the selected coordinate set eventually coincides with the support of the limit {point} when the usual nondegeneracy condition fails}. We prove that finite support identification still holds under a suitable thresholding rule determined by the positive power-type KL exponent of the objective at the limit.
	
	In this paper, we consider composite problems of the form
	\begin{equation}
		\min_{\mathbf{x}\in\mathbb{R}^n}
		F(\mathbf{x})
		:=
		g(\mathbf{x})+h(\mathbf{x}),
		\label{Problem}
	\end{equation}
	where
	$g:\mathbb{R}^{n}\rightarrow\overline{\mathbb{R}}
	:=\mathbb{R}\cup\{+\infty\}$
	is a proper, lower semicontinuous (l.s.c.),
	possibly nonconvex and nonsmooth function
	that is prox-bounded with threshold $\lambda_g>0$
	and whose proximal mapping $\prox_{\gamma g}$ can be computed efficiently for $0<\gamma<\lambda_g$.
	We use the following standing assumptions.
	
\begin{assumption}\phantomsection\label{Ass1}
		\begin{enumerate}
			\item[(i)]
			For every nonzero $\mathbf x\in\operatorname{dom}F$, let
			\[
			\mathcal M_{\mathbf x}=\{\mathbf z\in\mathbb R^n:
			\operatorname{supp}(\mathbf z)=\operatorname{supp}(\mathbf x)\}.
			\]
			
			The restriction of $F$ to $\mathcal M_{\mathbf x}$ is finite-valued and twice continuously differentiable in a relative neighborhood of $\mathbf x$. For every nonempty $\mathcal S\subseteq \operatorname{supp}(\mathbf x)$, let $\nabla_{\mathcal S}F$ and $\nabla^2_{\mathcal S}F$ denote the corresponding coordinate gradient and principal Hessian of this restriction, respectively. The mapping $\nabla^2_{\mathcal S}F$ is strictly continuous at $\mathbf x$ relative to $\mathcal M_{\mathbf x}$, with the constant
			$
			\operatorname{lip}_{\mathcal M_{\mathbf x}}
			\nabla^2_{\mathcal S}F(\mathbf x);
			$
			see Definition~\ref{def:Lipschitz}.

				\item[(ii)] $h:\mathbb{R}^{n}\rightarrow\mathbb{R}$
				is continuously differentiable with an $L$-Lipschitz continuous gradient for some $L\geq0$.
			
			\item[(iii)]
			The objective function is bounded below, i.e.,
		$
			\underline F
			:=
			\inf_{\mathbf{x}\in\mathbb{R}^{n}}
			F(\mathbf{x})
			>
			-\infty.
			$
		\end{enumerate}
	\end{assumption}

We illustrate the local smoothness condition in
Assumption~\ref{Ass1}(i) with two models and give bounds for their local restricted Hessian Lipschitz constant.

\begin{enumerate}
\item[(a)] \textbf{Sparse recovery.}
	For $A\in\mathbb{R}^{m\times n}$,
$\mathbf{b}\in\mathbb{R}^m$, and $\lambda>0$,
consider problem~\eqref{Problem} with $h({\h x})=\frac12\|A\mathbf{x}-\mathbf{b}\|_2^2$.
\begin{enumerate}
	\item[(i)] \textbf{$L_1/L_2$ regularization.}
		Let
	$
	g(\mathbf{x})
	=\lambda{\|\mathbf{x}\|_1}/{\|\mathbf{x}\|_2};
	$
	see \citet{Tao20}.
	We define
	$\|\mathbf{0}\|_1/\|\mathbf{0}\|_2=1$
	to make $g$ lower semicontinuous at the
	origin.
	For $\mathbf{x}\neq\mathbf{0}$ and
	$\emptyset\neq\mathcal{S}
	\subseteq\supp(\mathbf{x})$,
	its restricted Hessian satisfies
	\[
	\mathrm{lip}_{\mathcal{M}_{\mathbf{x}}}
	\nabla^2_{\mathcal{S}}F(\mathbf{x})
	\leq
	\frac{36\lambda\sqrt{\|\mathbf{x}\|_0}}
	{\|\mathbf{x}\|_2^3}.
	\]
	\item[(ii)] \textbf{$L_p^p$ regularization.}
	Let
	$
	g(\mathbf{x})=\lambda\|\mathbf{x}\|_p^p, 0<p<1;
	$
	see \cite{xu2012half}.
	For $\mathbf{x}\neq\mathbf{0}$ and
	$\emptyset\neq\mathcal{S}
	\subseteq\supp(\mathbf{x})$,
	its restricted Hessian satisfies
	\[
	\mathrm{lip}_{\mathcal{M}_{\mathbf{x}}}
	\nabla^2_{\mathcal{S}}F(\mathbf{x})
	\leq
	\frac{\lambda p(1-p)(2-p)}
	{(\min_{i\in\mathcal{S}}|x_i|)^{3-p}}.
	\]
\end{enumerate}

		\item[(b)] \textbf{$L_1-\mu L_2$-regularized logistic regression.}
			Let the rows of $A\in\mathbb{R}^{m\times n}$ be
			$\mathbf{a}_i^\top$, and let $b_i\in\{-1,1\}$.
			For $\lambda>0$ and $0<\mu\leq1$, consider \citep{qin20191}
			\begin{equation}
				\min_{(\mathbf{z},u)\in\mathbb{R}^n\times\mathbb{R}}
				\lambda(\|\mathbf{z}\|_1-\mu\|\mathbf{z}\|_2)
				+\frac1m\sum_{i=1}^m
				\log\bigl(1+e^{-b_i(\mathbf{a}_i^\top\mathbf{z}+u)}\bigr).
				\label{eq:intro_sparse_logistic}
			\end{equation}
			Here $u$ is an unpenalized intercept. Set $\mathbf{x}=(\mathbf{z},u)$,
			take the regularization term as $g$, and take the logistic loss as $h$.
			For $\mathbf{z}\neq0$, write {${\mathcal{S}_z}=\supp(\mathbf{z})$.}
			For $\emptyset\neq\mathcal{S}\subseteq\supp(\mathbf{x})$,
			{\[
				\begin{aligned}
				\mathrm{lip}\,\nabla^2_{\mathcal S}F(\x)
				&\le
				\begin{cases}
					M_{\mathrm{log}}+\frac{2\lambda\mu}{\sqrt{3}\|\z\|_2^2},&\z\ne0,\\
					{\frac{1}{6\sqrt3}},&\z=0, u\ne0.
				\end{cases},\\
				M_{\mathrm{log}}&=
				\frac{1}{6\sqrt3m}\sum_{i=1}^m
				(\|\mathbf a_{i,\mathcal{S}_z}\|_2^2+1)^{3/2}.
				\end{aligned}
				\]}
	\end{enumerate}
	Table~\ref{tab:regularizers} lists several commonly used sparsity-promoting regularizers that satisfy Assumption~\ref{Ass1}(i).
	\begin{table}[htbp]
		\centering
		\caption{Examples of regularization functions \(g\) for \(\mathbf{x}\neq\mathbf{0}\) and \(\emptyset\neq\mathcal{S}\subseteq\operatorname{supp}(\mathbf{x})\).}
		\label{tab:regularizers}
		%\small
		\setlength{\tabcolsep}{6pt}
		\renewcommand{\arraystretch}{1.5}
		\begin{tabular}{@{}llc@{}}
			\toprule
			Regularizer
			& $g(\mathbf{x}),\;\lambda >0 \;(\lambda_{1,2}>0)$
			& $\operatorname{lip}_{\mathcal M_{\mathbf{x}}}\nabla_{\mathcal S}^2 g(\mathbf{x})$ \\
			\midrule
			
			$L_0$ \cite{blumensath2009iterative}
			& $\lambda\|\mathbf{x}\|_0$
			& $0$ \\
			
			$L_1$ \cite{tibshirani1996regression}
			& $\lambda\|\mathbf{x}\|_1$
			& $0$ \\

			Elastic Net \cite{zou2005regularization}
			& $\lambda_1\|\mathbf{x}\|_1+\frac{\lambda_2}{2}\|\mathbf{x}\|^2$
			& $0$ \\
			
			Exclusive Lasso \cite{lin2024highly}
			& $\lambda\|\mathbf{x}\|_1^2$
			& $0$ \\		
			
			$L_1+\mu L_2\ (\mu>0)$
			\cite{simon2013sparse}
			& $\lambda(\|\mathbf{x}\|_1+\mu\|\mathbf{x}\|_2)$
			& $\leq\frac{2\lambda\mu}{\sqrt{3}\|\mathbf{x}\|^2}$ \\
			
			$L_p^p$ $(0<p<1)$ \cite{wu2023regularized,xu2012half}
			& $\lambda\|\mathbf{x}\|_p^p$
			& $\frac{\lambda p(1-p)(2-p)}{({\min_{i\in\mathcal S}}|x_i|)^{3-p}}$ \\
			
			$L_1-\mu L_2\ (0<\mu\leq 1)$ \cite{lou2018fast}
			& $\lambda(\|\mathbf{x}\|_1-\mu\|\mathbf{x}\|_2)$
			& $\leq\frac{2\lambda\mu}{\sqrt{3}\|\mathbf{x}\|^2}$\\
			
			$L_1/L_2$ \cite{Tao20}
			& $\lambda\frac{\|\mathbf{x}\|_1}{\|\mathbf{x}\|_2}$
			& $\leq\frac{36\lambda\sqrt{\|\mathbf{x}\|_0}}{\|\mathbf{x}\|^3}$ \\
			
			Log-sum \cite{gong2013general}
			& $\lambda\sum_i\log(1+|x_i|/\varepsilon),\ \varepsilon>0$
			& $  \frac{2\lambda}{(\varepsilon+{\min_{i\in\mathcal S}}|x_i|)^3}$\\
			
			Transformed $L_1$ \cite{zhang2017minimization}
			& $\lambda\sum_i\frac{(a+1)|x_i|}{a+|x_i|}, a>0$
			& $\frac{6\lambda a(a+1)}{(a+{\min_{i\in\mathcal S}}|x_i|)^4}$ \\
			
			%			Group Lasso \cite{yuan2006model}
			%			& $\lambda\sum_\ell\|\mathbf{x}_{G_\ell}\|_2$
			%			& $\frac{6\lambda}{r_G^2}$ \\
			%			
			%			Sparse Group Lasso \cite{simon2013sparse}
			%			& $\lambda_1\|\mathbf{x}\|_1+\lambda_2\sum_\ell\|\mathbf{x}_{G_\ell}\|_2$
			%			& $\frac{6\lambda_2}{r_G^2}$ \\
			\bottomrule
		\end{tabular}
	\end{table}

\subsection{Related Work and Motivation}

Proximal-gradient methods combine a gradient evaluation of $h$
with a proximal evaluation of $g$, without requiring prior
knowledge of the relevant support \citep{beck2017first}. To accelerate these methods, \citet{themelis2018forward} propose
ZeroFPR, which uses quasi-Newton corrections based on the
forward--backward residual and a nonmonotone line search.
For sparse composite problems, hybrid methods exploit
the local smoothness of the objective on a fixed support.
They combine first-order updates with second-order steps on selected
subspaces.

Several algorithms use this structure.
For {$L_p^p$} regularization, {$0<p<1$}, HpgSRN uses signed-support
identification and nonzero-magnitude checks before taking
subspace regularized Newton steps \citep{wu2023regularized}.
IReNA combines reweighted {$L_1$} updates with inexact subspace
Newton steps for separable concave penalties, imposing additional
gradient restrictions beyond power penalties
\citep{wang2024alternating}.
FaRSA-Group combines group proximal-gradient and
conjugate-gradient-based subspace Newton steps, with analysis
under strong convexity and a strict dual margin on inactive
groups \citep{curtis2022subspace}.
For $L_1/L_2$ minimization, HAFAM combines $\mathrm{ADMM}_p$
or $\mathrm{ADMM}_p^{+}$ with globalized semismooth Newton
steps, using a relative-interior condition for identification
\citep{tao2024partlysmooth}.
Optimality-based coordinate selection is also used in
reduced-space methods for {$L_1$}-regularized convex problems
\citep{chen2017reduced}.

For these methods, the choice of subspace determines whether the
selected coordinates eventually coincide with the limiting support. Finite identification is commonly established through
partial smoothness and nondegeneracy
\citep{lewis2002active,lewis2013partial}.
Related results concern active constraints under partial
smoothness and prox-regularity \citep{hare2004identifying},
proximal linearized methods \citep{lewis2016proximal}, and
connections with Newton-type acceleration
\citep{lewis2021active}.
For forward--backward methods, \citet{liang2017activity}
establish finite active-manifold identification and local
linear convergence, while \citet{nutini2019active} derive
identification-complexity bounds for proximal-gradient methods.

Manifold-based methods also use identification to accelerate the iteration. For convex partly smooth
regularizers, ISQA$^{+}$ combines successive quadratic
approximation with semismooth Newton or Riemannian quasi-Newton
steps \citep{lee2023accelerating}. ManAcc combines
proximal-gradient and Riemannian Newton-type steps for
potentially nonconvex partly smooth regularizers
\citep{bareilles2023newton}.
Both use a nondegeneracy condition to establish identification,
so these guarantees do not cover limits at which that condition
fails. Identification without nondegeneracy is nevertheless
possible for particular penalties: \citet{tao2026complexity}
establish identification-complexity bounds for $L_0$ and $L_p$
regularization. This motivates the study of support-selection
rules that avoid nondegeneracy in the composite setting
considered here.

{Besides selecting a suitable subspace, the second-order method should also handle the potential nonconvexity of the resulting model. Cubic regularization offers a natural way to address this challenge.
Under suitable smoothness assumptions and inexact subproblem-solution criteria,}
cubic-regularized Newton methods for smooth nonconvex optimization
achieve an $\mathcal O(\varepsilon^{-3/2})$ iteration bound
for a gradient tolerance $\varepsilon$ and a curvature tolerance
proportional to $\sqrt{\varepsilon}$
\citep{nesterov2006cubic,cartis2011adaptive,cartis2011adaptive2}.
Cubic Newton steps on random coordinate subspaces have also
been studied by \citet{zhao2024cubic}.
These results motivate combining cubic regularization with
a coordinate-selection rule designed to identify the limiting
support.

Our method uses full-space proximal-gradient steps and inexact
cubic-regularized Newton steps on a coordinate subspace. A threshold
guided by the Kurdyka--\L ojasiewicz property determines the selected
coordinates. We study when this coordinate set identifies the limiting
support without a nondegeneracy assumption. We also establish
whole-sequence convergence to a first-order stationary point and derive
an iteration complexity bound for reaching an approximate restricted
second-order stationary point once the support has been identified and
the line search accepts unit step sizes.
	\subsection{Main results}
We propose the Hybrid Proximal-Gradient and Subspace Cubic-Newton (HPG-SCN) method and its adaptive variant (AHPG-SCN). At iteration $k$, a proximal-gradient step produces $\mathbf{x}^k$, and the method selects the coordinates
\[
\mathcal S^k:=\{i:|x_i^k|\geq\varepsilon_k\}.
\]
It then computes an inexact cubic-regularized Newton direction on this coordinate subspace and uses backtracking on the original objective to accept a step.
Our main results are as follows.
\begin{enumerate}
\item We establish sufficient decrease and a relative-error estimate
for the composite iteration. If the iterates are bounded and the
objective has the Kurdyka--\L ojasiewicz (KL) property, the iterate
sequence has finite length and converges to a first-order stationary
point. A threshold rule guided by the KL exponent ensures finite
support identification without a nondegeneracy assumption.
A residual-based logarithmic threshold gives the same conclusion
when the KL exponent is unknown.

\item Under an additional continuity assumption on the restricted
Hessian, the line search eventually accepts full cubic steps.
With bounded cubic regularization parameters and damping chosen
for the target tolerances, Theorem~\ref{thm:local_complexity}
gives an explicit outer-iteration count after support identification
and entry into the full-step regime.

\item {
	We develop AHPG-SCN, which adaptively updates the cubic regularization parameter by comparing actual and predicted reductions. Numerical experiments investigate support identification without nondegeneracy and compare HPG-SCN and AHPG-SCN with several existing methods, demonstrating the computational efficiency of the proposed algorithms.
}
\end{enumerate}

\paragraph{Organization.}
Section~\ref{sec:prelim} introduces notation and preliminary results.
Section~\ref{sec:algorithm} presents HPG-SCN and the Lanczos solver
for the cubic subproblem. Section~\ref{sec:global} gives the global
convergence analysis. Section~\ref{sec:finite} analyzes finite support
identification for the power and logarithmic thresholds.
Section~\ref{sec:local} studies restricted second-order stationarity
and bounds the outer-iteration count in the full-step regime.
Section~\ref{sec:numerical} introduces AHPG-SCN and reports the
numerical experiments. Section~\ref{sec:conclusion} summarizes the
results. Appendices~\ref{app:LHLC} and~\ref{sec:boundness_sequence}
give bounds for the local restricted Hessian Lipschitz modulus and
sufficient conditions for boundedness, respectively.

\section{Preliminaries}\label{sec:prelim}
	Vector norms are Euclidean and matrix norms are spectral unless stated otherwise. We write $\mathbb N=\{0,1,2,\ldots\}$ and $[m]=\{1,\ldots,m\}$ for a positive integer $m$.
	For any $r>0$, let $\cB({\widehat{\mathbf{x}}},r):=\{\mathbf{x}:\|\mathbf{x}-{\widehat{\mathbf{x}}}\|\leq r\}$.
	For an extended-real-valued function $f:\mathbb R^n\to(-\infty,+\infty]$, define $\operatorname{dom}f:=\{\mathbf{x}:f(\mathbf{x})<+\infty\}$. The function is proper if $\operatorname{dom}f\neq\varnothing$.
	It is lower semicontinuous if
	$
	f(\overline{\mathbf x})
	\leq
	\liminf_{\mathbf x\to\overline{\mathbf x}}f(\mathbf x)
	\;
	\text{for every }\overline{\mathbf x}\in\mathbb R^n.$
	The symbol $\sharp(\mathcal{D})$ denotes the cardinality of $\mathcal{D}$.
	$\lambda_{\min}(\cdot)$ denotes the smallest eigenvalue, and ${\text{ker}}(\cdot)$ denotes the null space. $\mathbb S^{n\times n}$ denotes the set of $n\times n$ symmetric matrices.
	\noindent For a proper, lower-semicontinuous function $f$, the proximal mapping with parameter $\alpha>0$ is defined by
	$
		\operatorname{prox}_{\alpha f}(\mathbf{z}) = \mathop{\arg\min}_{\mathbf{x}\in\mathbb{R}^{n}}\left(  f(\mathbf{x}) + \frac{1}{2\alpha}\|\mathbf{x} - \mathbf{z}\|^2\right).
		$
	\begin{definition}
		Let $f:\mathbb{R}^n\to\overline{\mathbb{R}}$ be a proper lower semicontinuous function, and let ${\h x}\in\operatorname{dom}f$.
		\begin{enumerate}
			\item[(a)] A vector $\mathbf v$ is a regular (Fr\'echet) subgradient of $f$ at ${\h x}$, written $\mathbf v\in\widehat{\partial}f({\h x})$, if
			$
			f({\h y})\ge f({\h x})+\langle \mathbf{v},{\h y}-{\h x}\rangle
			+o(\|{\h y}-{\h x}\|).
			$ The little-$o$ relation is understood as ${\h y}\to{\h x}$.
			
		\item[(b)] A vector $\mathbf v$ belongs to the limiting subdifferential $\partial f(\mathbf{x})$ if there exist sequences
			${\h x}_k\to{\h x}$ and $\mathbf{v}_k\to\mathbf{v}$ such that
			${\h x}_k\in\operatorname{dom}f$, $\mathbf{v}_k\in\widehat{\partial}f({\h x}_k)$, and
			$
			f({\h x}_k)\to f({\h x}).
			$
		\end{enumerate}
	\end{definition}
	We say $\mathbf{x}\in\operatorname{dom} f$ is a first-order stationary point of $f$ if $\mathbf{0}\in\partial f(\mathbf{x})$.
We use the KL property and exponent as in \citet{attouch2009convergence}, and strict continuity as in \citet[Definition~9.1]{rockafellar1998variational}.
	\begin{definition}
		\label{def:KL}
		A proper function $f: \mathbb R^n\rightarrow {\overline{\mathbb R}}$ is said to satisfy the {\it Kurdyka-\L ojasiewicz (KL) property} at a point ${\widehat{\mathbf{x}}}\in {\text{\rm dom }}\partial f$ if there exist a constant $\eta\in(0,+\infty]$, an open neighborhood
		$U$ of ${\widehat{\mathbf{x}}},$ and a concave and continuous function $\phi:\;[0,\eta)\rightarrow[0,+\infty)$ with the properties
		\begin{itemize}
			\item[(i)] $\phi(0)=0$;
			\item[(ii)] $\phi$ is continuously differentiable on $(0,\eta)$ and $\phi'{(\cdot)}>0$ on $(0,\eta)$;
			\item[(iii)] for every $\mathbf{x}\in U$ with $f({\widehat{\mathbf{x}}})< f(\mathbf{x})< f(\widehat{\mathbf{x}}) +\eta$, it holds
			\begin{equation*}
				\phi'(f(\mathbf{x})-f(\widehat{\mathbf{x}})){\text{\rm{dist}}}(\mathbf{0},\partial f({\mathbf{x}}))\ge 1.
			\end{equation*}
		\end{itemize}
		The class of desingularization functions $\phi:[0,\eta)\rightarrow[0,+\infty)$ satisfying (i)-(iii) is denoted by $\Phi_\eta$.
	\end{definition}

	\begin{definition}\label{def:KL_exponent}
		Let $f:\mathbb{R}^n\to(-\infty,+\infty]$ be proper and lower semicontinuous, and
		let $\overline \x\in\mathrm{dom}(\partial f)$ be a point satisfying $0\in\partial f(\overline \x)$.
		We say that $f$ satisfies the \emph{KL property at  $\overline \x$ with exponent
			$\theta\in[0,1)$} if there exist constants $\xi>0$, $\varepsilon>0$ and $\nu>0$ such that
		\begin{equation*}\label{eq:KL_basic}
				\mathrm{dist}({{\bf 0}},\partial f(\x))\ge\nu\,(f(\x)-f(\overline \x))^{\theta},
		\end{equation*}
		whenever $\|\x-\overline \x\|<\varepsilon$ and
			$0<f(\x)-f(\overline \x)<\xi$.
	\end{definition}
	
	\begin{definition}\label{def:Lipschitz}
		Let $F:D\subseteq\mathbb{R}^n\to\mathbb{R}^{m\times s}$ be a single-valued mapping, and let $X\subseteq D$. We say that $F$ is \emph{strictly continuous} at $\overline{\x}$ relative to $X$ if $\overline{\x}\in X$ and the quantity
		\begin{equation*}
			\operatorname{lip}_{X}F(\overline{\x})
			:=
			\limsup_{\substack{\x,\x'\to\overline{\x}\\
					\x,\x'\in X,\ \x\neq\x'}}
			\frac{\|F(\x')-F(\x)\|}{\|\x-\x'\|}
		\end{equation*}
		is finite. This quantity is also called the Lipschitz constant of $F$ at $\overline{\x}$ relative to $X$.
	\end{definition}
	
	Throughout the paper, we use
\(\operatorname{lip}_{\mathcal M}\nabla^2 F(\x)\) to denote
the Lipschitz constant of \(\nabla^2 F\) at \(\x\) relative to
\(\mathcal M\). When no ambiguity arises, we write
\(\operatorname{lip}\nabla^2 F(\x)\).
Proposition~\ref{thm:cubic_solution} states the global optimality
conditions for the cubic subproblem. For a proof, see
\citet[Proposition~1]{nesterov2006cubic} and \citet{jia2022solving}. %Appendix~\ref{prop1proof} gives the proof.

\begin{proposition}\label{thm:cubic_solution}
			Let $\widehat H\in\mathbb S^{n\times n}$,
			$\mathbf g\in\mathbb R^n$, and $\sigma>0$.
			Consider the following cubic-regularized Newton problem:
			\begin{equation}\label{prob:cubic_newton_sub}
				\min_{\mathbf d\in\mathbb R^n}
				\left\{
				\langle\mathbf g,\mathbf d\rangle
				+\frac12\langle\widehat H\mathbf d,\mathbf d\rangle
				+\frac{\sigma}{6}\|\mathbf d\|^3
				\right\}.
			\end{equation}
			Let $\lambda_1$ denote the smallest eigenvalue of
			$\widehat H$.
			Then there exists $\lambda^\star\geq 0$ such that
			every global minimizer $\mathbf d^*$ satisfies
			\begin{equation}\label{eq:cubic-sol}
				(\widehat H+\lambda^\star I)\mathbf d^*=-\mathbf g,
				\quad
				\widehat H+\lambda^\star I\succeq0,
				\quad
				\lambda^\star=\frac{\sigma}{2}\|\mathbf d^*\|.
			\end{equation}
			Moreover, for any unit eigenvector $\mathbf v_1$
			associated with $\lambda_1$, the following expression
			gives a global minimizer $\mathbf d^*$ of
			\eqref{prob:cubic_newton_sub}:
			\begin{equation*}\label{eq:cubic-sol-rep}
				\mathbf d^*=
				\begin{cases}
					-(\widehat H+\lambda^\star I)^{-1}\mathbf g,
					& \text{if }\lambda^\star>-\lambda_1,\\[4pt]
					-(\widehat H+\lambda^\star I)^\dagger\mathbf g
					+\alpha\,\mathbf v_1,
					& \text{if }\lambda^\star=-\lambda_1,
				\end{cases}
			\end{equation*}
			where, in the second case, $\alpha\in\mathbb R$
			satisfies
			\begin{equation*}
				\bigl\|
				-(\widehat H+\lambda^\star I)^\dagger\mathbf g
				+\alpha\,\mathbf v_1
				\bigr\|
				=\frac{2\lambda^\star}{\sigma}.
			\end{equation*}
			Here, $(\cdot)^\dagger$ denotes the Moore--Penrose
			pseudoinverse.
	\end{proposition}

	\section{The Hybrid Subspace Cubic-Newton Method}\label{sec:algorithm}

Algorithm~\ref{alg:apgcn} combines a full-space proximal-gradient step
with an inexact cubic-regularized Newton step on a selected coordinate
subspace. At each iteration, the proximal-gradient step generates
$\mathbf x^k$ from $\mathbf y^{k-1}$. The subspace step then updates
selected coordinates to obtain $\mathbf y^k$.
{Both steps employ line searches to ensure sufficient decrease in the original objective function $F$.}

The proximal-gradient line search starts from
\begin{equation}\label{eq:adaptive_initial_stepsize}
			\alpha_{k,0}=
			\begin{cases}
				\min\left\{\dfrac{\|\mathbf{y}^{k-1}-\mathbf{y}^{k-2}\|}
				{\|\nabla h(\mathbf{y}^{k-1})-\nabla h(\mathbf{y}^{k-2})\|},\ \overline{\alpha}\right\}
				& \begin{aligned}
				  &\text{if }k\geq 2\text{ and}\\[-1pt]
				  &\nabla h(\mathbf{y}^{k-1})\neq\nabla h(\mathbf{y}^{k-2})
				  \end{aligned}\\
				\alpha_0,
				& \text{otherwise},
			\end{cases}
		\end{equation}
		where $0<\alpha_0 <\overline{\alpha}<\lambda_g$.
The trial stepsize is reduced by a factor $q_1\in(0,1)$ until the
proximal-gradient point satisfies the prescribed decrease condition.
For positive trial stepsizes, the bound
$\alpha_{k,l_k}\leq\overline{\alpha}<\lambda_g$ ensures that the
proximal subproblems are well-defined.
At the accepted point $\mathbf x^k$, we select the coordinates
\[
    \mathcal S^k
    :=\{\,i:|x_i^k|\geq\varepsilon_k\,\},
    \qquad
    n_k:=\sharp(\mathcal S^k),
\]
where the thresholds $\varepsilon_k$ are chosen according to the
rules in Section~\ref{sec:finite}.
When $\mathcal S^k\neq\emptyset$,
Assumption~\ref{Ass1}(i) ensures that the restriction of $F$ to the
affine slice through $\mathbf x^k$, with all coordinates outside
$\mathcal S^k$ fixed, is twice continuously differentiable near
$\mathbf x^k$.
We denote its gradient and Hessian at $\mathbf x^k$ by
\[
    \mathbf g^k
    :=\nabla_{\mathcal S^k}F(\mathbf x^k),
    \qquad
    H^k
    :=\nabla^2_{\mathcal S^k}F(\mathbf x^k).
\]
The subspace step $\mathbf d^k\in\mathbb R^{n_k}$ is an approximate
solution of the cubic-regularized Newton subproblem
\begin{equation}\label{subd}
    \min_{\mathbf d\in\mathbb R^{n_k}}
    \left\{
        \langle\mathbf g^k,\mathbf d\rangle
        +\frac12
        \langle(H^k+\eta I)\mathbf d,\mathbf d\rangle
        +\frac{\sigma_k}{6}\|\mathbf d\|^3
    \right\}.
\end{equation}
Here, $\eta>0$ provides quadratic regularization, and $\sigma_k$
is chosen using a bound on the local variation of the restricted
Hessian, as specified in Algorithm~\ref{alg:apgcn}.

To specify the required accuracy of the subspace step, define the
first-order residual
\[
    \mathbf r^k
    :=
    \mathbf g^k+(H^k+\eta I)\mathbf d^k
    +\frac{\sigma_k}{2}\|\mathbf d^k\|\mathbf d^k.
\]
We require $\mathbf d^k$ to satisfy
\begin{subequations}\label{cubic_inexact_conditions}
    \begin{align}
        \langle\mathbf r^k,\mathbf d^k\rangle
        &\leq 0,
        \label{inexact1}\\
        \left\langle
        \left(
            H^k+\eta I
            +\frac{\sigma_k}{2}\|\mathbf d^k\|I
        \right)\mathbf d^k,
        \mathbf d^k
        \right\rangle
        &\geq 0,
        \label{inexact2}\\
        \|\mathbf r^k\|
        &\leq\kappa_g\|\mathbf d^k\|^2,
        \label{inexact3}\\
        H^k+\eta I+\frac{\sigma_k}{2}\|\mathbf d^k\|I
        &\succeq-\kappa_H\|\mathbf d^k\|I,
        \label{inexact4}
    \end{align}
\end{subequations}
where $\kappa_g,\kappa_H\geq0$.

A global minimizer of~\eqref{subd} satisfies
\eqref{inexact1}--\eqref{inexact4}, with equality
in~\eqref{inexact1}; see
\cite{cartis2011adaptive,cartis2011adaptive2}.
Condition~\eqref{inexact3} bounds the first-order residual by the
squared step norm, as in inexact regularization methods
\cite{birgin2017worst}.
Condition~\eqref{inexact4} controls the negative curvature of the
shifted Hessian.
Conditions~\eqref{inexact1}--\eqref{inexact2} are used to establish
sufficient decrease in Section~\ref{sec:global}, whereas
\eqref{inexact3}--\eqref{inexact4} are needed for the local complexity
analysis in Section~\ref{sec:local}.

A backtracking line search along $\mathbf d^k$ then determines
$\mathbf y^k$, with all coordinates outside $\mathcal S^k$ fixed.
The acceptance condition contains both quadratic and cubic decrease
terms. If $\mathcal S^k=\emptyset$, we set
$\mathbf y^k=\mathbf x^k$.
Algorithm~\ref{alg:apgcn} gives the complete update.

\begin{algorithm}[htbp]
    \caption{Hybrid proximal-gradient and subspace cubic-Newton
    method (HPG-SCN) for~\eqref{Problem}}
    \label{alg:apgcn}
    \begin{algorithmic}[1]
        \Require{
            $\mathbf x^0=\mathbf y^0\in\operatorname{dom}F$,
            $0<\alpha_0<\overline{\alpha}<\lambda_g$,
            {$c,\ c_1>0$,}
            $\eta>0$,
            {$q,\ q_1\in(0,1)$,}
            $\underline{\omega}>1$,
            $\varepsilon_k>0$,
            {$\kappa_g,\ \kappa_H\geq0$,}
            \texttt{MaxIt}.
        }

        \For{$k=1,2,\ldots,\texttt{MaxIt}$}

            \State \textbf{Proximal-gradient step.}
            \State Set $\alpha_{k,0}$ according
            to~\eqref{eq:adaptive_initial_stepsize}.

            \For{$l_k=0,1,2,\ldots$}
                \State Compute
                \[
                    \mathbf x^{k,l_k}
                    \in
                    \prox_{\alpha_{k,l_k}g}
                    \bigl(
                        \mathbf y^{k-1}
                        -\alpha_{k,l_k}
                         \nabla h(\mathbf y^{k-1})
                    \bigr).
                \]
                \If{$
                    F(\mathbf x^{k,l_k})
                    \leq
                    F(\mathbf y^{k-1})
                    -c_1
                    \|\mathbf x^{k,l_k}-\mathbf y^{k-1}\|^2
                $}
                    \State Set
                    $\mathbf x^k=\mathbf x^{k,l_k}$ and
                    $\alpha_k=\alpha_{k,l_k}$;
                    \textbf{break}.
                \Else
                    \State Set
                    $\alpha_{k,l_k+1}=q_1\alpha_{k,l_k}$.
                \EndIf
            \EndFor

            \State \textbf{Coordinate selection.}
            \State Set
            $\mathcal S^k
            :=\{\,i:|x_i^k|\geq\varepsilon_k\,\}$
            and $n_k:=\sharp(\mathcal S^k)$.

            \If{$\mathcal S^k\neq\emptyset$}

                \State \textbf{Subspace cubic-Newton step.}
                \State Compute
                $\mathbf g^k
                =\nabla_{\mathcal S^k}F(\mathbf x^k)$
                and
                $H^k
                =\nabla^2_{\mathcal S^k}F(\mathbf x^k)$.

                \State Choose $M(\mathbf x^k)>0$ such that
                \[
                    M(\mathbf x^k)
                    \geq
                    \underline{\omega}\,
                    \mathrm{lip}_{\mathcal M_{\mathbf x^k}}
                    \nabla^2_{\mathcal S^k}F(\mathbf x^k).
                \]
                \State Choose
                $\sigma_k\geq2(M(\mathbf x^k)+c)/3$.

                \State Compute an approximate solution
                $\mathbf d^k$ of~\eqref{subd}
                satisfying
                \eqref{inexact1}--\eqref{inexact4}.

                \State \textbf{Backtracking line search.}
                \State Set $\beta_{k,0}=1$.

                \For{$j_k=0,1,2,\ldots$}
                    \State Define
                    \[
                        \mathbf y^{k,j_k}_{\mathcal S^k}
                        =
                        \mathbf x^k_{\mathcal S^k}
                        +\beta_{k,j_k}\mathbf d^k,
                        \qquad
                        \mathbf y^{k,j_k}_{(\mathcal S^k)^c}
                        =
                        \mathbf x^k_{(\mathcal S^k)^c}.
                    \]
                    \If{$
                        F(\mathbf y^{k,j_k})
                        \leq
                        F(\mathbf x^k)
                        -\frac{\eta}{4}
                         \|\mathbf y^{k,j_k}-\mathbf x^k\|^2
                        -\frac{c}{6}
                         \|\mathbf y^{k,j_k}-\mathbf x^k\|^3
                    $}
                        \State Set
                        $\mathbf y^k=\mathbf y^{k,j_k}$
                        and $\beta_k=\beta_{k,j_k}$;
                        \textbf{break}.
                    \Else
                        \State Set
                        $\beta_{k,j_k+1}=q\beta_{k,j_k}$.
                    \EndIf
                \EndFor

            \Else
                \State Set $\mathbf y^k=\mathbf x^k$.
            \EndIf

            \If{the termination criterion is satisfied}
                \State \textbf{break}.
            \EndIf

        \EndFor

        \State \textbf{Output:} $\mathbf y^k$.
    \end{algorithmic}
\end{algorithm}

\subsection{Solving the Cubic Subproblem via a Restarted Lanczos Algorithm}
Methods for solving cubic-regularization subproblems include
Lanczos methods \cite{cartis2011adaptive,jia2022solving},
gradient descent \cite{carmon2019gradient},
first-order methods applied to convex reformulations
\cite{jiang2021accelerated,jiang2022cubic},
generalized eigenvalue methods \cite{lieder2020solving},
and methods based on approximate secular equations
\cite{gao2022approximate}.
We use the symmetric Lanczos process in
\citet[Algorithm~1]{jia2022solving}, with the orthogonal
continuation after exact breakdown described at the end of their
Section~4.1. Their pseudocode stops at breakdown. If breakdown occurs before the
target dimension is reached, we retain the existing
basis $Q_i$ and choose a unit vector
$\mathbf q_{i+1}\in\operatorname{ker}(Q_i^\top)$ to start a new
Lanczos block. This is distinct from their nested restarting
method in Section~4.2: we increase the target dimension until
\eqref{inexact1}--\eqref{inexact4} hold.

For a symmetric $A\in\mathbb R^{N\times N}$, a unit vector
$\mathbf q$, and $1\leq r\leq N$, let
$[Q_r,T_r]=\texttt{Lanczos}(A,\mathbf q,r)$ denote this process
continued until $r$ basis vectors have been generated, i.e.,
\[
Q_r^\top Q_r=I,
\qquad
T_r=Q_r^\top A Q_r,
\]
where $T_r$ is symmetric and tridiagonal, with zero off-diagonal
entries at block boundaries.
Applying this procedure to $A=H^k+\eta I$ and using
$\mathbf d=Q_r\mathbf u$ gives the reduced cubic problem
\begin{equation}\label{reduced_cubic_subproblem}
\min_{\mathbf u\in\mathbb R^r}
\left\{
\langle Q_r^\top\mathbf g^k,\mathbf u\rangle
+\frac12\mathbf u^\top T_r\mathbf u
+\frac{\sigma_k}{6}\|\mathbf u\|^3
\right\}.
\end{equation}
The identity $\|Q_r\mathbf u\|=\|\mathbf u\|$ shows that the reduced problem is the
restriction of~\eqref{subd} to $\operatorname{Range}(Q_r)$.
Its objective is continuous and coercive, so a global minimizer
exists.

Algorithm~\ref{alg:lanczos_cubic} solves the reduced problem and
increases its dimension until the resulting step satisfies
\eqref{inexact1}--\eqref{inexact4}. The additional $1$ in the
dimension update ensures a strict increase before reaching $n_k$,
even when $s=1$.

\begin{algorithm}[!t]
    \caption{A Lanczos method for solving~\eqref{subd}}
    \label{alg:lanczos_cubic}
    \begin{algorithmic}[1]
        \Require Integers $n_k\geq1$ and $1\leq r_{k,0}\leq n_k$,
        $s\geq1$, $\mathbf g^k\in\mathbb R^{n_k}$, a symmetric
        $H^k\in\mathbb R^{n_k\times n_k}$,
        {$\eta,\ \sigma_k>0$, and $\kappa_g,\ \kappa_H\geq0$.}
        \If{$\mathbf g^k\neq\mathbf0$}
            \State Set $\mathbf q^k=\mathbf g^k/\|\mathbf g^k\|$.
        \Else
            \State Choose any unit vector
            $\mathbf q^k\in\mathbb R^{n_k}$.
        \EndIf
        \For{$j=1,2,\ldots$}
            \State Set
            $r_{k,j}=\min\{\lceil s r_{k,j-1}\rceil+1,n_k\}$.
            \State Compute
            $[Q_{k,j},T_{k,j}]
            =\texttt{Lanczos}(H^k+\eta I,\mathbf q^k,r_{k,j})$.
            \State Compute a global minimizer $\mathbf u^{k,j}$
            of~\eqref{reduced_cubic_subproblem} with
            $(r,Q_r,T_r)=(r_{k,j},Q_{k,j},T_{k,j})$.
            \State Set $\mathbf d^{k,j}=Q_{k,j}\mathbf u^{k,j}$
            and compute
            \[
            \mathbf r^{k,j}
            =\mathbf g^k+(H^k+\eta I)\mathbf d^{k,j}
            +\frac{\sigma_k}{2}\|\mathbf d^{k,j}\|\mathbf d^{k,j}.
            \]
            \If{conditions~\eqref{inexact1}--\eqref{inexact4}
                hold for $\mathbf d^{k,j}$}
                \State Set
                {$\mathbf d^k=\mathbf d^{k,j}$},
                %$\mathbf u^k=\mathbf u^{k,j}$,
                %$Q_k=Q_{k,j}$, and $T_k=T_{k,j}$.
                \State \Return {$\mathbf d^k$}.
            \EndIf
        \EndFor
    \end{algorithmic}
\end{algorithm}

\begin{theorem}\label{thm:lanczos_termination}
   {In exact arithmetic, the Lanczos process described above} generates
    $r$ orthonormal vectors for every $1\leq r\leq n_k$.
    If every reduced cubic problem is solved globally,
    Algorithm~\ref{alg:lanczos_cubic} terminates after finitely
    many inner iterations.
\end{theorem}
\begin{proof}
    Within each Lanczos block, the recurrence generates
    orthonormal vectors until breakdown. At breakdown, the
    accumulated subspace is invariant under the symmetric
    matrix $H^k+\eta I$, so its orthogonal complement is also
    invariant. A unit vector in this complement exists whenever
    fewer than $r\leq n_k$ vectors have been generated, and
    starting a new block from it preserves orthogonality.
    Since $s\geq1$, the dimension update reaches $n_k$ after
    finitely many inner iterations unless a step is accepted
    earlier. At full dimension, $Q_{k,j}$ is square and
    orthogonal, so the reduced problem is equivalent
    to~\eqref{subd}. Its global minimizer has zero residual and
    satisfies~\eqref{eq:cubic-sol}; hence, all four acceptance
    conditions hold.
\end{proof}

\begin{remark}
    An exact global minimizer of the reduced
    problem automatically satisfies~\eqref{inexact1}
    and~\eqref{inexact2}. Let
    $\lambda_{k,j}:=\frac{\sigma_k}{2}\|\mathbf d^{k,j}\|$,
    Proposition~\ref{thm:cubic_solution} gives
    \[
    Q_{k,j}^\top\mathbf r^{k,j}=\mathbf0,
    \qquad
    T_{k,j}+\lambda_{k,j}I\succeq0.
    \]
    Since $\mathbf d^{k,j}=Q_{k,j}\mathbf u^{k,j}$, it follows that
    \[
    \langle\mathbf r^{k,j},\mathbf d^{k,j}\rangle=0,
    \]
    and
    \[
    (\mathbf d^{k,j})^\top
    (H^k+\eta I+\lambda_{k,j}I)\mathbf d^{k,j}
    =
    (\mathbf u^{k,j})^\top
    (T_{k,j}+\lambda_{k,j}I)\mathbf u^{k,j}
    \geq0.
    \]
    Thus, only~\eqref{inexact3} and~\eqref{inexact4} require
    separate verification.
    Condition~\eqref{inexact3} is checked using the full residual.
    A numerical certificate for~\eqref{inexact4} can be obtained
    from a lower bound $\ell_k\leq\lambda_{\min}(H^k+\eta I)$
    satisfying
    $\ell_k+(\sigma_k/2+\kappa_H)\|\mathbf d^{k,j}\|\geq0$.
\end{remark}

	\section{Convergence Analysis}\label{sec:global}
	In this section, we establish the global convergence result of Algorithm~\ref{alg:apgcn}.
	\subsection{Descent Property}
		We first show that the line search for
		the proximal-gradient step terminates
		after finitely many trials.
	\begin{theorem}\label{thm:PG_linesearch}
		Assume that $F$ defined in \eqref{Problem}
		satisfies Assumption~\ref{Ass1} and $\alpha_{k,0}$ is chosen according to
		\eqref{eq:adaptive_initial_stepsize}.
		Then, the procedure of finding the smallest
		$l_k\in\{0,1,2,\ldots\}$ such that
		\begin{equation}\label{eq:PG_linesearch}
			F(\mathbf{x}^{k,l_k})
			\leq
			F(\mathbf{y}^{k-1})
			-c_1\|\mathbf{x}^{k,l_k}
			-\mathbf{y}^{k-1}\|^2
		\end{equation}
		terminates after finitely many trials.
		Consequently, the proximal-gradient step
		in Algorithm~\ref{alg:apgcn} is well-defined.
		Moreover, the accepted step sizes satisfy
		$ 0<\underline{\alpha} \leq\alpha_k\leq\overline{\alpha}\;(k\geq1)$ with
		$
		\underline{\alpha}:=
		\min\left\{
		\alpha_0,
		\frac{q_1}{L+2c_1}
		\right\}>0.
		$
	\end{theorem}
	
	\begin{proof}
		By the bound $0<\alpha_{k,l_k}\leq\alpha_{k,0}\leq\overline{\alpha}<\lambda_g$ and the definition of
		$\mathbf{x}^{k,l_k}$, we have
		\begin{equation*}
		\begin{aligned}
		&g(\mathbf{x}^{k,l_k})
			+\frac{1}{2\alpha_{k,l_k}}
			\left\|
			\mathbf{x}^{k,l_k}
			-\mathbf{y}^{k-1}
			+\alpha_{k,l_k}
			\nabla h(\mathbf{y}^{k-1})
			\right\|^2\\
		&\qquad\leq
			g(\mathbf{y}^{k-1})
			+\frac{\alpha_{k,l_k}}{2}
			\|\nabla h(\mathbf{y}^{k-1})\|^2.
		\end{aligned}
		\end{equation*}
		Moreover, since $\nabla h$ is $L$-Lipschitz continuous, the descent lemma gives
		\begin{equation*}
			h(\mathbf{x}^{k,l_k})
			\leq
			h(\mathbf{y}^{k-1})
			+\left\langle
			\nabla h(\mathbf{y}^{k-1}),
			\mathbf{x}^{k,l_k}
			-\mathbf{y}^{k-1}
			\right\rangle +\frac{L}{2}
			\|\mathbf{x}^{k,l_k}
			-\mathbf{y}^{k-1}\|^2.
		\end{equation*}
		Combining the above two inequalities,
		we obtain
		\begin{equation*}\label{eq:PG_descent_linesearch}
				F(\mathbf{x}^{k,l_k})
				\leq
				F(\mathbf{y}^{k-1})-\left(
				\frac{1}{2\alpha_{k,l_k}}
				-\frac{L}{2}
				\right)
				\|\mathbf{x}^{k,l_k}
				-\mathbf{y}^{k-1}\|^2.
		\end{equation*}
		Since $\alpha_{k,l_k}=q_1^{l_k}\alpha_{k,0}$ and $0<q_1<1$, there exists a {smallest finite index $\widehat{l}_k$} such that
		$ \alpha_{k,\widehat{l}_k} \leq{1}/{(L+2c_1)}.
		$
		Therefore, \eqref{eq:PG_linesearch}
		holds for $l_k=\widehat{l}_k$,
		and the proximal-gradient backtracking
		procedure terminates after finitely
		many trials.
If $l_k=0$, then $\alpha_k=\alpha_{k,0}$.
By \eqref{eq:adaptive_initial_stepsize}, the initial trial is either $\alpha_0$
or, when the gradient difference is nonzero, at least
$\min\{1/L,\overline{\alpha}\}$. In the latter case $L>0$, and this minimum
is at least $\underline{\alpha}$.
If $l_k>0$, the preceding trial was rejected, so
$\alpha_{k,l_k-1}>1/(L+2c_1)$. Thus
$\alpha_k=q_1\alpha_{k,l_k-1}>q_1/(L+2c_1)\geq\underline{\alpha}$.
	\end{proof}	
For subspace updates, we adapt the cubic upper bound from \citep[Lemma~1]{nesterov2006cubic} to restricted derivatives and subspace directions. An additional term accounts for error from variation in the restricted Hessian.
	
	\begin{lemma}\label{lem:inexact_cubic_upper}
		%Suppose that $F$ defined in \eqref{Problem} satisfies Assumption~\ref{Ass1}.
		Let $\x,\mathbf{d}\in\R^n$ and $\M:=\{\z\in\R^{n}:\supp(\z)=\supp(\x) \}$. Suppose $\supp(\mathbf{d})\subseteq \supp(\mathbf{x})$ and $ \supp(\mathbf{x} + t\mathbf{d}) = \supp(\x)$ for all $t\in[0,1]$.
		Assume that the restriction of $F$ to $\M$ is finite and twice continuously differentiable in a relative neighborhood of the segment $\{\x+t\mathbf d:0\leq t\leq1\}$. Assume $\mathcal{S}$ is an index set satisfying $\supp(\mathbf{d})\subseteq\mathcal{S}\subseteq \supp(\mathbf{x})$ and there exist constants $M\ge 0$ and $a\ge 0$ such that for all $t\in[0,1]$, the function $F$ satisfies
		\begin{equation}\label{eq:hess_var}
			\bigl\|\nabla_{\mathcal{S}}^2 F(\x+t\mathbf{d})-\nabla_{\mathcal{S}}^2 F(\x)\bigr\|
			\le M\,t\|\mathbf{d}\|+a.
		\end{equation}
		Then, for all $t\in[0,1]$, it holds that
		\begin{equation}\label{ineq:grad_var_bound}
			\bigl\|\nabla_{\mathcal{S}}F(\x+t\mathbf{d})-\nabla_{\mathcal{S}}F(\x)
			- t\,\nabla_{\mathcal{S}}^2F(\x)\mathbf{d}_{\mathcal{S}}\bigr\|
			\le a\,t\|\mathbf{d}\|+\frac{M}{2}\,t^2\|\mathbf{d}\|^2,
		\end{equation}
		where $\mathbf{d}_{\mathcal{S}}\in\mathbb{R}^{\sharp(\mathcal{S})}$
		denotes the subvector of $\mathbf{d}$ consisting of the components
		indexed by $\mathcal{S}$.
		
		\noindent Moreover,
		\begin{equation}\label{ineq:cubic_upper}
			F(\x+\mathbf{d})\le F(\x)
			+\langle \nabla_{\mathcal{S}}F(\x),\mathbf{d}_{\mathcal{S}}\rangle
			+\frac12\langle \nabla_{\mathcal{S}}^2F(\x)\mathbf{d}_{\mathcal{S}},\mathbf{d}_{\mathcal{S}}\rangle
			+\frac{a}{2}\|\mathbf{d}\|^2+\frac{M}{6}\|\mathbf{d}\|^3.
		\end{equation}
	\end{lemma}
	
	\begin{proof}
		For $t\in[0,1]$, define
		$
		\psi(t):=\nabla_{\mathcal{S}}F(\x+t\mathbf{d})\in\R^{\sharp(\mathcal{S})}.
		$
		Since $\supp(\x+t\mathbf{d})=\supp(\x)$ for all $t\in[0,1]$,  $\supp(\mathbf{d})\subseteq\mathcal{S}\subseteq \supp(\mathbf{x})$
		and (\ref{eq:hess_var}) holds for $t\in[0,1]$,
		$\psi$ is differentiable on $[0,1]$, with
		\[
		\psi'(t)=\nabla_{\mathcal{S}}^{2}F(\x+t\mathbf{d})\,\mathbf{d}_{\mathcal{S}}, \quad t\in[0,1].
		\]
		Then, we obtain
		\begin{equation*}
			\begin{aligned}
			\psi(t)-\psi(0)-t\psi'(0)
			&=\int_{0}^{t}\bigl(\psi'(u)-\psi'(0)\bigr)\,du\\
			&=\int_{0}^{t}\bigl(\nabla_{\mathcal{S}}^{2}F(\x+u\mathbf{d})-\nabla_{\mathcal{S}}^{2}F(\x)\bigr)
			\,\mathbf{d}_{\mathcal{S}}\,du .
			\end{aligned}
		\end{equation*}
		Taking norms and using $\|\mathbf{d}_{\mathcal{S}}\| = \|\mathbf{d}\|$, together with
		\eqref{eq:hess_var}, yields
		\begin{equation*}
			\begin{aligned}
				\bigl\|\psi(t)-\psi(0)-t\psi'(0)\bigr\|
				&\le \int_{0}^{t}\bigl\|\nabla_{\mathcal{S}}^{2}F(\x+u\mathbf{d})
				-\nabla_{\mathcal{S}}^{2}F(\x)\bigr\|\,\|\mathbf{d}\|\,du \\
				&\le \int_{0}^{t}\bigl(Mu\|\mathbf{d}\|+a\bigr)\|\mathbf{d}\|\,du
				= a\,t\|\mathbf{d}\|+\frac{M}{2}t^{2}\|\mathbf{d}\|^{2},
			\end{aligned}
		\end{equation*}
		which proves \eqref{ineq:grad_var_bound}.
		Next, define $\phi(t)=F(\x+t\mathbf{d})$. Then, $\phi$ is
		differentiable on $[0,1]$ and satisfies
		$
		\phi'(t)=\langle \nabla_{\mathcal{S}}F(\x+t\mathbf{d}),\mathbf{d}_{\mathcal{S}}\rangle.
		$
		Consequently,
		\begin{equation*}
			\begin{aligned}
				&F(\x+\mathbf{d})-F(\x)-\langle \nabla_{\mathcal{S}}F(\x),\mathbf{d}_{\mathcal{S}}\rangle = \phi(1) - \phi(0) - \phi'(0)\\
				& = \int_{0}^{1}(\phi'(t) - \phi'(0))\,dt\\
				&=\int_{0}^{1}\bigl\langle \nabla_{\mathcal{S}}F(\x+t\mathbf{d})
				-\nabla_{\mathcal{S}}F(\x),\mathbf{d}_{\mathcal{S}}\bigr\rangle\,dt \\
				&=\int_{0}^{1}\Bigl\langle t\nabla_{\mathcal{S}}^{2}F(\x)\mathbf{d}_{\mathcal{S}},
				\mathbf{d}_{\mathcal{S}}\Bigr\rangle\,dt\\
				&\quad+\int_{0}^{1}\Bigl\langle \nabla_{\mathcal{S}}F(\x+t\mathbf{d})
				-\nabla_{\mathcal{S}}F(\x)
				-t\nabla_{\mathcal{S}}^{2}F(\x)\mathbf{d}_{\mathcal{S}},\mathbf{d}_{\mathcal{S}}\Bigr\rangle\,dt \\
				&\le\frac 12\bigl\langle \nabla_{\mathcal{S}}^{2}F(\x)\mathbf{d}_{\mathcal{S}},\mathbf{d}_{\mathcal{S}}\rangle\\
				&\quad+\int_{0}^{1}\Bigl|\Bigl\langle \nabla_{\mathcal{S}}F(\x+t\mathbf{d})
				-\nabla_{\mathcal{S}}F(\x)
				-t\nabla_{\mathcal{S}}^{2}F(\x)\mathbf{d}_{\mathcal{S}},\mathbf{d}_{\mathcal{S}}\Bigr\rangle\Bigr|\,dt \\
				&\le\frac12\bigl\langle \nabla_{\mathcal{S}}^{2}F(\x)\mathbf{d}_{\mathcal{S}},\mathbf{d}_{\mathcal{S}}\rangle\\
				&\quad+\int_{0}^{1}\bigl\|\nabla_{\mathcal{S}}F(\x+t\mathbf{d})
				-\nabla_{\mathcal{S}}F(\x)
				-t\nabla_{\mathcal{S}}^{2}F(\x)\mathbf{d}_{\mathcal{S}}\bigr\|\,\|\mathbf{d}_{\mathcal{S}}\|\,dt \\
				&\le \frac12\bigl\langle \nabla_{\mathcal{S}}^{2}F(\x)\mathbf{d}_{\mathcal{S}},\mathbf{d}_{\mathcal{S}}\rangle+\int_{0}^{1}\Bigl(a\,t\|\mathbf{d}\|
				+\frac{M}{2}t^{2}\|\mathbf{d}\|^{2}\Bigr)\|\mathbf{d}\|\,dt\\
				&=\frac12\bigl\langle \nabla_{\mathcal{S}}^{2}F(\x)\mathbf{d}_{\mathcal{S}},\mathbf{d}_{\mathcal{S}}\rangle+\frac{a}{2}\|\mathbf{d}\|^{2}+\frac{M}{6}\|\mathbf{d}\|^{3}.
			\end{aligned}
		\end{equation*}
		
		\noindent Combining the above estimates establishes \eqref{ineq:cubic_upper} and completes the
		proof.
	\end{proof}

The subspace line search for \(\beta_k\) terminates after finitely many trials and gives a descent estimate.

%The following Theorem shows that the line search for $\beta_{k}$ in Algorithm \ref{alg:apgcn} will terminate in finite step.
%The following theorem establishes the finite termination of the line search for \(\beta_k\) and the descent property.
\begin{theorem}
\label{Newton_decrease}
Assume that $F$ defined in \eqref{Problem} satisfies Assumption~\ref{Ass1}.
{Assume also that $M(\mathbf x^k)>0$ and
\[
M(\mathbf x^k)\geq\underline\omega\,
\mathrm{lip}_{\mathcal M_{\mathbf x^k}}\nabla^2_{\mathcal S^k}F(\mathbf x^k),
\qquad \underline\omega>1.
\]}
At iteration $k$, suppose that $\mathcal{S}^k\neq\emptyset$ and $c>0$. Let $\mathbf{d}^{k}$ be an inexact solution of the cubic subproblem satisfying Conditions~(\ref{inexact1}) and (\ref{inexact2}).
Assume that
\begin{equation}
	\label{sigmak}
	\sigma_k\geqslant \frac{2(M(\mathbf{x}^{k})+c)}{3}.
\end{equation}
Then, the procedure of finding the smallest $j_k\in\mathbb N$ such that
\begin{equation}
	\label{back}
	F(\mathbf{y}^{k,j_k})
	\leq
	F(\mathbf{x}^{k})
	-\frac{\eta}{2}\|\mathbf{y}^{k,j_k}-\mathbf{x}^{k}\|^{2}
	-\frac{c}{6}\|\mathbf{y}^{k,j_k}-\mathbf{x}^{k}\|^{3}
\end{equation}
terminates after finitely many trials. %Since \eqref{back} implies the weaker acceptance condition with quadratic coefficient $\eta/4$ used in Algorithm~\ref{alg:apgcn}, its backtracking procedure also terminates.
\end{theorem}

\begin{proof}
Let $\widehat{\mathbf d}^k\in\mathbb R^n$ be the zero extension of
$\mathbf d^k$, so $\|\widehat{\mathbf d}^k\|=\|\mathbf d^k\|$.
 Since $M(\mathbf x^k)>\mathrm{lip}_{\mathcal M_{\mathbf x^k}}\nabla^2_{\mathcal S^k}F(\mathbf x^k)$, strict continuity supplies $\delta_k>0$ such that
\[
\|\nabla^2_{\mathcal S^k}F(\mathbf z)-\nabla^2_{\mathcal S^k}F(\mathbf x^k)\|
\leq M(\mathbf x^k)\|\mathbf z-\mathbf x^k\|
\quad\text{for }\mathbf z\in\mathcal B(\mathbf x^k,\delta_k)\cap\mathcal M_{\mathbf x^k}.
\]
Because $\mathcal S^k\subseteq\operatorname{supp}(\mathbf x^k)$, there is $\widehat\beta_k\in(0,1]$ such that $\mathbf x^k+t\beta\widehat{\mathbf d}^k\in\mathcal B(\mathbf x^k,\delta_k)\cap\mathcal M_{\mathbf x^k}$ for all $t\in[0,1]$ and $0<\beta\leq\widehat\beta_k$. Lemma~\ref{lem:inexact_cubic_upper}, with $a=0$, therefore gives
\begin{equation}
	\label{eq:cubic_upper}
	F(\mathbf{x}^k+\beta\widehat{\mathbf d}^k)
	\leq
	F(\mathbf{x}^{k})
	+\beta\langle \mathbf{g}^{k},\mathbf{d}^{k}\rangle
	+\frac{\beta^{2}}{2}\langle H^{k}\mathbf{d}^{k},\mathbf{d}^{k}\rangle
	+\frac{M(\mathbf{x}^{k})\beta^{3}}{6}\|\mathbf{d}^{k}\|^{3}.
\end{equation}
Recall that $\mathbf{r}^{k}=\mathbf{g}^{k}+(H^{k}+\eta I)\mathbf{d}^{k}+\frac{\sigma_k}{2}\|\mathbf{d}^{k}\|\mathbf{d}^{k}$.
By Condition~(\ref{inexact1}) we have $\langle\mathbf{r}^{k},\mathbf{d}^{k}\rangle\leq0$, and hence
\begin{equation*}
	\langle\mathbf{g}^{k},\mathbf{d}^{k}\rangle
	\leq
	-\langle(H^{k}+\eta I)\mathbf{d}^{k},\mathbf{d}^{k}\rangle
	-\frac{\sigma_k}{2}\|\mathbf{d}^{k}\|^{3}.
\end{equation*}
Moreover, Condition~(\ref{inexact2}) gives
\begin{equation*}
	\langle H^{k}\mathbf{d}^{k},\mathbf{d}^{k}\rangle
	\geq
	-\eta\|\mathbf{d}^{k}\|^{2}
	-\frac{\sigma_k}{2}\|\mathbf{d}^{k}\|^{3}.
\end{equation*}
Substituting the above two inequalities into \eqref{eq:cubic_upper}, we obtain
\begin{align*}
	F(\mathbf{x}^k+\beta\widehat{\mathbf d}^k) &\leq
	F(\mathbf{x}^{k})
	-\eta\beta\|\mathbf{d}^{k}\|^2
	-\frac{\sigma_k\beta}{2}\|\mathbf{d}^{k}\|^3\\
	&\quad+\left(-\beta+\frac{\beta^2}{2}\right)
	\langle H^k\mathbf{d}^{k},\mathbf{d}^{k}\rangle
	+\frac{M(\mathbf{x}^{k})\beta^3}{6}\|\mathbf{d}^{k}\|^3 \\
	&\leq
	F(\mathbf{x}^{k})
	-\frac{\eta\beta^2}{2}\|\mathbf{d}^{k}\|^2
	-\left(\frac{\sigma_k\beta^2}{4}
	-\frac{M(\mathbf{x}^{k})\beta^3}{6}\right)
	\|\mathbf{d}^{k}\|^3.
\end{align*}
By \eqref{sigmak} and $\beta\in(0,1]$, we have
\begin{equation*}
	\frac{\sigma_k\beta^2}{4}
	-\frac{M(\mathbf{x}^{k})\beta^3}{6}
	\geq
	\frac{\beta^3}{6}
	\left(\frac{3\sigma_k}{2}-M(\mathbf{x}^{k})\right)
	\geq
	\frac{c\beta^3}{6}.
\end{equation*}
Therefore, for every $\beta\in(0,\widehat{\beta}_k]$, it holds that
\begin{equation*}
	F(\mathbf{x}^{k}+\beta\widehat{\mathbf d}^k)
	\leq
	F(\mathbf{x}^{k})
	-\frac{\eta\beta^2}{2}\|\mathbf{d}^{k}\|^2
	-\frac{c\beta^3}{6}\|\mathbf{d}^{k}\|^3.
\end{equation*}
Since $\mathbf{y}^{k,j_k}=\mathbf{x}^{k}+\beta_{k,j_k}\widehat{\mathbf d}^k$ and $\beta_{k,j_{k}+1}=\beta_{k,j_k}q$, there exists a finite $\widehat{j}_k$ such that \eqref{back} holds for $j_k=\widehat{j}_k$.
The test in Algorithm~\ref{alg:apgcn} uses the smaller quadratic coefficient $\eta/4$, so it also terminates after finitely many trials.
\end{proof}

The following theorem establishes sufficient decrease at each outer iteration.
\begin{theorem}
\label{lem:descent}
Suppose  Assumption~\ref{Ass1} holds.  Let the sequences
$\{\mathbf{x}^{k}\}$ and $\{\mathbf{y}^k\}$ be generated by Algorithm~\ref{alg:apgcn}.
Let $\mathbf x^0=\mathbf y^0$.
If $\sigma_{k}$ satisfies \eqref{sigmak}, then there exists a
constant $\widetilde C > 0$ such that the following
property holds for $k\ge1$:
\begin{equation}\label{ineq:x_sufficient_descent}
	F(\mathbf{x}^{k})
	\;\leq\;
	F(\mathbf{x}^{k-1})
	- \widetilde C\|\mathbf{x}^{k}-\mathbf{y}^{k-1}\|^2
	- \widetilde C\|\mathbf{x}^{k-1}-\mathbf{y}^{k-1}\|^2.
\end{equation}
Consequently,
\begin{equation}\label{adjent}
	\lim_{k\to\infty}\|\mathbf{x}^{k}-\mathbf{x}^{k-1}\|
	=
	\lim_{k\to\infty}\|\mathbf{x}^{k}-\mathbf{y}^{k-1}\|
	=
	\lim_{k\to\infty}\|\mathbf{x}^{k-1}-\mathbf{y}^{k-1}\|
	=0.
\end{equation}
Furthermore, there exists a constant $\widehat{F}$ such that
$
\lim_{k\to\infty} F(\mathbf{x}^{k}) = \widehat{F}.
$
\end{theorem}

\begin{proof}
From Theorem \ref{thm:PG_linesearch}, we have
\begin{equation}\label{ieq:PG_descent}
	F(\mathbf{x}^{k}) \le F(\mathbf{y}^{k-1})
	- c_1\|\mathbf{x}^{k}-\mathbf{y}^{k-1}\|^2 .
\end{equation}
For $k=1$, inequality~\eqref{ieq:PG_descent} and $\mathbf x^0=\mathbf y^0$ directly give~\eqref{ineq:x_sufficient_descent}. For $k\geq2$, we consider two cases: (a) $\mathcal{S}^{k-1} \neq \emptyset$ and
(b) $\mathcal{S}^{k-1} = \emptyset$. \\
In case (a), the cubic line search gives $
F(\mathbf{y}^{k-1})
\le
F(\mathbf{x}^{k-1})
-\frac{\eta}{4}\|\mathbf{y}^{k-1}-\mathbf{x}^{k-1}\|^{2}
-\frac{c}{6}\|\mathbf{y}^{k-1}-\mathbf{x}^{k-1}\|^{3}.
$
Combining this inequality with  (\ref{ieq:PG_descent}), we obtain
\begin{align*}
	&F(\mathbf{x}^{k}) \le F(\mathbf{x}^{k-1})
	- c_1\|\mathbf{x}^{k}-\mathbf{y}^{k-1}\|^2
	- \frac{c}{6}\|\mathbf{x}^{k-1}-\mathbf{y}^{k-1}\|^3
	- \frac{\eta}{4}\|\mathbf{x}^{k-1}-\mathbf{y}^{k-1}\|^2 \notag \\
	&\le F(\mathbf{x}^{k-1})
	- c_1\|\mathbf{x}^{k}-\mathbf{y}^{k-1}\|^2
	- \frac{\eta}{4}\|\mathbf{x}^{k-1}-\mathbf{y}^{k-1}\|^2 \notag
\end{align*}
\noindent Let $\widetilde C := \min\left\{c_1, \frac{\eta}{4}\right\}$. Then
(\ref{ineq:x_sufficient_descent}) follows.
Case (b). It leads to $\y^{k-1} = \x^{k-1}$.  (\ref{ineq:x_sufficient_descent}) also holds.
Rearranging (\ref{ineq:x_sufficient_descent}) and summing the resulting inequalities yields
$$
\widetilde C \sum_{k=1}^{\infty}
\left(
\|\mathbf{x}^{k}-\mathbf{y}^{k-1}\|^2
+ \|\mathbf{x}^{k-1}-\mathbf{y}^{k-1}\|^2
\right)
\le F(\mathbf{x}^0) - \underline{F} < \infty .
$$
Hence, \eqref{adjent} holds.
Since $\{F(\mathbf{x}^{k})\}$ is nonincreasing and bounded below, there exists a constant $\widehat{F}$ such that
$\lim_{k\to\infty} F(\mathbf{x}^{k}) = \widehat{F}.$
\end{proof}

\subsection{Global Convergence}
%We assume boundedness to prove convergence of the entire sequence.
Boundedness gives accumulation points and is a hypothesis of the whole-sequence convergence theorem.

\begin{assumption}
\label{Ass2}
The sequence $\{\mathbf{x}^k\}$ generated from Algorithm~\ref{alg:apgcn} is bounded.
\end{assumption}
Coercivity of $F$ implies Assumption~\ref{Ass2}. Appendix~\ref{sec:boundness_sequence} gives boundedness conditions for the ratio-regularized least-squares and regularized logistic models.

\begin{lemma}\label{lem:critical}
Suppose that Assumptions~\ref{Ass1} and~\ref{Ass2} hold. Let $\{\mathbf{x}^{k}\}$ and $\{\mathbf{y}^{k}\}$ be generated from Algorithm~\ref{alg:apgcn}. Let $\h x^0=\h y^0$.
Let $\overline{\mathcal X}$ denote the set of accumulation points
of $\{\mathbf x^k\}$.
If $\sigma_{k}$ satisfies \eqref{sigmak},
then the following assertions hold:
\begin{itemize}
	\item[(i)] The sequence $\{\mathbf{y}^k\}$ is bounded and admits the same set of accumulation points $\overline{\mathcal{X}}$.
	\item[(ii)] $F(\overline{\mathbf{x}})= \widehat{F}$ for all $\overline{\mathbf{x}}\in\overline{\mathcal{X}}$,
	where $\widehat{F}$ is defined in Theorem~\ref{lem:descent}.
	\item[(iii)] Every accumulation point $\overline{\mathbf{x}}$ of $\{\mathbf{x}^k\}$ and $\{\mathbf{y}^k\}$ is a
	critical point of problem~\eqref{Problem}, i.e.,
	\begin{equation}\label{critical}
		{{\bf 0}} \in\partial g(\overline{\mathbf{x}}) + \nabla h(\overline{\mathbf{x}}).\end{equation}
\end{itemize}
\end{lemma}
\begin{proof}

By \eqref{adjent}, $\|\mathbf{x}^k-\mathbf{y}^k\|\to0$,
which implies that $\{\mathbf{y}^k\}$ is bounded and shares
the same accumulation points as $\{\mathbf{x}^k\}$.
This establishes (i).
To prove (ii), let $\mathbf{x}^{k_j}\to\overline{\mathbf{x}}$
be any convergent subsequence.
By \eqref{adjent}, we also have
$\mathbf{y}^{k_j-1}\to\overline{\mathbf{x}}$.
The lower semicontinuity of $F$ then yields
$F(\overline{\mathbf{x}})\leq\widehat F<+\infty$.
 Setting $\mathbf{z}^{j}:=\mathbf{y}^{k_j-1}-\alpha_{k_j}\nabla h(\mathbf{y}^{k_j-1})$, we have
$g(\mathbf{x}^{k_j})+\frac{1}{2\alpha_{k_j}}
	\|\mathbf{x}^{k_j}-\mathbf{z}^{j}\|^2
	\leq g(\overline{\mathbf{x}})+\frac{1}{2\alpha_{k_j}}
	\|\overline{\mathbf{x}}-\mathbf{z}^{j}\|^2.
$
Since $\{\mathbf z^j\}$ is bounded,
\[
\left|\|\mathbf x^{k_j}-\mathbf z^j\|^2
-\|\overline{\mathbf x}-\mathbf z^j\|^2\right|
\le
\|\mathbf x^{k_j}-\overline{\mathbf x}\|
\bigl(\|\mathbf x^{k_j}-\mathbf z^j\|
+\|\overline{\mathbf x}-\mathbf z^j\|\bigr)
\to0.
\]
{Since $\alpha_{k_j}\ge\underline\alpha>0$, we obtain
$\limsup_{j\to\infty}g(\mathbf x^{k_j})
\le g(\overline{\mathbf x})$.
By the continuity of $h$, it follows that
$\widehat F\le F(\overline{\mathbf x})$.}
Consequently, $F(\overline{\mathbf{x}})=\widehat F$,
which completes the proof of (ii).
Finally, the proximal optimality condition implies

\begin{equation}\label{Temp}
\begin{aligned}
\mathbf{v}^k&:=\frac{\mathbf{y}^{k-1}-\mathbf{x}^k}{\alpha_k}
+\nabla h(\mathbf{x}^k)-\nabla h(\mathbf{y}^{k-1})
\in\partial F(\mathbf{x}^k),\\
\|\mathbf{v}^k\|&\leq(\alpha_k^{-1}+L)
\|\mathbf{x}^k-\mathbf{y}^{k-1}\|.
\end{aligned}
\end{equation}
Along the above subsequence, $\mathbf{v}^{k_j}\to\mathbf{0}$ and
$F(\mathbf{x}^{k_j})\to F(\overline{\mathbf{x}})$.
The defining closedness property of the limiting subdifferential gives
$\mathbf{0}\in\partial F(\overline{\mathbf{x}})
=\partial g(\overline{\mathbf{x}})+\nabla h(\overline{\mathbf{x}})$,
proving (iii).
\end{proof}

The residual estimate below follows from \eqref{Temp}.
\begin{lemma}
\label{thm:subgradient_bound}
Suppose that Assumptions~\ref{Ass1} and~\ref{Ass2} hold.
Let $\{\mathbf{x}^{k}\}$ and $\{\mathbf{y}^{k}\}$ be generated from Algorithm~\ref{alg:apgcn}.
If $\sigma_{k}$ satisfies \eqref{sigmak}, then the following holds:
\begin{equation}\label{ineq:dist_upper}
\begin{aligned}
	\mathrm{dist}\bigl({{\bf 0}},\partial F(\mathbf{x}^{k})\bigr)
	&\le \left(\frac{1}{\alpha_k}+L\right)
	\|\mathbf{x}^k - \mathbf{y}^{k-1}\|\\
	&\le \left(\frac{1}{\underline{\alpha}}+L\right)
	\bigl(\|\mathbf{x}^{k} - \mathbf{y}^{k-1}\|
	+ \|\mathbf{x}^{k-1} - \mathbf{y}^{k-1}\|\bigr).
\end{aligned}
\end{equation}

\end{lemma}
%\begin{proof}
%The proximal optimality condition gives
%\[
%\frac{\mathbf{y}^{k-1}-\mathbf{x}^k}{\alpha_k}
%+\nabla h(\mathbf{x}^k)-\nabla h(\mathbf{y}^{k-1})
%\in\partial F(\mathbf{x}^k).
%\]
%Taking norms and using the $L$-Lipschitz continuity of $\nabla h$
%together with $\alpha_k\ge\underline\alpha$ yields
%\eqref{ineq:dist_upper}.
%\end{proof}
To prove convergence of the entire sequence generated by Algorithm~\ref{alg:apgcn}, we use the KL inequality to show that the sequence has finite length. The proof follows \cite[Theorem~4]{LiPong15}, \cite[Theorem~3.4]{BCN19}, and \cite[Theorem~6.11]{boct2025full}.

%\begin{proof}
%The proximal optimality condition gives
%\[
%\frac{\mathbf{y}^{k-1}-\mathbf{x}^k}{\alpha_k}
%+\nabla h(\mathbf{x}^k)-\nabla h(\mathbf{y}^{k-1})
%\in\partial F(\mathbf{x}^k).
%\]
%Taking norms and using the $L$-Lipschitz continuity of $\nabla h$
%together with $\alpha_k\ge\underline\alpha$ yields
%\eqref{ineq:dist_upper}.
%\end{proof}
%Next, we establish the global convergence of Algorithm~\ref{alg:apgcn}, following the arguments in \cite[Theorem~4]{LiPong15}, \cite[Theorem~3.4]{BCN19}, and \cite[Theorem~6.11]{boct2025full}.

\begin{theorem}
\label{thm:global_convegence}
Suppose that Assumptions~\ref{Ass1} and~\ref{Ass2} hold, and let $\{\mathbf{x}^{k}\}$ and $\{\mathbf{y}^{k}\}$ be generated from Algorithm~\ref{alg:apgcn}. Let $\h x^0=\h y^0$. Assume that $F$ satisfies the Kurdyka--\L ojasiewicz (KL) property at every point in $\dom(\partial F)$. If $\sigma_k$ satisfies \eqref{sigmak}, then
$$
\sum_{k=1}^{\infty} \|\mathbf{x}^{k} - \mathbf{x}^{k-1}\| \leq\sum_{k=1}^{\infty} \left(||{\mathbf{x}}^{k} - {\mathbf{y}}^{k-1} || + ||{\mathbf{x}}^{k-1} - {\mathbf{y}}^{k-1} ||\right) < \infty.$$
Consequently, the sequences  $\{\mathbf{x}^{k}\}$ and $\{\mathbf{y}^{k}\}$  both converge to the same stationary point $\overline{\mathbf{x}}$ satisfying \eqref{critical}.
\end{theorem}

\begin{proof}
Define $\mu_k:=F(\mathbf{x}^k)-\widehat F\geq0$.
%and let
%$\Omega:=\overline{\mathcal X}$ be the accumulation-point set.
Lemma \ref{lem:critical} shows that $\overline{\mathcal X}$ is nonempty and compact with $F\equiv\widehat F$ on $\overline{\mathcal X}$ and $\operatorname{dist}(\mathbf{x}^k,\overline{\mathcal X})\to0$.
A finite collection of KL neighborhoods covers the compact set
$\overline{\mathcal X}$. Restrict their desingularizing functions to a
common interval and let $\phi$ be their sum. The functions have positive
derivatives, so $\phi'$ dominates each individual derivative. For all
sufficiently large $k$ with $\mu_{k-1}>0$, the iterate $\mathbf x^{k-1}$
lies in at least one of these neighborhoods and $\mu_{k-1}$ lies in the
common interval. Hence
\begin{equation}\label{ineq:KL_used}
	\phi'(\mu_{k-1})\,
	\operatorname{dist}(\mathbf{0},\partial F(\mathbf{x}^{k-1}))\geq1.
\end{equation}
If $\mu_K=0$ for some $K$, \eqref{ineq:x_sufficient_descent} and
monotonicity force
$\mathbf{x}^k=\mathbf{y}^k=\mathbf{x}^K$ for every $k\geq K$,
and the theorem follows immediately. Otherwise every $\mu_k$ is positive.
In this case \eqref{ineq:dist_upper} and \eqref{ineq:KL_used} ensure
that the denominators below are nonzero for all sufficiently large $k$.
For such $k$, the concavity of $\phi$ gives
\begin{equation*}
	\begin{aligned}
		&\phi(\mu_{k-1})-\phi(\mu_k) \ge \phi'(\mu_{k-1}) (\mu_{k-1}-\mu_k )\\
		& \ge \widetilde{C}\phi'(\mu_{k-1})\left(\|\mathbf{x}^k - \mathbf{y}^{k-1}\|^2 +\|\mathbf{x}^{k-1} - \mathbf{y}^{k-1}\|^2\right) \\
		& \ge \frac{\|\mathbf{x}^k - \mathbf{y}^{k-1}\|^2 +\|\mathbf{x}^{k-1} - \mathbf{y}^{k-1}\|^2 }{\widehat M \left(\|\mathbf{x}^{k-1} - \mathbf{y}^{k-2}\|+ \|\mathbf{x}^{k-2} - \mathbf{y}^{k-2}\|\right)},
	\end{aligned}
\end{equation*}
where $\widehat M = (1/\underline{\alpha} + L)/\widetilde{C}$. Rearranging the above inequality yields
\begin{equation*}
	\begin{aligned}
		&	\|\mathbf{x}^k - \mathbf{y}^{k-1}\| +\|\mathbf{x}^{k-1} - \mathbf{y}^{k-1}\|\\
		& \leq \sqrt{2\widehat M(\phi(\mu_{k-1}) - \phi(\mu_{k}))\left(\|\mathbf{x}^{k-1} - \mathbf{y}^{k-2}\|+ \|\mathbf{x}^{k-2} - \mathbf{y}^{k-2}\|\right)}\\
		& \leq {\widehat M}(\phi(\mu_{k-1}) - \phi(\mu_{k})) + \frac{1}{2}\left(\|\mathbf{x}^{k-1} - \mathbf{y}^{k-2}\|+ \|\mathbf{x}^{k-2} - \mathbf{y}^{k-2}\|\right).
	\end{aligned}
\end{equation*}
So,
\begin{equation}\label{ineq:KL_telescope}
	\begin{aligned}
		2\bigl(\|\mathbf{x}^k-\mathbf{y}^{k-1}\|
		+\|\mathbf{x}^{k-1}-\mathbf{y}^{k-1}\|\bigr)
		&\leq 2\widehat M\bigl(\phi(\mu_{k-1})-\phi(\mu_k)\bigr)\\
		&\quad+\|\mathbf{x}^{k-1}-\mathbf{y}^{k-2}\|
		+\|\mathbf{x}^{k-2}-\mathbf{y}^{k-2}\|.
	\end{aligned}
\end{equation}
Summing \eqref{ineq:KL_telescope} and using $\phi\geq0$ bounds the partial sums uniformly. Hence
%\begin{align*}
%		\sum_{k = K_0 + 2}^{K}\|\mathbf{x}^k - \mathbf{x}^{k-1}\| +\|\mathbf{x}^{k-1} - \mathbf{y}^{k-1}\| \leq 4M\phi(r_{K_0+1}) + \left(\|\mathbf{x}^{K_0+1} - \mathbf{x}^{K_0}\| + \|\mathbf{x}^{K_0} - \mathbf{y}^{K_0}\| \right).
%	\end{align*}
%	Let $K \to \infty$, we obtain
\[
\sum_{k=1}^{\infty} \left(||{\mathbf{x}}^{k} - {\mathbf{y}}^{k-1} || + ||{\mathbf{x}}^{k-1} - {\mathbf{y}}^{k-1} ||\right) < \infty.
\]
Thus, $\{\mathbf{x}^k\}$  converges to some $\overline{\mathbf{x}}$ satisfying (\ref{critical}). Moreover, $\{\y^k\}$ also converges to $\overline{\mathbf{x}}$.
\end{proof}

%	{
%	Next, we review the standard error bound results under different KL exponents follows from the proof of \citep[Theorem 2]{attouch2009convergence}.
\begin{theorem}\label{thm:KL_rate}
Suppose that the assumptions in Theorem~\ref{thm:global_convegence} hold.
Let $\{\x^k\}$ and $\{\mathbf{y}^{k}\}$ be the sequences generated from  Algorithm~\ref{alg:apgcn}. Let $\h x^0=\h y^0$.
Assume that $F$ satisfies the KL property at the limit point
$\overline{\mathbf{x}}$ with exponent $\theta\in(0,1)$ and
desingularizing function $\phi(s)=s^{1-\theta}/[\nu(1-\theta)]$, where $\nu>0$ is the constant in Definition~\ref{def:KL_exponent}.
Then the following statements hold:
\begin{enumerate}
	\item[(i)] For $\theta\in(0,1)$, there exists a constant $\widehat{C} > 0$ such that, for all sufficiently large $k$,
	\begin{equation*}
		\|\x^k - \overline{\x}\|  \le \widehat C\left(\|\mathbf{x}^{k} - \mathbf{y}^{k-1}\| + \|\mathbf{x}^{k-1} - \mathbf{y}^{k-1}\| \right)^{\displaystyle{\min\left(1,\frac{1-\theta}{\theta}\right)}}.
	\end{equation*}
	\item[(ii)] If $\theta \in (0,\tfrac{1}{2}]$,  there exist constants $\widehat{C}_1 > 0$ and $\rho \in (0,1)$ such that, for all sufficiently large $k$,
	\begin{equation*}
		\|\x^k - \overline{\x}\| \le \widehat{C}_1 \rho^k.
	\end{equation*}
	\item[(iii)] If $\theta \in (\tfrac{1}{2},1)$,  there exists a constant $\widehat{C}_2 > 0$ such that, for all sufficiently large $k$,
	\[
	\|\x^k - \overline{\x}\| \le \widehat{C}_2 \, k^{-\frac{1-\theta}{2\theta-1}}.
	\]
\end{enumerate}

\end{theorem}
%}
\begin{proof}
If the iterates are eventually constant, all conclusions are immediate.
Otherwise, restrict to a sufficiently large index range on which the KL
inequality at $\overline{\mathbf{x}}$ holds and
\[
s_k:=\|\mathbf{x}^k-\mathbf{y}^{k-1}\|
+\|\mathbf{x}^{k-1}-\mathbf{y}^{k-1}\|\leq1.
\]
(i) Apply the argument for \eqref{ineq:KL_telescope} on this tail with this $\phi$. Combining \eqref{ineq:dist_upper} and \eqref{ineq:KL_used}, with
$k$ replaced by $k+1$, and using
$\phi'(\mu_k)=\nu^{-1}\mu_k^{-\theta}$ gives
\begin{equation}
\left(\frac{1}{\underline{\alpha}}+L\right)
\bigl(\|\mathbf{x}^{k} - \mathbf{y}^{k-1}\|
+ \|\mathbf{x}^{k-1} - \mathbf{y}^{k-1}\|\bigr)\ge  \mathrm{dist}(0,\partial F(\mathbf{x}^{k})) \ge \nu\mu_k^\theta.	\label{keyineq}
\end{equation}

\noindent Define $\Delta_k:=\sum_{p=k+1}^{\infty}s_p$. Summing \eqref{ineq:KL_telescope} from $p=k+1$ onward yields
\begin{equation*}
\begin{aligned}
	&2\Delta_k = 2\sum_{p=k+1}^{+\infty} \left(\|\mathbf{x}^p-{\mathbf{y}}^{p-1}\|+\|\mathbf{x}^{p-1}-\mathbf{y}^{p-1}\|\right)\nn\\
	&\le 2\widehat M\phi(\mu_k)+
	%\underbrace{\sum_{p=k+1}^{+\infty}\left(
		%\|\mathbf{x}^p-{\mathbf{y}}^{p-1}\|+\|\mathbf{x}^{p-1}-\mathbf{y}^{p-1}\|\right)}_{\Delta_k}
	\Delta_k+\|\mathbf{x}^k-\mathbf{y}^{k-1}\|+\|\mathbf{x}^{k-1}-\mathbf{y}^{k-1}\|\end{aligned}
	\end{equation*}
	Combining~\eqref{keyineq} with the preceding bound yields
	\begin{eqnarray}\label{keyin}	
&&\Delta_k\le 2\widehat M\phi(\mu_k)+\|\mathbf{x}^k-\mathbf{y}^{k-1}\|+\|\mathbf{x}^{k-1}-\mathbf{y}^{k-1}\|\nn\\
&&\le \frac{2\widehat M}{\nu(1-\theta)}\biggl(\frac{1}{\nu}(\frac{1}{\underline{\alpha}}+L)[\|\mathbf{x}^{k} - \mathbf{y}^{k-1}\|
+ \|\mathbf{x}^{k-1} - \mathbf{y}^{k-1}\|]\biggr)^{\frac{1-\theta}{\theta}}\nn\\
&&+[\|\mathbf{x}^{k} - \mathbf{y}^{k-1}\|
+ \|\mathbf{x}^{k-1} - \mathbf{y}^{k-1}\|].
\end{eqnarray}
Let
$
\widehat C = \frac{2\widehat M}{\nu(1-\theta)}\left(\nu^{-1}\Big(\tfrac{1}{\underline{\alpha}} + L\Big)\right)^{\frac{1-\theta}{\theta}} + 1.
$
The triangle inequality gives
$
\|\mathbf{x}^k-\overline{\mathbf{x}}\|
\leq\sum_{p=k}^{\infty}\|\mathbf{x}^{p+1}-\mathbf{x}^p\|
\leq\Delta_k.
$
Since $s_k\leq1$, \eqref{keyin} implies
$\Delta_k\leq\widehat C s_k^{\min\{1,(1-\theta)/\theta\}}$,
which proves (i).
\\
%If $\theta \in (0,\tfrac{1}{2}]$, then $\frac{1-\theta}{\theta} \ge 1$.
%Hence, the desired result follows directly by noting that
%$\|\x^k - \overline{\x}\| < \Delta_k$.
%If $\theta \in (\tfrac{1}{2},1)$, then $\frac{1-\theta}{\theta} \le 1$,
%and the same conclusion still holds.
(ii)
For $\theta \in (0,\tfrac{1}{2}]$, it follows from (\ref{keyin}) that
$
\Delta_k \le \widehat C \, (\Delta_{k-1} - \Delta_k).
$
The remainder of the proof proceeds along the same lines as
\citep[Theorem 2(ii)]{attouch2009convergence}.\\
\noindent
(iii)
From~\eqref{keyin}, we obtain
$
\Delta_k^{\frac{\theta}{1-\theta}}
\le \widehat C^{\frac{\theta}{1-\theta}} (\Delta_{k-1} - \Delta_k).
$
The desired result then follows by applying the argument of
\citep[Theorem 2(iii)]{attouch2009convergence}.
\end{proof}

\section{Finite Identification via Subspace Selection}
\label{sec:finite}
We specify thresholds $\varepsilon_k$ for selecting the index set
$\mathcal S^k$ to ensure finite identification. Here, finite
identification means that there exists an integer $K$ such that
$\mathcal S^k=\overline{\mathcal S}$ for all $k\geq K$, where
$\overline{\mathcal S}$ denotes the support of the limit point
of the sequence $\{\h{x}^k\}$. %\cite{hare2004identifying}.

\begin{lemma}
\label{lem:identify_general}
Suppose that the thresholds $\varepsilon_k$ in
Algorithm~\ref{alg:apgcn} are positive and satisfy
$\lim_{k\to\infty}\varepsilon_k=0$.
Let $\{\x^k\}$ and $\{\mathbf{y}^{k}\}$ be generated from  Algorithm~\ref{alg:apgcn}. Let $\h x^0=\h y^0$.
Suppose the hypotheses of Theorem~\ref{thm:global_convegence} hold. Let $\overline{\x}$ be the limit point of $\{\x^k\}$, and define $\overline{\mathcal{S}} = \supp(\overline{\x})$ and $\overline{n} = \|\overline{\mathbf{x}}\|_0$. Then,
\begin{itemize}
\item[(i)] there exists an index $K_0$ such that
$\overline{\mathcal{S}} \subseteq \mathcal{S}^{k},~ \forall k\ge K_{0};$
\item[(ii)] for all $k\ge K_0$,
$\overline n\leq n_k\leq
\min\{\lfloor\|\mathbf{x}^k-\overline{\mathbf{x}}\|^2/
\varepsilon_k^2\rfloor+\overline n,\|\mathbf{x}^k\|_0\}$.
\end{itemize}

\end{lemma}

\begin{proof}
%(i) The proof is similar to \citep[Theorem 4]{li2015global} and \citep[Theorem 6.11]{boct2025full}.\\
(i) If $\overline{\mathcal{S}} = \emptyset$, take $K_0=1$. Suppose $\overline{\mathcal{S}} \neq \emptyset$ and denote $\rho_{\min} := \min\{|\overline{x}_i|:i\in \overline{\mathcal{S}}\} > 0$.
Choose a common index $K_0$ such that $\varepsilon_k\le \rho_{\min}/2$ and $
\|\mathbf{x}^k - \overline{\mathbf{x}}\| \le \rho_{\min}/2$ when $k\ge K_0$.
Thus,
\begin{equation*}
|x_i^k| \ge\;|\overline{x}_i| - |x_i^k - \overline{x}_i|\ge \rho_{\min} - \|\mathbf{x}^k - \overline{\mathbf{x}}\| \ge\frac{\rho_{\min}}{2}\ge \varepsilon_k,
\end{equation*}
which implies $i\in\mathcal{S}^k$.
(ii) Let $k\geq K_0$. By assumption, $\varepsilon_k>0$.
Let  $T_{k}:=\mathcal{S}^{k}\setminus\overline{\mathcal S}$.
Because $\overline{x}_i=0$ and $|x_i^k|\ge\varepsilon_k$ for $i\in T_k$,
$\|{\mathbf{x}}^k-\overline {\mathbf{x}}\|^2\ge \sum_{i\in T_{k}} |x_i^k|^2 \ge \sharp(T_{k})\varepsilon^2_k=(n_k-\overline n)\varepsilon_k^2. $
The assertion follows directly.
\end{proof}

%We first consider the case where $F$ satisfies KL property with exponent $\theta\in(0,1)$, see Definition \ref{def:KL_exponent}.

\begin{theorem}
\label{thm:Finite_identify_KL_general}
Suppose that all the assumptions of
Theorem~\ref{thm:global_convegence} hold, and let
$\{\mathbf{x}^k\}$ and $\{\mathbf{y}^k\}$ be generated by
Algorithm~\ref{alg:apgcn}, with $\mathbf{x}^0=\mathbf{y}^0$.
Let $\overline{\mathbf{x}}$ be its limit and
$\overline{\mathcal S}:=\operatorname{supp}(\overline{\mathbf{x}})$.
Assume that $F$ has the KL property at $\overline{\mathbf{x}}$
with exponent $\theta\in(0,1)$. Define

\begin{eqnarray}
	&&s_k:=\|\mathbf{x}^{k}-\mathbf{y}^{k-1}\|
	+\|\mathbf{x}^{k-1}-\mathbf{y}^{k-1}\|,\nonumber\\
	&&a_k:=\max\{s_k,{C_a}/({k+1})\}>0,\;
	\nu:=\min\{1,(1-\theta)/\theta\},\label{nu}
\end{eqnarray}
where $C_a >0$.
For any $C_{\varepsilon}>0$ and $0<\gamma<\nu$, take
$
\varepsilon_k:=C_{\varepsilon} a_k^\gamma.
$
Then $\mathcal S^k=\overline{\mathcal S}$ for all sufficiently large $k$.
\end{theorem}
\begin{proof}
We have $a_k>0$ and $a_k\to0$ by \eqref{adjent}; hence
$\varepsilon_k>0$ and $\varepsilon_k\to0$.
Lemma~\ref{lem:identify_general} gives
$\overline{\mathcal S}\subseteq\mathcal S^k$ eventually.
By Theorem~\ref{thm:KL_rate}(i), for all sufficiently large $k$,
$
\|\mathbf{x}^k-\overline{\mathbf{x}}\|
\leq\widehat C s_k^\nu\leq\widehat C a_k^\nu.
$
Consequently, the  bound in Lemma~\ref{lem:identify_general} gives
$
0\leq n_k-\overline n\leq
\left\lfloor\frac{\widehat C^2}{C_{\varepsilon}^2}
a_k^{2(\nu-\gamma)}\right\rfloor=0
$
for all sufficiently large $k$, proving the claim.
\end{proof}

%Under the same assumptions, the logarithmic threshold below also ensures finite identification and can be evaluated without knowing $\theta$.
Under the same assumptions, the logarithmic threshold also identifies the support in finitely many iterations. Its computation requires no value for $\theta$.
\begin{eqnarray*}\label{eq:unknown_KL_threshold}
\varepsilon_k:=\frac{C_{\varepsilon}}{\log(e+a_k^{-1})^{q_{\varepsilon}}},
\quad C_{\varepsilon}>0,\ q_{\varepsilon} > 0,
\end{eqnarray*}
\noindent Here $a_k$ is defined in \eqref{nu}.
Indeed, $\varepsilon_k>0$ and $\varepsilon_k\to 0$. For $\nu$ defined in (\ref{nu}), we have
\[
\frac{\|\mathbf{x}^k-\overline{\mathbf{x}}\|}{\varepsilon_k}
\leq\frac{\widehat C}{C_{\varepsilon}}
{a_k^\nu[\log(e+a_k^{-1})]^{q_{\varepsilon}}}\longrightarrow0.
\]
Consequently, the selected set $\mathcal S^k$ eventually equals $\operatorname{supp}(\overline{\mathbf x})$ without knowledge of the KL exponent. The floor ${C_a}/({k+1})$ in $a_k$ ensures that $a_k>0$, and hence $\varepsilon_k>0$.

\subsection{Finite Identification without Nondegeneracy: A Degenerate Example}
\label{sec:degenerate-example}

The following example shows that proximal-gradient methods may fail to
identify the support of the limit in finitely many iterations when the
nondegeneracy condition fails, even if the objective is partly smooth.
In the same example, thresholding with $\varepsilon_k$ chosen as in
Theorem~\ref{thm:Finite_identify_KL_general} recovers the limiting support
in finitely many iterations. For simplicity, we apply the selection rule
within the proximal-gradient framework, setting
$\mathbf{y}^k=\mathbf{x}^k$ for all $k$.

\begin{example}
	Consider the problem
	\[
	\min_{\mathbf{x}\in\mathbb{R}^2}
	F(\mathbf{x}):=\|\mathbf{x}\|_1+\frac12\|\mathbf{x}-\mathbf{b}\|^2,
	\qquad \mathbf{b}=(3,1)^\top.
	\]
	Let $g(\mathbf{x}):=\|\mathbf{x}\|_1$ and
	$h(\mathbf{x}):=\frac12\|\mathbf{x}-\mathbf{b}\|^2$.
	Then $\nabla h$ is $1$-Lipschitz continuous, and the unique
	minimizer is $\overline{\mathbf{x}}=(2,0)^\top$.
	{Moreover, $F$ is partly smooth \citep{lewis2002active}
	at every $\mathbf{x}\in\mathbb{R}^2$ relative to $\mathcal{M}_{\mathbf{x}} :=\{\mathbf{z}:\supp(\mathbf{z})=\supp(\mathbf{x})\}$,
	however, the nondegeneracy condition fails at $\overline{\mathbf{x}}$ since}
	\[
	-\nabla h(\overline{\mathbf{x}})=(1,1)^\top
	\notin\ri\partial g(\overline{\mathbf{x}})
	=\{1\}\times(-1,1).
	\]
	
	We first consider the standard PG method.
	For any $\alpha\in(0,1)$ and
	$\mathbf{x}^0\in(1,3)\times(0,1)$, let $r:=1-\alpha$.
	Soft thresholding gives
	$x_1^k-2=r(x_1^{k-1}-2)$ and $x_2^k=rx_2^{k-1}$.
	Since these updates preserve $(1,3)\times(0,1)$, we obtain
	\[
	\mathbf{x}^k
	=\prox_{\alpha g}\bigl(\mathbf{x}^{k-1}
	-\alpha\nabla h(\mathbf{x}^{k-1})\bigr)
	=\begin{pmatrix}
		2+r^k(x_1^0-2)\\
		r^k x_2^0
	\end{pmatrix},
	\qquad k\ge1.
	\]
	Thus, although $\mathbf{x}^k\to\overline{\mathbf{x}}$,
	we have $\supp(\mathbf{x}^k)=\{1,2\}$ for every finite $k$.
	Hence, PG does not identify the support of the limit in finitely many iterations.
	
	We next apply a thresholding strategy to the PG iterates.
	Since $F$ is strongly convex, it has KL exponent $\theta=1/2$
	at $\overline{\mathbf{x}}$.
	Consider a proximal-gradient variant of the power threshold in Theorem~\ref{thm:Finite_identify_KL_general}.
	For any $\gamma\in(0,1)$, let
	\[
	\mathcal{S}^k:=\{i\in\{1,2\}:|x_i^k|\ge\varepsilon_k\},
	\qquad k\ge1,
	\]
	and $C>0$, we obtain
	\[
	\varepsilon_k=C^\gamma r^{\gamma(k-1)}>0,
	\qquad
	\frac{\|\mathbf{x}^k-\overline{\mathbf{x}}\|^2}{\varepsilon_k^2}
	=C^{2(1-\gamma)}r^{2k(1-\gamma)+2\gamma}
	\longrightarrow0.
	\]
	Since $\varepsilon_k>0$ for every $k\ge1$,
	no additional positive floor is needed.
	Moreover, since
	$|x_2^k|\le\|\mathbf{x}^k-\overline{\mathbf{x}}\|$,
	$\varepsilon_k\to0$, and $x_1^k\to2$,
	there exists $K$ such that
	\[
	|x_2^k|<\varepsilon_k<|x_1^k|,
	\qquad
	\mathcal{S}^k=\{1\}=\supp(\overline{\mathbf{x}}),
	\qquad k\ge K.
	\]
	Thus, the thresholding strategy identifies the limiting support
	in finitely many iterations, even though the PG iterates do not.
\end{example}

\section{Restricted Local Stationarity}\label{sec:local}
{
Throughout this section, we assume that the conditions of
Theorem~\ref{thm:global_convegence} hold and that there exists an index $K_1(\ge K_0)$
such that
$\mathcal S^k=\overline{\mathcal S}:=\supp(\overline{\x})$
for all $k\ge K_1$.
Theorem~\ref{thm:Finite_identify_KL_general} gives this eventual equality of supports.
{We restrict the analysis in this section to $\overline{\mathcal S}\neq\emptyset$. If $\overline{\mathcal S}=\emptyset$, the assumed
identification implies that $\mathcal S^k=\emptyset$
for all sufficiently large $k$. Consequently, $\mathbf y^k=\mathbf x^k$, and no cubic refinement is performed.}
For $k\geq K_1$, let $\widehat{\mathbf d}^k\in\mathbb R^n$ denote the zero extension of $\mathbf d^k$ from $\overline{\mathcal S}$.
For $k\ge K_1$, we denote
\begin{equation*}
\g^k=\nabla_{\overline{\mathcal S}}F(\x^k),
\quad
H^k=\nabla^2_{\overline{\mathcal S}}F(\x^k).
\end{equation*}
%When appearing in full-space expressions, $\mathbf d^k$
%is understood as the zero extension of the direction computed
%on $\overline{\mathcal S}$.
We also require continuity of the restricted Hessian near the limit.

{
\begin{assumption}\label{ass:partial_Hessian_continuous}
For each $\x\in\dom\partial F$ with
$\mathcal S=\supp(\x)$,
there exists $\delta_{\x}>0$ such that
$\nabla^2_{\mathcal S}F(\z)$ is well-defined
and continuous on $\cB(\x,\delta_{\x})$.
\end{assumption}
{The next theorem proves eventual acceptance of full cubic}
steps, that is, $\beta_k=1$ for all sufficiently large $k$.}

\begin{theorem}\label{lemma:beta_eq1}
Under the hypotheses stated at the beginning of this section and
Assumption~\ref{ass:partial_Hessian_continuous}, the following assertions
hold.
\begin{enumerate}
\item[(i)] The sequences $\{\g^k\}$, $\{H^k\}$ and
$\{\mathbf d^k\}$ are bounded for $k\ge K_1$.
\item[(ii)] There is a constant $\underline\beta>0$ such that
$\beta_k\ge\underline\beta$ for all sufficiently large $k$.
\item[(iii)] $\mathbf d^k\to0$.
{\item[(iv)] {There is an index $\widehat K\ge K_1$ such that, for every $k\geq\widehat K$, $\beta_k=1$ and the following bound holds for all $t\in[0,1]$}}
	\begin{equation}\label{key2}
		\|\nabla^2_{\overline{\mathcal S}}F(\x^k+t\widehat{\mathbf d}^k)
		-\nabla^2_{\overline{\mathcal S}}F(\x^k)\|
		\le\frac{\eta}{2}+M(\overline{\x})t\|\mathbf d^k\|,
	\end{equation}
	{where any $M(\overline{\x})\ge0$ is valid, including}
	$M(\overline{\x})=(\underline\omega+1)
	\mathrm{lip}\,\nabla^2_{\overline{\mathcal S}}F(\overline{\x})/2$.
	\end{enumerate}
\end{theorem}
\begin{proof}
{Choose a sufficiently small ball centered at \(\overline{\mathbf{x}}\) and contained in \(\cB(\overline{\mathbf{x}},\delta_{\overline{\mathbf{x}}})\)}, such that all
coordinates in $\overline{\mathcal S}$ remain nonzero throughout
the ball. By continuity, the restricted Hessian is bounded there
by some constant $C_H>0$.
{Since $\x^k\to\overline{\x}$, there exists $r>0$ such that, for all sufficiently large $k$, the line segment joining $\x^k$ and $\x^k+r\mathbf d$ lies in this ball for any unit vector $\mathbf d$ supported on $\overline{\mathcal S}$.}
For any $\g^k\neq\mathbf0$, choosing the direction $-\g^k/\|\g^k\|$ and applying Taylor's theorem yields
\begin{equation*}
\underline F
\le F\left(\x^k-r\frac{\g^k}{\|\g^k\|}\right)
\le F(\x^k)-r\|\g^k\|+\frac{C_Hr^2}{2},
\end{equation*}
{where the direction is understood to be {extended by zero outside} $\overline{\mathcal S}$.}
Since $\{F(\x^k)\}$ is bounded above, it follows that
$\{\g^k\}$ is bounded. Enlarging the bound if necessary
to cover the finitely many remaining indices, we may choose
$C_g>0$ such that $\|\g^k\|\le C_g$ for all $k$.
{Now, let $v_k=\|\mathbf d^k\|$. By the definition of the
residual $\mathbf r^k$ and Condition~\eqref{inexact1},
we obtain
\begin{equation*}
	\begin{aligned}
		\langle\g^k,\mathbf d^k\rangle
		+\langle H^k\mathbf d^k,\mathbf d^k\rangle
		+\eta v_k^2+\frac{\sigma_k}{2}v_k^3
		=\langle\mathbf r^k,\mathbf d^k\rangle
		\le 0.
	\end{aligned}
\end{equation*}
Hence the Cauchy--Schwarz inequality and the bounds
$\|\g^k\|\le C_g$, $\|H^k\|\le C_H$ yield
\[
\begin{aligned}
	\frac{\sigma_k}{2}v_k^3
	&\le-\langle\g^k,\mathbf d^k\rangle
	-\langle H^k\mathbf d^k,\mathbf d^k\rangle
	-\eta v_k^2\\
	&\le|\langle\g^k,\mathbf d^k\rangle|
	+|\langle H^k\mathbf d^k,\mathbf d^k\rangle|\\
	&\le C_gv_k+C_Hv_k^2.
\end{aligned}
\]
Since $\sigma_k\ge2c/3$, for $v_k>0$ the estimate gives
$
\frac{c}{3}v_k^2\le C_g+C_Hv_k.
$
Hence, there exists $C_d >0$ such that $v_k\leq C_d$, which also holds trivially when $v_k=0$. This proves (i).}
By continuity of $\nabla^2_{\overline{\mathcal S}}F$ at $\overline{\x}$, there exists $\rho>0$ such that
\[
\|\nabla^2_{\overline{\mathcal S}}F(\z)
-\nabla^2_{\overline{\mathcal S}}F(\mathbf w)\|\le\eta/2
\quad(\z,\mathbf w\in\cB(\overline{\x},\rho)),
\]
and the coordinates indexed by $\overline{\mathcal{S}}$ remain nonzero on this ball. Let $\widehat\beta=\min\{1,\rho/(2C_d)\}$.
For all sufficiently large $k$, $\x^k\in\cB(\overline{\x},\rho/2)$,
and the segments $\x^k+t\beta\widehat{\mathbf d}^k$, $t\in[0,1]$,
remain in the ball whenever $0<\beta\le\widehat\beta$. The complementary
coordinates are fixed, so the support is unchanged along each segment.
Lemma~\ref{lem:inexact_cubic_upper}, with $a=\eta/2$ and $M=0$, and
Conditions~\eqref{inexact1}--\eqref{inexact2} therefore yield
\[
\begin{aligned}
F(\x^k+\beta\widehat{\mathbf d}^k)
&\le F(\x^k)-\frac{\eta\beta^2}{4}v_k^2
-\frac{\sigma_k\beta^2}{4}v_k^3\le F(\x^k)-\frac{\eta\beta^2}{4}v_k^2
-\frac{c\beta^3}{6}v_k^3.
\end{aligned}
\]
Thus every such $\beta$ satisfies the line-search test. Backtracking
gives $\beta_k\ge q\widehat\beta$, proving (ii).
{Since $\|\mathbf{d}^k\| = \|\x^k-\y^k\|/\beta_k \leq \|\x^k-\y^k\|/\underline{\beta} \to0$, assertion (iii) follows.}
{Since $\mathbf d^k\to\mathbf0$ and $\mathbf x^k\to\overline{\mathbf x}$,
the full trial segment lies in the ball for all sufficiently large $k$.
The same estimate gives $\beta_k=1$, and the Hessian variation is at most
$\eta/2$ throughout the segment. This proves \eqref{key2} for any
nonnegative choice of $M(\overline{\mathbf x})$.}
\end{proof}

\begin{definition}
For a twice continuously differentiable function $\widetilde F$,
a point $\x$ is an $(\varepsilon_g,\varepsilon_H)$-approximate
second-order stationary point if
\begin{eqnarray}\label{apprsec}
\|\nabla\widetilde F(\x)\|\le\varepsilon_g,
\qquad
\lambda_{\min}(\nabla^2\widetilde F(\x))
\ge-\sqrt{\varepsilon_H}.
\end{eqnarray}
For the nonsmooth objective in \eqref{Problem}, replace $\nabla\widetilde F$ and $\nabla^2\widetilde F$ in
\eqref{apprsec} with $\nabla_{\cal S} F$ and $\nabla^2_{\cal S} F$ to define approximate second-order stationarity on
${\cal S}$.
\end{definition}
\begin{corollary}\label{thm:approximate_second_order_stationary}
Suppose that the hypotheses of Theorem~\ref{lemma:beta_eq1} hold and
$\sup_k\sigma_k=:\overline\sigma<\infty$. Then
\[
\nabla_{\overline{\mathcal S}}F(\overline{\x})=0,\qquad
\nabla^2_{\overline{\mathcal S}}F(\overline{\x})\succeq-\eta I.
\]
\end{corollary}
\begin{proof}
{Let $v_k:=\|\mathbf d^k\|$.
By conditions (\ref{inexact3}) and (\ref{inexact4}), for all sufficiently large $k$,
\[
\|\mathbf g^k\|
\leq (C_H+\eta)v_k
+\left(\frac{\overline{\sigma}}{2}+\kappa_g\right)v_k^2,
\qquad
H^k\succeq
-\left[
\eta+
\left(\frac{\overline{\sigma}}{2}+\kappa_H\right)v_k
\right]I.
\]
{Theorem~\ref{lemma:beta_eq1}(iii)} gives $v_k\to0$, and hence
$\mathbf g^k\to0$.
By the continuity of the restricted Hessian,
$H^k\to\nabla_{\overline{\mathcal S}}^2
F(\overline{\mathbf x})$, which yields the asserted
Hessian inequality.}
{To justify the gradient conclusion, fix a unit vector $\mathbf v$ supported on $\overline{\mathcal S}$.} For sufficiently small $t>0$,
Taylor's theorem gives
$
F(\x^k+t\mathbf v)\le F(\x^k)+t\langle\g^k,\mathbf v_{\overline{\mathcal S}}\rangle
+C_Ht^2/2.
$
The global convergence analysis gives $F(\x^k)\to F(\overline{\x})$.
Using lower semicontinuity and $\g^k\to0$, we obtain
$F(\overline{\x}+t\mathbf v)\le F(\overline{\x})+C_Ht^2/2$.
Divide by $t$ and let $t\downarrow0$. Repeating the argument with
$-\mathbf v$ proves $\nabla_{\overline{\mathcal S}}F(\overline{\x})=0$.
\end{proof}

Let $\g^k_{\y}=\nabla_{\overline{\mathcal S}}F(\y^k)$ and
$H^k_{\y}=\nabla^2_{\overline{\mathcal S}}F(\y^k)$.
The following theorem provides an upper bound on the
number of outer iterations required to reach an
approximate second-order stationary point on the
identified subspace, starting from the iteration
at which the selected subspace has stabilized and
full cubic steps are always accepted.
%For a fixed run, the iteration bound starts once the selected subspace
%has stabilized and the line search accepts full steps.%The index at which this regime begins and the constants in the bound may depend on the run.

\begin{theorem}\label{thm:local_complexity}
Suppose that the assumptions of Corollary~\ref{thm:approximate_second_order_stationary}
hold. Define
\[
A=M(\overline{\x})+\overline\sigma+2\kappa_g,\qquad
B=\frac32M(\overline{\x})+\overline\sigma+\kappa_g+\kappa_H.
\]
Given $\varepsilon_g,\varepsilon_H>0$, suppose that the parameter $\eta$ satisfies
\[
0<\eta\le\min\left\{\frac{\sqrt{A\varepsilon_g}}{3},
\frac{A\sqrt{\varepsilon_H}}{3B}\right\}.
\]
Let $T$ be the first $k\ge\widehat K$ for which
$\|\g^k_{\y}\|\le\varepsilon_g$ and
$\lambda_{\min}(H^k_{\y})\ge-\sqrt{\varepsilon_H}$.
Then $T$ is finite and
\begin{equation}\label{eq:local_iteration_count}
T-\widehat K\le
\left\lceil\frac{6(F(\x^{\widehat K})-\underline F)}
{c\min\{\sqrt{\varepsilon_g/A},\sqrt{\varepsilon_H}/B\}^{3}}
\right\rceil.
\end{equation}
\end{theorem}
\begin{proof}
{By Theorems~\ref{thm:global_convegence} and~\ref{lemma:beta_eq1}, we have}
$\x^k\to\overline{\x}$, $\mathbf d^k\to0$,
and $\beta_k=1$ for all sufficiently large $k$.
{Recall $v_k=\|\mathbf d^k\|$.}
{{By \eqref{key2} and conditions~\eqref{inexact3}--\eqref{inexact4}, we obtain}}
{
\begin{equation*}\label{ieq:gy_bound}
	\begin{aligned}
	\|\mathbf g_y^k\|
	&\leq
	\left(
	\frac{M(\overline{\mathbf x})+\overline{\sigma}}{2}
	+\kappa_g
	\right)v_k^2
	+\frac{3\eta}{2}v_k,\\
	H_y^k&\succeq
	-\left[
	\frac{3\eta}{2}
	+
	\left(
	M(\overline{\mathbf x})
	+\frac{\overline{\sigma}}{2}
	+\kappa_H
	\right)v_k
	\right]I.
	\end{aligned}
\end{equation*}
}
%	and
%	\begin{equation}\label{ieq:Hy_bound}
%	
%	\end{equation}
If $v_k\le3\eta/A$, these estimates imply
\[
\|\g^k_{\y}\|\le9\eta^2/A\le\varepsilon_g,\quad
H^k_{\y}\succeq-(3\eta B/A)I\succeq-\sqrt{\varepsilon_H}I,
\]
{which implies $\mathbf{y}^k$ is already an $(\varepsilon_g,\varepsilon_H)$-approximate second-order stationary point
on
$\overline{\cal S}:={\text{supp}}(\overline{\h x})$.}
Otherwise $\eta<Av_k/3$, so
\[
\|\g^k_{\y}\|<Av_k^2,\quad H^k_{\y}\succeq-Bv_kI.
\]
Consequently, every iteration not yet satisfying the target tolerances satisfies
$
v_k>\delta:=\min\{\sqrt{\varepsilon_g/A},\sqrt{\varepsilon_H}/B\}.
$
Use the gradient bound when the gradient test fails or the Hessian
bound when the curvature test fails.
{Since the full Newton step is accepted for all $k\ge\widehat K$,
the subsequent proximal-gradient and {subspace cubic Newton} steps yield
\begin{equation*}
	F(\x^{k+1})\le F(\y^k)
	\le F(\x^k)-\frac c6v_k^3
	< F(\x^k)-\frac c6\delta^3.
\end{equation*}
Since $F$ is bounded below, the target tolerances must be satisfied
after finitely many iterations. Summing the above inequality from
$k=\widehat K$ to $T-1$ yields \eqref{eq:local_iteration_count}.}
\end{proof}

\section{Numerical Results}\label{sec:numerical}

AHPG-SCN is an adaptive variant of HPG-SCN (see Algorithm~\ref{Ag2}). It sets the cubic regularization parameter using the ratio of actual to predicted reduction, so it can dispense with estimating the local Hessian Lipschitz modulus.
\begin{itemize}
%\item We employ an adaptive line-search strategy for the proximal-gradient step, with an initial stepsize estimated from local curvature information using~\eqref{eq:adaptive_initial_stepsize} and then reduced until a sufficient decrease condition is satisfied (see Lines~2--10). This strategy avoids the need to prescribe a conservative fixed stepsize based on the global Lipschitz constant.

{\item AHPG-SCN takes a cubic step only if $1\le n_k\le n_{\max}$.}

{\item The ratio $\rho_k=(F(\x^k)-F(\y^k))/(F(\x^k)-m_k(\beta_k\mathbf d^k))$ compares actual and predicted reductions \citep{cartis2011adaptive}. The method decreases $\sigma_k$ if $\rho_k\ge\eta_2$, and rejects the step and increases $\sigma_k$ if $\rho_k\le\eta_1$, subject to the prescribed bounds. For $0<\eta_1<\rho_k<\eta_2<1$, $\sigma_k$ stays fixed.}
\end{itemize}

\begin{algorithm}[htbp]
	\footnotesize
	\setlength{\baselineskip}{8.8pt}
	\caption{Adaptive hybrid proximal-gradient and subspace cubic-Newton
		method (AHPG-SCN) for~\eqref{Problem}}
	\label{Ag2}
	\begin{algorithmic}[1]
		\Require{
			$\mathbf x^0=\mathbf y^0\in\operatorname{dom}F$,
			$0<\alpha_0<\overline{\alpha}<\lambda_g$,
			{$c,\ c_1>0$,}
			$\eta>0$,
			{$q,\ q_1\in(0,1)$,}
			$\varepsilon_k>0$,
			$0<\eta_1<\eta_2<1$,
			$0<\gamma_1<1<\gamma_2$,
			$0<\sigma_{\min}\leq\sigma_1\leq\sigma_{\max}$,
			$n_{\max}>0$,
			{$\kappa_g,\ \kappa_H\geq0$,}
			\texttt{InnerMaxIt},
			\texttt{MaxIt}.
		}
		
		\For{$k=1,2,\ldots,\texttt{MaxIt}$}
		
		\State \textbf{Proximal-gradient step.}
		\State Set $\alpha_{k,0}$ according
		to~\eqref{eq:adaptive_initial_stepsize}.
		
		\For{$l_k=0,1,2,\ldots$}
		\State Compute
		$
		\mathbf x^{k,l_k}
		\in
		\prox_{\alpha_{k,l_k}g}
		\bigl(
		\mathbf y^{k-1}
		-\alpha_{k,l_k}
		\nabla h(\mathbf y^{k-1})
		\bigr)$.
		\If{$
			F(\mathbf x^{k,l_k})
			\leq
			F(\mathbf y^{k-1})
			-c_1
			\|\mathbf x^{k,l_k}-\mathbf y^{k-1}\|^2
			$}
		\State Set
		$\mathbf x^k=\mathbf x^{k,l_k}$ and
		$\alpha_k=\alpha_{k,l_k}$;
		\textbf{break}.
		\Else
		\State Set
		$\alpha_{k,l_k+1}=q_1\alpha_{k,l_k}$.
		\EndIf
		\EndFor
		
		\State \textbf{Coordinate selection.}
		\State Set
		$\mathcal S^k
		:=\{\,i:|x_i^k|\geq\varepsilon_k\,\}$
		and $n_k:=\sharp(\mathcal S^k)$.
		
		\If{$\mathcal S^k\neq\emptyset$ and
			$n_k\leq n_{\max}$}
		
		\State \textbf{Subspace cubic-Newton step.}
		\State Compute
		$\mathbf g^k
		=\nabla_{\mathcal S^k}F(\mathbf x^k)$
		and
		$H^k
		=\nabla^2_{\mathcal S^k}F(\mathbf x^k)$.
		
		\State Compute an approximate solution
		$\mathbf d^k$ of~\eqref{subd}
		satisfying
		\eqref{inexact1}--\eqref{inexact4}.
		
		\State \textbf{Backtracking line search.}
		\State Set $\beta_{k,0}=1$ and
		$\texttt{Tag}=0$.
		
		\For{$j_k=0,1,2,\ldots,\texttt{InnerMaxIt}$}
		\State Define
		$
		\mathbf y^{k,j_k}_{\mathcal S^k}
		=
		\mathbf x^k_{\mathcal S^k}
		+\beta_{k,j_k}\mathbf d^k,
		\quad
		\mathbf y^{k,j_k}_{(\mathcal S^k)^c}
		=
		\mathbf x^k_{(\mathcal S^k)^c}$.
		\If{$
			F(\mathbf y^{k,j_k})
			\leq
			F(\mathbf x^k)
			-\frac{\eta}{4}
			\|\mathbf y^{k,j_k}-\mathbf x^k\|^2
			-\frac{c}{6}
			\|\mathbf y^{k,j_k}-\mathbf x^k\|^3
			$}
		\State Set
		$\mathbf y^k=\mathbf y^{k,j_k}$,
		$\beta_k=\beta_{k,j_k}$,
		and $\texttt{Tag}=1$;
		\textbf{break}.
		\Else
		\State Set
		$\beta_{k,j_k+1}=q\beta_{k,j_k}$.
		\EndIf
		\EndFor
		
		\If{$\texttt{Tag}=1$}
		
		\State \textbf{Adaptive regularization update.}
		\State Define
		$
			m_k(\mathbf d)
			=F(\mathbf x^k)
			+\langle\mathbf g^k,\mathbf d\rangle+\frac12
			\langle(H^k+\eta I)\mathbf d,
			\mathbf d\rangle
			+\frac{\sigma_k}{6}\|\mathbf d\|^3$.
		
		\If{$\mathbf d^k=0$}
		\State Set
		$\mathbf y^k=\mathbf x^k$ and
		$\sigma_{k+1}=\sigma_k$.
		\Else
		\State Compute
		$
		\rho_k
		=
		(F(\mathbf x^k)-F(\mathbf y^k))/(F(\mathbf x^k)-m_k(\beta_k\mathbf d^k)).
		$
		
		\If{$\rho_k\geq\eta_2$}
		\State Set
		$\sigma_{k+1}
		=\max\{\sigma_{\min},
		\gamma_1\sigma_k\}$.
		
		\ElsIf{$\rho_k\leq\eta_1$}
		\State Set
		$\sigma_{k+1}
		=\min\{\sigma_{\max},
		\gamma_2\sigma_k\}$
		and $\mathbf y^k=\mathbf x^k$.
		
		\Else
		\State Set
		$\sigma_{k+1}=\sigma_k$.
		\EndIf
		\EndIf
		
		\Else
		\State Set
		$\sigma_{k+1}
		=\min\{\sigma_{\max},
		\gamma_2\sigma_k\}$
		and $\mathbf y^k=\mathbf x^k$.
		\EndIf
		
		\Else
		\State Set
		$\mathbf y^k=\mathbf x^k$ and
		$\sigma_{k+1}=\sigma_k$.
		\EndIf
		
		\If{the termination criterion is satisfied}
		\State \textbf{break}.
		\EndIf
		
		\EndFor
		
		\State \textbf{Output:} $\mathbf y^k$.
	\end{algorithmic}
\end{algorithm}
The convergence and complexity analysis covers HPG-SCN. We leave the convergence analysis of AHPG-SCN for future work.
We compare HPG-SCN and AHPG-SCN with ZeroFPR, ManAcc, and HpgSRN. We tuned the parameters of HPG-SCN and AHPG-SCN by hand and list the settings for each method below.
{\begin{itemize}
		\item \textbf{ZeroFPR.}
		We follow the algorithm and BFGS update given by (5.9) in
		\citep{themelis2018forward}, using the publicly available
		MATLAB implementation from the ForBES package\footnote{\url{https://github.com/kul-optec/ForBES}} with its default algorithmic parameters.
		\item \textbf{ManAcc.}
		We follow Algorithms~1 and 3 of \citep{bareilles2023newton}.
		Since the authors' implementation is available in
		Julia\footnote{\url{https://github.com/GillesBareilles/NewtonRiemannAccel-ProxGrad}}, we reimplement the method in MATLAB. We use the truncated Newton variant and retain the original parameter settings.
		The only modification is the use of the shifted Hessian $H^k+\eta I$, with $\eta=10^{-3}$, to improve numerical stability.
		\item \textbf{HpgSRN.}
		We follow Algorithm~3.2 of \citep{wu2023regularized}, using the authors'
		publicly available MATLAB code\footnote{\url{https://github.com/YuqiaWU/HpgSRN}} with its default algorithmic parameters.
		\item \textbf{HPG-SCN.}
		We set $\alpha_0=0.95/L$,
		$\overline{\alpha}=10^6$,
		$c=c_1=10^{-3}$, $\eta=10^{-10}$,
		$q=q_1=0.5$, and $\kappa_g=\kappa_H=10^{-2}$.
		The cubic regularization parameter is chosen as
		$\sigma_k=0.67(M(\mathbf{x}^k)+c)$.
		The Lanczos solver uses $s=1.5$ and
		$r_{k,0}=\lceil\sqrt{n_k}\rceil$.
		The curvature estimate $M(\mathbf{x}^k)$ and
		the threshold parameters $C_\varepsilon$ and
		$q_\varepsilon$ are specified for each problem.
		\item \textbf{AHPG-SCN.}
		We use the same settings as HPG-SCN, except
		that the cubic regularization parameter is
		updated adaptively.
		We set $\eta_1=0.1$, $\eta_2=0.9$,
		$\gamma_1=0.5$, $\gamma_2=2$,
		$\sigma_{\min}=10^{-3}$,
		$\sigma_1=1$, $\sigma_{\max}=10^3$,
		and $\texttt{InnerMaxIt}=50$.
		The maximum subspace dimension $n_{\max}$
		is specified in each experiment.
\end{itemize}}

We chose these solvers because each combines proximal-gradient-type steps with second-order or quasi-Newton acceleration.
ZeroFPR employs quasi-Newton directions based on the forward--backward envelope, while ManAcc combines proximal-gradient steps with Riemannian Newton-type acceleration on an identified manifold. For the sparse models considered here, the manifold is determined by the current support. HpgSRN combines proximal-gradient steps with subspace regularized Newton steps and is specifically designed for nonconvex $L_p^p$-regularized optimization.
As in \cite{wu2023regularized}, we define
\begin{equation*}
\operatorname{Res}(\mathbf{x})
:=
\gamma_{\mathrm{res}}
\operatorname{dist}_{\infty}
\left(
\mathbf{x},
\operatorname{prox}_{g/\gamma_{\mathrm{res}}}
\left(
\mathbf{x}
-\gamma_{\mathrm{res}}^{-1}
\nabla h(\mathbf{x})
\right)
\right),
\end{equation*}
where $\gamma_{\mathrm{res}}=\frac{L}{0.95}$, and $\operatorname{dist}_{\infty}$ denotes the distance induced by the infinity norm.
All solvers stop as soon as any of the following conditions holds.
\begin{equation*}
\mathrm{Res}(\mathbf{x})<10^{-8},\quad
k\ge\texttt{MaxIt},\quad\text{or}\quad
\text{CPU time}>200~\text{seconds},
\end{equation*}
where ${\texttt{MaxIt}}=10^5$.
We report means over 20 independent trials.

\subsection{Known KL Exponent}
\label{sec:toy_example}
{We test the two threshold rules on the following problem.}
\begin{equation*}
\min_{\mathbf{x}\in\mathbb{R}^n}
F(\mathbf{x})
:=
\sum_{i=1}^n
\left(\log\left(1+(x_i-b_i)^4\right)-x_i\right)
+\|\mathbf{x}\|_1,
\end{equation*}
where $\mathbf{b}\in\{0,1\}^n$ contains both zero and nonzero components. The objective $F$ is coercive and admits a unique
global solution $\overline{\mathbf{x}}=\mathbf{b}$, which is also its unique critical point.
Taking $g(\mathbf{x})=\|\mathbf{x}\|_1$ and $h(\mathbf{x})=\sum_{i=1}^n
(\log(1+(x_i-b_i)^4)-x_i)$, we can verify that $F$ satisfies Assumption~\ref{Ass1} with $L<4$ and $\operatorname{lip}_{\mathcal M_{\mathbf{x}}} \nabla_{\mathcal S}^2F(\mathbf{x})<9$
for any $\mathbf{x}\neq\mathbf{0}$ and $\emptyset\neq\mathcal S\subseteq\operatorname{supp}(\mathbf{x})$.
Moreover, $F$ has the KL property at $\overline{\mathbf{x}}$ with exponent $\theta=3/4$, while the standard nondegeneracy condition
$-\nabla h(\overline{\mathbf{x}}) \in\operatorname{ri}\partial g(\overline{\mathbf{x}})$ {fails (see Proposition~\ref{prop:toy_example}).}
We compare the sets selected by the two HPG-SCN threshold rules with the supports of the proximal-gradient iterates. The proximal-gradient baseline uses the update below with a line search (see Section~10.3.3 of \cite{beck2017first}).
	\begin{equation*}
		\x^k \in \underset{\mathbf{x}}{\arg\min}
		\left\{
		g(\mathbf{x}) +
		\frac{1}{2\alpha_k}
		\left\|\mathbf{x} -
		\left(\x^{k-1}-\alpha_k\nabla h(\x^{k-1})\right)\right\|^2
		\right\},
	\end{equation*}
	where $\alpha_k$ is {chosen by the same line-search procedure as in} HPG-SCN.

{We set $n=50$ and construct $\mathbf{b}$ by setting its first $10$ components to one and the remaining components to zero so that $\overline{\mathcal S}=[10]$. We compare PG with HPG-SCN equipped with two thresholding strategies:
(i) the KL-exponent-guided power-type threshold
$\varepsilon_k =0.5(\max\{s_k,10^{-8}/(k+1)\})^\gamma$, where $s_k$ is defined in
Theorem~\ref{thm:Finite_identify_KL_general},
with $\gamma = 1/6$;
and (ii) the adaptive threshold
$\varepsilon_k
=2/\log^2(e+\max\{s_k,10^{-8}/(k+1)\}^{-1})$.}

%{To illustrate that finite identification may fail in the absence of the nondegeneracy condition numerically, we consider the standard proximal-gradient (PG) method \citep{beck2017first} as a baseline, whose iterates are generated by

{For both methods, we set $\alpha_0=1/5$.
	For HPG-SCN, we choose $\sigma_k=0.67(M(\mathbf{x}^k)+c)$ with $M(\mathbf{x}^k) = 9$.
	The initial point is independently sampled from the uniform distribution on $[-3,3]$.}
%We set $\texttt{RTol}=10^{-8},\texttt{MaxIt}=100000$ and $\texttt{CTol} = 200$ in (\ref{stoppping}).}

{We record $n_k$ for HPG-SCN, $\|\mathbf{x}^k\|_0$ for PG (also denoted by $n_k$) and the Jaccard index \cite{liu2024jaccard} $ J_k= {\sharp(\mathcal S^k\cap\overline{\mathcal S})}/{\sharp(\mathcal S^k\cup\overline{\mathcal S})},$ where $\mathcal S^k=\operatorname{supp}(\mathbf{x}^k)$
for PG. Note that $J_k=1$ if and only if $\mathcal S^k=\overline{\mathcal S}$.}

\begin{figure}[htpb]
\centering
\begin{tabular}{@{}c@{\hspace{.025\linewidth}}c@{}}
%\toprule
%$(\text{ODCT},10,20)$ & $(\text{Gauss},10,0.3)$ \\
%\midrule
\includegraphics[width=.48\linewidth]{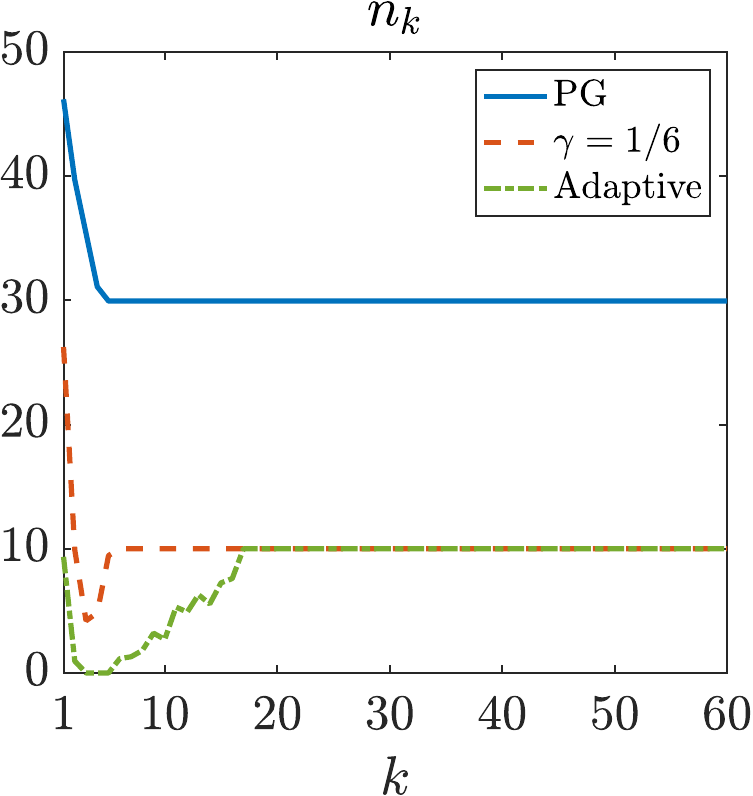} &
\includegraphics[width=.48\linewidth]{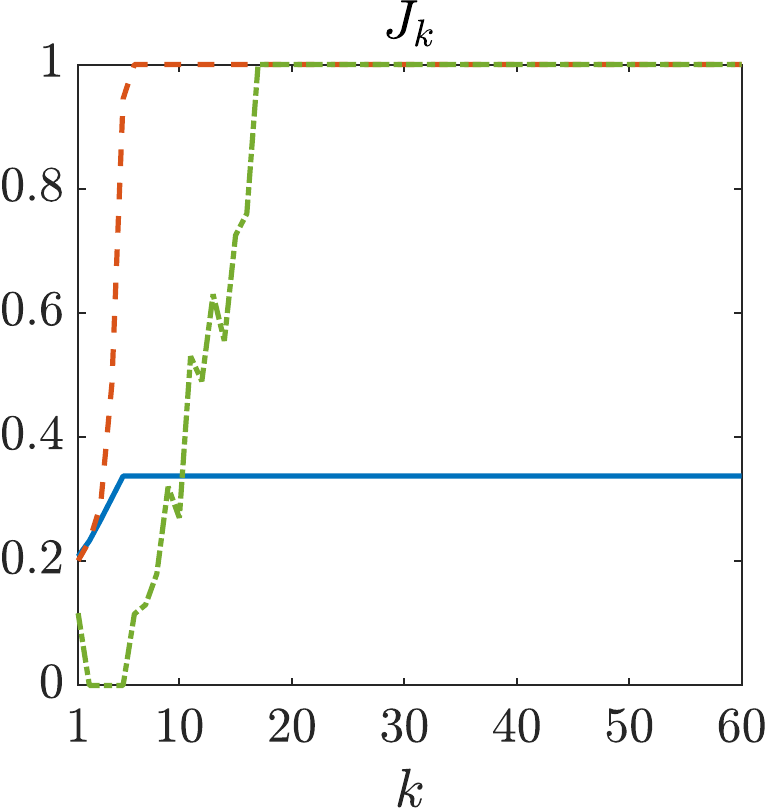} %& \includegraphics[scale=.36]{APG_CN_figure/toy_example_threshold.eps}
%\\
%\bottomrule
\end{tabular}
\caption{Averaged $n_k$ and Jaccard index.}
\label{fig:toy_example}
\end{figure}

{Figure~\ref{fig:toy_example} reports the mean selected-set size and Jaccard index over the first $60$ iterations.
For both the power-type and logarithmic thresholds, the mean selected-set size reaches $10$ and the mean Jaccard index reaches $1$.
This behavior is consistent with Theorem~\ref{thm:Finite_identify_KL_general} and its
logarithmic-threshold counterpart, which give finite identification
of $\mathcal S^k=\overline{\mathcal S}$ without a nondegeneracy assumption.
PG retains extra nonzero coordinates throughout the reported iterations
and therefore does not identify $\overline{\mathcal S}$.}

\subsection{Unknown KL Exponent}\label{sec:recovery}
\subsubsection{Sparse Signal Recovery}
{We consider the ratio-regularized recovery problem}
%proposed in \citep{Tao20}:
\begin{equation*}
\min_{\x\in\mathbb{R}^n}
\lambda\frac{\|\x\|_1}{\|\x\|_2}
+ \frac{1}{2}\|A\x-\mathbf{b}\|_2^2.
\label{prob:L1oL2_recovery}
\end{equation*}
%where $A\in\mathbb{R}^{m\times n}$ with $m\ll n$, $\mathbf{b}\in\mathbb{R}^{m}$, and $\lambda>0$.
%The objective function in (\ref{eq:intro_sparse_recovery}) is semi-algebraic on $\mathbb{R}^n\backslash\{\mathbf{0}\}$ and hence satisfies KL property \cite{attouch2010proximal}.
{With the convention $\|\mathbf 0\|_1/\|\mathbf 0\|_2=1$, the objective is proper, lower semicontinuous, and semialgebraic. It therefore satisfies the KL property \cite[Theorem 14]{attouch2010proximal}.}

We consider two types of sensing matrices.
(i) Gaussian matrix (Gauss).
{The rows of $A$ are sampled independently from $\mathcal N(0,\Sigma)$, where $\Sigma=(1-r)I_n+r\mathbf1_n\mathbf1_n^\top$.}
(ii) Oversampled DCT matrix (ODCT). Here $A = [{\bf a}_1,{\bf a}_2,\ldots,{\bf a}_n]\in {\mathbb R}^{m\times n}$
with each column $$ {\bf a}_j:=\frac{1}{\sqrt{m}}\cos \left( \frac{2\pi {\bf w}j}{G}\right),\;\;j=1,\ldots,n,$$
where ${\bf w}\in{\mathbb R}^m$ is a random vector with entries uniformly distributed on $[0,1]$ and
$G\in {\mathbb R}_+$ controls the coherence. A larger value of $G$ yields a more coherent matrix $A$.
We set $m=64$, $n=1024$. The ground truth $\x^*$ is $K$-sparse with
$K \in \{2,4,6,8,10\}.$ {The support of $\mathbf x^*$ is drawn uniformly among the $K$-element subsets of $\{1,\ldots,n\}$ and its nonzero entries are generated by}
\[
\texttt{x(supp) = sign(randn(K,1)).*1e2.*randn(K,1)}.
\]
We form the observation vector as $\mathbf{b}=A\x^*$.
We use a different set of configurations for each type of sensing matrix.
ODCT: $G \in \{15,20,25\}$,
Gaussian: $r \in \{0.1,0.2,0.3\}$.
We denote them by $(\mathrm{ODCT},K,G)$ and $(\mathrm{Gauss},K,r)$, respectively.

{We set $\lambda=10^{-2}$.
For HPG-SCN, we set $M(\mathbf{x}^k) = 256\lambda\sqrt{\|\x^k\|_0}/\|\x^k\|^3$ (Theorem~\ref{LHbound:L1oL2} permits the choice $\underline{\omega}=256/36>1$). For both HPG-SCN and AHPG-SCN we choose the threshold as $\varepsilon_k=C_{\varepsilon}/\log(e+\max\{s_k,10^{-8}/(k+1)\}^{-1})^{q_{\varepsilon}}$, with $C_{\varepsilon} = 150$ and $q_{\varepsilon} = 5.5$. We also set $n_{\max}=n$.}
All methods are initialized at the same point, obtained by applying
ADMM to the $\ell_1$-regularized least-squares problem with
$\lambda=10^{-3}$.
\begin{figure}[t]
\centering
\begin{tabular}{@{}c@{\hspace{.025\linewidth}}c@{}}
%\toprule
%ODCT & Gauss \\
%\midrule
%\includegraphics[scale=.32]{APG_CN_figure/perf_profile_cup_ODCT.eps} &
%\includegraphics[scale=.32]{APG_CN_figure/perf_profile_cup_Gauss.eps} \\
\includegraphics[width=.48\linewidth]{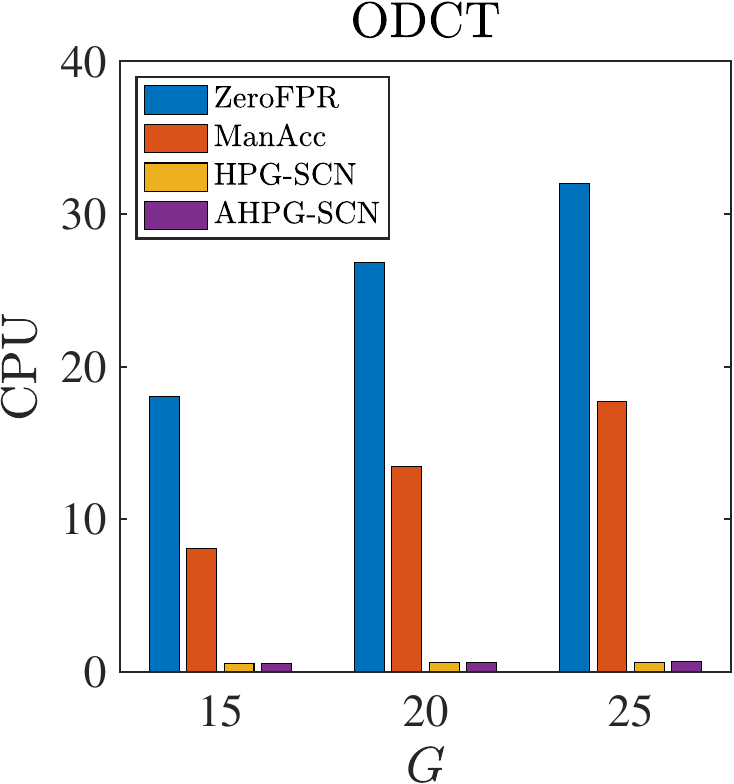} &
\includegraphics[width=.48\linewidth]{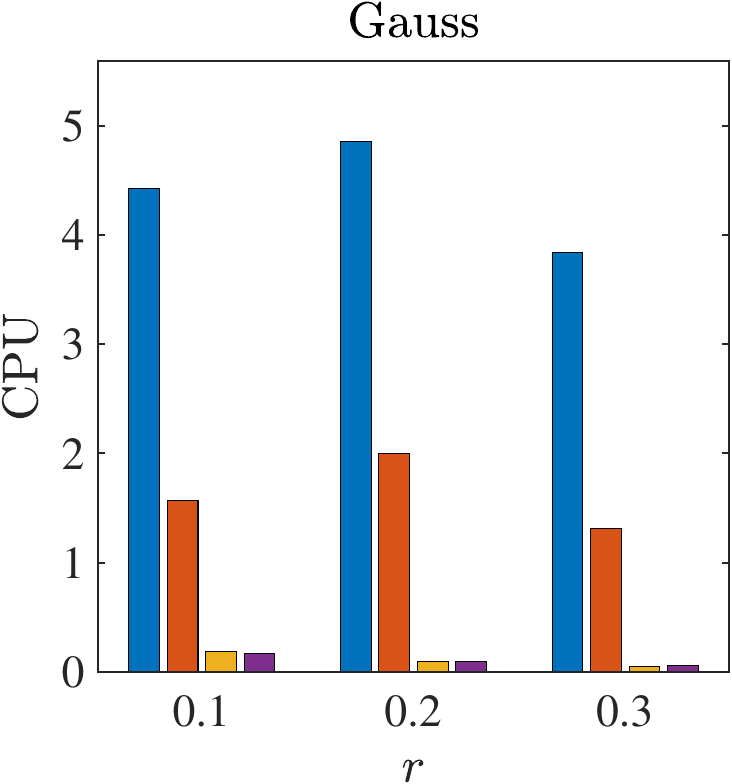} \\
%\bottomrule
\end{tabular}
	\caption{Mean CPU times (seconds) by $G$ for ODCT (left) and by $r$ for Gaussian matrices (right).}
\label{fig:performance_profile}
\end{figure}

%	\begin{figure}[t]
%		\centering
%		\caption*{\textbf{Comparison of Objective, Error and KKT versus Iteration ($\times 10^4$)}}%
%		\vspace{1mm}
%		\begin{subfigure}{1}
%			\includegraphics[width=0.1\linewidth]{APG_CN_figure/ODCT_Iteration_obj.eps}
%			\caption{Obj vs Iter $(\times 10^4)$}
%		\end{subfigure}
%		\begin{subfigure}{2}
%			\includegraphics[width=0.1\linewidth]{APG_CN_figure/ODCT_Iteration_obj.eps}
%			\caption{Obj vs Iter $(\times 10^4)$}
%		\end{subfigure}
%	\begin{subfigure}{3}
%		\includegraphics[width=0.1\linewidth]{APG_CN_figure/ODCT_Iteration_obj.eps}
%		\caption{Obj vs Iter $(\times 10^4)$}
%	\end{subfigure}
%		
%		\vspace{-1mm}
%		\caption{$r_k$ and $\text{rank}(X_k)$ for RA-PFBS-L$_1$/L$_2$, RA-PFBS-MCP, RA-PFBS-SCAD and RA-PFBS-Lp.}
%		\label{fig:odct_iterations}
%	\end{figure}
{%
%Figure~\ref{fig:emperical_CDF_kkt} presents the empirical cumulative distributions of the KKT residuals obtained by the four algorithms on 600 test instances generated using two types of sensing matrices, ODCT and Gauss. A point $(x,y)$ on a curve indicates that the corresponding algorithm achieves a KKT residual no greater than $x$ on $y\%$ of the test instances. As shown in the figure, both HPG-SCN and AHPG-SCN achieve KKT residuals below $10^{-12}$ on the vast majority of test instances for both ODCT and Gauss, with HPG-SCN generally attaining higher stationarity accuracy. In contrast, ZeroFPR and ManAcc typically require substantially larger thresholds to achieve a comparable proportion of instances. These results demonstrate that HPG-SCN and AHPG-SCN consistently attain substantially higher stationarity accuracy than ZeroFPR and ManAcc.

{Figure~\ref{fig:performance_profile} reports the mean CPU times for ODCT (left) and Gaussian (right) sensing matrices. HPG-SCN and AHPG-SCN require substantially less CPU time than ZeroFPR and ManAcc at every plotted value of $G$ and $r$. For ODCT, the CPU times of ZeroFPR and ManAcc increase with $G$, whereas those of HPG-SCN and AHPG-SCN remain relatively stable. In particular, when $G=25$, the proposed methods require approximately $2\%$ of the CPU time of ZeroFPR and $3$--$5\%$ of that of ManAcc. For Gaussian matrices, the proposed methods also maintain low CPU times across different values of $r$. For example, when $r=0.1$, their CPU times are approximately $4$--$5\%$ of that of ZeroFPR and $11$--$13\%$ of that of ManAcc.}

%Figure \ref{fig:performance_profile} presents the performance profile for all tested algorithms for the same
%450 test instances for ODCT (left) and Gauss (right), respectively. In each plot, the point where the `percentage' axis intersects with the curve indicates the performance ratio of the current solver compared to the best solver. From the figures, we observe that in approximately $80\%$ of the cases, {HPG-SCN} requires the least CPU time, whereas in about $20\%$ of the cases, {AHPG-SCN} achieves the smallest CPU time. Moreover, ManAcc and PG are approximately $50$ times slower than HPG-SCN in nearly $50\%$ of the test instances.
%%		Regarding the number of iterations, ManAcc and PG require roughly $10$ times and $1000$ times more iterations than AHPG-SCN in nearly half of the experiments.
%These observations clearly demonstrate the superior efficiency of HPG-SCN and AHPG-SCN.

%For this problem, HPG-SCN performs slightly better than AHPG-SCN. This may be attributed to the relatively tight estimates of the local Lipschitz constants of $\nabla h$ and $\nabla_{\mathcal{S}}^2 F$, which already yield effective proximal-gradient stepsizes and cubic regularization parameters. Consequently, the adaptive strategies in AHPG-SCN provide limited additional benefits while potentially introducing extra computational overhead.

%\subsubsection[L1/2-Regularized Least-Squares Regression]{$L_{1/2}$-Regularized Least-Squares Regression}
%\label{sec:lp_regression}

Next, we solve the $L_{1/2}^{1/2}$-regularized least-squares problem.
{\begin{equation*}
\min_{\mathbf{x}\in\mathbb{R}^{n}}
F(\mathbf{x})
:=\frac12\|A\mathbf{x}-\mathbf{b}\|_2^2
+\lambda\|\mathbf{x}\|_{\frac{1}{2}}^{\frac{1}{2}}
\end{equation*}}
{The objective function $F$ is continuous and semialgebraic and hence satisfies the KL property \cite[Theorem 14]{attouch2010proximal}.}
We test this model on four real-world datasets.
\begin{itemize}

\item \texttt{winequality-red}\footnote{\url{https://archive.ics.uci.edu/dataset/186/wine+quality}}: The red-wine subset of the Wine
Quality dataset, which predicts wine quality from physicochemical
measurements. It contains $m = 1599$ observations and $n_0=11$ input
features.

\item \texttt{pyrim}\footnote{\url{https://www.csie.ntu.edu.tw/~cjlin/libsvm/}}: A quantitative structure--activity
relationship regression dataset with $m = 74$ observations and
$n_0= 27$ input features.

\item \texttt{airfoil}\footnote{\url{https://archive.ics.uci.edu/dataset/291/airfoil+self+noise}}: The Airfoil Self-Noise dataset, which
predicts the scaled sound pressure level from aerodynamic
measurements. It contains $m = 1503$ observations and $n_0=5$ input
features.

\item \texttt{mpg}\footnote{\url{https://www.csie.ntu.edu.tw/~cjlin/libsvm/}}: An automobile fuel-consumption dataset
with $m = 392$ observations and $n_0 = 7$ input features.
\end{itemize}

{Because the original datasets have more observations than features, following \citep{huang2010predicting,wu2023regularized},
we expand the feature space using polynomial terms to obtain the expanded dataset $A\in\mathbb{R}^{m\times n}$
with $m<n$. Table~\ref{tab:lp_realdata_comparison} reports the expanded dimensions.}
%Let $Z\in\mathbb R^{m\times n_0}$ be the original feature matrix.
%We choose the smallest positive integer $r$ such that
%$n:=\binom{n_0+r}{r}>m$ and construct
%$A\in\mathbb R^{m\times n}$ by}
%\[
%A_{i,\boldsymbol{\nu}}
%=\prod_{j=1}^{n_0}Z_{ij}^{\nu_j},
%\quad
%\boldsymbol{\nu}\in\mathbb{N}^{n_0},
%\quad
%\sum_{j=1}^{n_0}\nu_j\le r,
%\quad
%i \in [m].
%\]
%The original and expanded dimensions are listed in Table~\ref{tab:lp_data_expansion}.
%
%\begin{table}[htbp]
%\centering
%\small
%\setlength{\tabcolsep}{6pt}
%\begin{tabular}{cccccc}
%\toprule
%Dataset
%& \multicolumn{2}{c}{Original}
%& Degree
%& \multicolumn{2}{c}{Expanded} \\
%\cmidrule(lr){2-3}
%\cmidrule(lr){5-6}
%& $m$ & $n_0$ & $r$ & $m$ & $n$ \\
%\midrule
%\texttt{winequality-red}  & 1599 & 11 & 5 & 1599 & 4368 \\
%\texttt{pyrim}            & 74   & 27 & 2 & 74   & 406  \\
%\texttt{airfoil}          & 1503 & 5  & 9 & 1503 & 2002 \\
%\texttt{mpg}              & 392  & 7  & 5 & 392  & 792  \\
%\bottomrule
%\end{tabular}
%\caption{Original and expanded dimensions of the datasets.}
%\label{tab:lp_data_expansion}
%\end{table}

{In this model, we set
$\lambda=10^{-4}\|A^\top\mathbf{b}\|_\infty.$
For HPG-SCN, we set $ M(\mathbf{x}^k) =1.1\,\frac{3\lambda}{8} \left(\min_{j\in\mathcal{S}^k}|x_j^k|\right)^{-5/2}$ (see Theorem \ref{LHbound:Lp} and it permits $\underline{\omega}=1.1>1$). For both HPG-SCN and AHPG-SCN, the threshold is given by
\[
\varepsilon_k={C_\varepsilon}/{\bigl[\log\bigl(e+\max\{s_k,10^{-8}/(k+1)\}^{-1}\bigr)\bigr]^{q_\varepsilon}},
\]
with $q_\varepsilon=3$ and $C_\varepsilon=0.3$. We also set $n_{\max}=n$.
The initial point is drawn from a standard Gaussian distribution
independently in each trial.}

\begin{table}[htbp]
	\centering
	\small
	\setlength{\tabcolsep}{6pt}
	\renewcommand{\arraystretch}{1.08}
	\caption{Mean results on four real-world datasets. An asterisk ($^{*}$) indicates that the method reached the $200$-second CPU time limit in at least one trial.}
	\label{tab:lp_realdata_comparison}
	\begin{tabular}{cccccc}
		\toprule
		Dataset
		& \multirow{2}{*}{Algorithm}
		& \multirow{2}{*}{Obj}
		& \multirow{2}{*}{NNZ}
		& \multirow{2}{*}{CPU}
		& \multirow{2}{*}{CPU Ratio} \\
		{$(m,n)$} & & & & & \\
		\midrule
		\multirow{5}{*}{\shortstack[c]{\texttt{winequality-red}\\{$(1599,4368)$}}}
		& ZeroFPR & 4.31e+02 & 339.30 & 1.96e+02$^*$ & 0.35\% \\
		& ManAcc & 5.35e+02 & 387.30 & 2.00e+02$^*$ & 0.35\% \\
		& HpgSRN & 3.12e+02 & 37.20 & 2.26e+01 & 3.06\% \\
		& HPG-SCN & 3.11e+02 & 38.90 & 6.90e-01 & 100.00\% \\
		& AHPG-SCN & 3.12e+02 & 30.60 & 8.26e-01 & 83.52\% \\
		\midrule
		\multirow{5}{*}{\shortstack[c]{\texttt{pyrim}\\{$(74,406)$}}}
		& ZeroFPR & 7.83e-02 & 24.70 & 6.80e+00 & 3.40\% \\
		& ManAcc & 8.34e-02 & 28.70 & 3.34e-01 & 69.24\% \\
		& HpgSRN & 1.85e-01 & 40.50 & 3.42e+01 & 0.68\% \\
		& HPG-SCN & 6.28e-02 & 29.90 & 3.49e-01 & 66.30\% \\
		& AHPG-SCN & 5.20e-02 & 30.30 & 2.31e-01 & 100.00\% \\
		\midrule
		\multirow{5}{*}{\shortstack[c]{\texttt{airfoil}\\{$(1503,2002)$}}}
		& ZeroFPR & 6.60e+03 & 38.00 & 3.94e+01 & 4.07\% \\
		& ManAcc & 7.12e+03 & 40.40 & 6.76e+00 & 23.75\% \\
		& HpgSRN & 6.89e+03 & 36.50 & 2.41e+00 & 66.45\% \\
		& HPG-SCN & 6.67e+03 & 40.00 & 1.60e+00 & 100.00\% \\
		& AHPG-SCN & 6.66e+03 & 40.30 & 1.88e+00 & 85.50\% \\
		\midrule
		\multirow{5}{*}{\shortstack[c]{\texttt{mpg}\\{$(392,792)$}}}
		& ZeroFPR & 7.80e+02 & 72.90 & 9.90e+00 & 1.87\% \\
		& ManAcc & 8.26e+02 & 67.70 & 7.48e-01 & 24.81\% \\
		& HpgSRN & 8.11e+02 & 71.80 & 2.96e+00 & 6.26\% \\
		& HPG-SCN & 8.13e+02 & 67.70 & 1.85e-01 & 100.00\% \\
		& AHPG-SCN & 7.96e+02 & 71.80 & 1.98e-01 & 93.56\% \\
		\bottomrule
	\end{tabular}
\end{table}

Table~\ref{tab:lp_realdata_comparison} {reports the mean objective value (Obj), number of nonzero entries (NNZ), CPU time, and CPU ratio. We define the CPU ratio by}
\begin{equation*}
	\mathrm{CPU\ Ratio}
	:=
	\frac{\text{minimum mean CPU time across methods}}
	{\text{mean CPU time for the method}}
	\times 100\%.
\end{equation*}

{As shown in the table, HPG-SCN achieves the lowest mean CPU time on \texttt{winequality-red}, \texttt{airfoil}, and \texttt{mpg}, while AHPG-SCN is the fastest on \texttt{pyrim}. In particular, on \texttt{winequality-red}, HPG-SCN and AHPG-SCN require only $0.69$ and $0.83$ seconds, respectively, corresponding to approximately $3.1\%$ and $3.7\%$ of the CPU time of HpgSRN, and less than $0.5\%$ of that of ZeroFPR and ManAcc. Both proposed methods also attain competitive objective values across the four datasets.}

\subsubsection{Logistic Regression}\label{sec:logistic}
For these experiments, we use the $L_1-\mu L_2$-regularized logistic regression model in \eqref{eq:intro_sparse_logistic}.
%\citep{qin20191}:
\begin{equation*}
\min_{(\z,u)\in\mathbb{R}^{n}\times\mathbb{R}} \frac{1}{m}\sum_{i = 1}^m\log(1 + e^{-b_i(\mathbf{a}^{\top}_i\z+u)}) + \lambda(\|\z\|_1 - \mu \|\z\|_2),\label{prob:L1mL2_logistic}
\end{equation*}
{The objective function is both continuous and definable on $\mathbb{R}^{n+1}$ and hence satisfies the KL property \cite[Theorem 14]{attouch2010proximal}.}

{We first test on synthetic datasets.
For the synthetic-data experiments, we set $\mu\in\{0.5,1\}$ and $\lambda=3\times10^{-3}$ in~\eqref{eq:intro_sparse_logistic}.}
Let
$
\widetilde A=[A,\mathbf e]
\in\mathbb R^{m\times(n+1)},
$
where $\mathbf e\in\mathbb R^m$ denotes the vector with all entries equal to one.
The logistic loss gradient is Lipschitz continuous with constant
$L=\|\widetilde A\|^2/(4m)$.
%For all these tested algorithms, we set $h(\x) = \phi(\x;A,\mathbf{b}),\ g(\x)=\lambda(\|\z\|_1 - \mu \|\z\|)$.
{For HPG-SCN, we take $M(\x^k)=1.1B_k$ (Theorem~\ref{LHbound:L1mL2} permits $\underline{\omega}=1.1>1$), where
$\mathcal S_z^k=\operatorname{supp}(\mathbf{z}^k)$ and
\[
\begin{aligned}
B_k&=
\begin{cases}
	M_{\mathrm{log}}+2\lambda\mu/(\sqrt{3}\|\z^k\|_2^2),&\z^k\ne0,\\
	{M_{\mathrm{log}}\;(=1/(6\sqrt3))},&\z^k=0\quad(u^k\ne0)
\end{cases},\\
M_{\mathrm{log}}&=
\frac{1}{6\sqrt3m}\sum_{i=1}^m
(\|\mathbf a_{i,{\mathcal S}_{z}^k}\|_2^2+1)^{3/2}.
\end{aligned}
\]
For both HPG-SCN and AHPG-SCN, we choose the threshold
$\varepsilon_k=C_{\varepsilon}/\log(e+\max\{s_k,10^{-8}/(k+1)\}^{-1})^{q_{\varepsilon}}$, with $C_{\varepsilon} = 3$ and $q_{\varepsilon} = 4$. For AHPG-SCN, we set $n_{\max}=n$, the number of features.}

{\begin{table}[htbp]
	\centering
	\small
	\setlength{\tabcolsep}{6pt}
	\renewcommand{\arraystretch}{1.12}
	\caption{Average CPU time and CPU Ratio for
		ZeroFPR, ManAcc, HPG-SCN, and AHPG-SCN
		with different $(\mu,n)$.}
	\label{tab:mean_logistic}
	\begin{tabular}{llccccc}
		\toprule
		\multicolumn{2}{c}{$\mu=0.5$}
		& $n=400$ & $n=800$ & $n=1200$
		& $n=1600$ & $n=2000$ \\
		\midrule
		\multirow{2}{*}{ZeroFPR}
		& CPU & 6.34e-01 & 3.37e+00 & 6.62e+00
		& 1.42e+01 & 2.11e+01 \\
		& CPU Ratio & 2.44\% & 0.73\% & 0.90\%
		& 0.60\% & 0.58\% \\
		\midrule
		\multirow{2}{*}{ManAcc}
		& CPU & 3.56e-01 & 1.07e+00 & 5.25e+00
		& 6.62e+00 & 1.50e+01 \\
		& CPU Ratio & 4.35\% & 2.31\% & 1.13\%
		& 1.29\% & 0.82\% \\
		\midrule
		\multirow{2}{*}{HPG-SCN}
		& CPU & 3.09e-02 & 4.66e-02 & 1.03e-01
		& 1.51e-01 & 2.01e-01 \\
		& CPU Ratio & 50.16\% & 53.00\% & 57.57\%
		& 56.49\% & 61.19\% \\
		\midrule
		\multirow{2}{*}{AHPG-SCN}
		& CPU & 1.55e-02 & 2.47e-02 & 5.93e-02
		& 8.53e-02 & 1.23e-01 \\
		& CPU Ratio & 100.00\% & 100.00\% & 100.00\%
		& 100.00\% & 100.00\% \\
		\midrule
		\multicolumn{2}{c}{$\mu=1$}
		& $n=400$ & $n=800$ & $n=1200$
		& $n=1600$ & $n=2000$ \\
		\midrule
		\multirow{2}{*}{ZeroFPR}
		& CPU & 6.01e-01 & 3.23e+00 & 6.61e+00
		& 1.47e+01 & 2.42e+01 \\
		& CPU Ratio & 2.36\% & 0.71\% & 0.95\%
		& 0.61\% & 0.61\% \\
		\midrule
		\multirow{2}{*}{ManAcc}
		& CPU & 7.36e-01 & 1.59e+00 & 8.57e+00
		& 1.68e+01 & 2.21e+01 \\
		& CPU Ratio & 1.93\% & 1.44\% & 0.74\%
		& 0.54\% & 0.67\% \\
		\midrule
		\multirow{2}{*}{HPG-SCN}
		& CPU & 2.82e-02 & 4.32e-02 & 1.00e-01
		& 1.56e-01 & 2.23e-01 \\
		& CPU Ratio & 50.35\% & 53.01\% & 63.00\%
		& 57.95\% & 66.37\% \\
		\midrule
		\multirow{2}{*}{AHPG-SCN}
		& CPU & 1.42e-02 & 2.29e-02 & 6.30e-02
		& 9.04e-02 & 1.48e-01 \\
		& CPU Ratio & 100.00\% & 100.00\% & 100.00\%
		& 100.00\% & 100.00\% \\
		\bottomrule
	\end{tabular}
\end{table}
}
{Table~\ref{tab:mean_logistic} reports the mean CPU times and CPU ratios of the four methods on synthetic datasets with $\mu\in\{0.5,1\}$ and $n\in\{400,800,1200,1600,2000\}$. AHPG-SCN achieves the lowest mean CPU time in all tested configurations. In particular, for $\mu=1$ and $n=2000$, AHPG-SCN requires only $0.148$ seconds, less than $1\%$ of the CPU time of ZeroFPR and ManAcc, while HPG-SCN takes $0.223$ seconds. Across all configurations, the CPU ratios of HPG-SCN range from $50.16\%$ to $66.37\%$, whereas those of ZeroFPR and ManAcc remain below $5\%$.}

We next consider two real-world datasets, each with a training set $(A,\mathbf{b})\in\mathbb{R}^{m\times n}\times\mathbb{R}^{m}$ and {an evaluation set} $(A_{\textrm{test}},\mathbf{b}_{\textrm{test}})\in\mathbb{R}^{m_{\textrm{te}}\times n}\times\mathbb{R}^{m_{\textrm{te}}}$.
\begin{itemize}

\item \texttt{duke}\footnote{\url{https://www.csie.ntu.edu.tw/~cjlin/libsvm/}}:
A binary breast-cancer classification dataset with 44 observations in total and $n=7129$ gene-expression features. {{We use all $44$ observations for both fitting and evaluation ($m_{\text{te}}=m=44$), so the reported accuracy is in-sample.}}

\item \texttt{news20}\footnote{\url{https://scikit-learn.org/stable/datasets/real_world.html}}:
The original dataset contains documents from 20 newsgroup categories. For binary classification, we retain only the baseball and hockey categories and randomly split them into training and test sets in an 8:2 ratio. Each sample is normalized to have unit 2-norm. The resulting datasets contain $m=1595$ training samples and $m_{\textrm{te}}=398$ test samples, with $n=23416$ features.
\end{itemize}

%To apply PG, we take $h(\x,u) = \phi(\x,u;A_{\textrm{train}},\mathbf{b}_{\textrm{train}}), g(\x)=\lambda(\|\x\|_1 - \mu \|\x\|)$, where $\lambda = 5\times10^{-2}$. We set $\alpha = 1$.
%
%To apply ManAcc, we take $f(\x,u) = \phi(\x,u;A_{\textrm{train}},\mathbf{b}_{\textrm{train}}), g(\x)=\lambda(\|\x\|_1 - \mu \|\x\|)$, where $5\times10^{-2}$. We set $\gamma = 1, m_1 = 10^{-4}, \theta = 0.5, \mathcal{M}_k = \{\z: \supp(\z) = \supp(\x^k)\}$. Moreover, to enhance the convergence, we also use $H^k + \eta I$ as the modified Hessian matrix at each iteration and we set $\eta = 10^{-3}$.
%
%To apply AHPG-SCN, we take $h(\x,u) = \phi(\x,u;A_{\textrm{train}},\mathbf{b}_{\textrm{train}}), g(\x)=\lambda(\|\x\|_1 - \mu \|\x\|)$, where $\lambda = 5\times10^{-2}$. We set $\alpha = 1, c= 10^{-3},\eta = 10^{-10}, q= 0.95, \eta_1 = 0.1,\eta_2 = 0.9, \gamma_1 = 0.5,\gamma_2 = 2,\sigma_{\min} = 2c,\sigma_{\max} = 100\sigma_{\min}$ and $\varepsilon_k = \tau_1\|(\x^{k} - \y^{k-1})/\alpha\|^{\tau_2}$ with $\tau_1 = 0.85,\tau_2 = 0.99$.

{We draw the initial point from a standard Gaussian distribution.}
%{For termination criterion (\ref{stoppping}), we set $\texttt{RTol} = 10^{-8},\texttt{CTol} = 200$ and $\texttt{MaxIt} = 100000$.}
{We compute the classification accuracy (TEAcc) as follows}
$$\textrm{TEAcc} = \displaystyle{\frac{\sharp\{i : \textrm{sign}(\mathbf{a}_{\textrm{test},i}^{\top}\z + u) = b_{\textrm{test},i}\}}{m_{\text{te}}}}.
$$

%In model (\ref{prob:L1mL2_logistic}), we set $(\mu,\lambda)=(0.5,5e-2)$ and $(\mu,\lambda)=(1.0,5e-2)$  for \texttt{a9a};
%and  $(\mu,\lambda)=(0.5,2e-2)$ and $(\mu,\lambda)=(1.0,2e-2)$  for \texttt{MNIST} by hand-tuning.
{In model (\ref{eq:intro_sparse_logistic}), we set $(\mu,\lambda,C_{\varepsilon},q_{\varepsilon})=(0.5,1.5\times10^{-2},8,8)$ and $(1.0,1.5\times10^{-2},8,8)$ for the \texttt{duke} dataset, and
$(\mu,\lambda,C_{\varepsilon},q_{\varepsilon})=(0.5,1\times10^{-3},30,4)$ and $(1.0,1\times10^{-3},30,4)$ for the \texttt{news20} dataset.
For AHPG-SCN, we set $n_{\max}=500$.
We retain the other algorithmic settings.}

\begin{figure}[htbp]
\centering
\begin{tabular}{@{}c@{\hspace{.025\linewidth}}c@{}}
%\toprule
%ODCT & Gauss \\
%\midrule
%\includegraphics[scale=.37]{APG_CN_figure/KKT_a9a_mu055.eps}
%& \includegraphics[scale=.37]{APG_CN_figure/KKT_MNIST_mu055.eps} \\
%\includegraphics[scale=.32]{APG_CN_figure/KKT_a9a_mu1.eps} &
%\includegraphics[scale=.32]{APG_CN_figure/KKT_MNIST_mu1.eps}
\includegraphics[width=.48\linewidth]{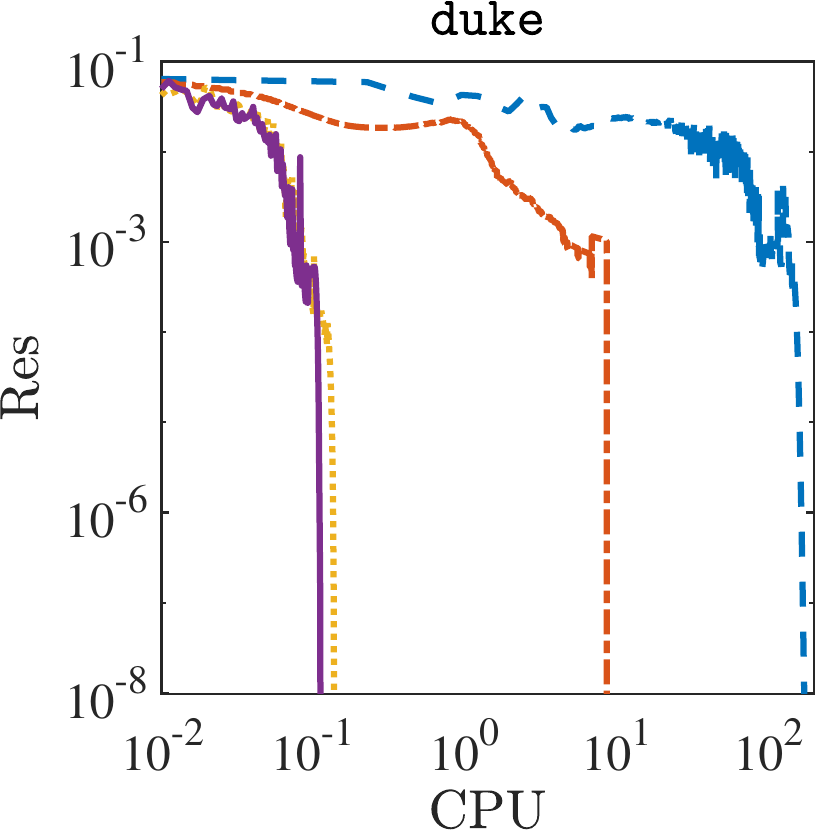}
& \includegraphics[width=.48\linewidth]{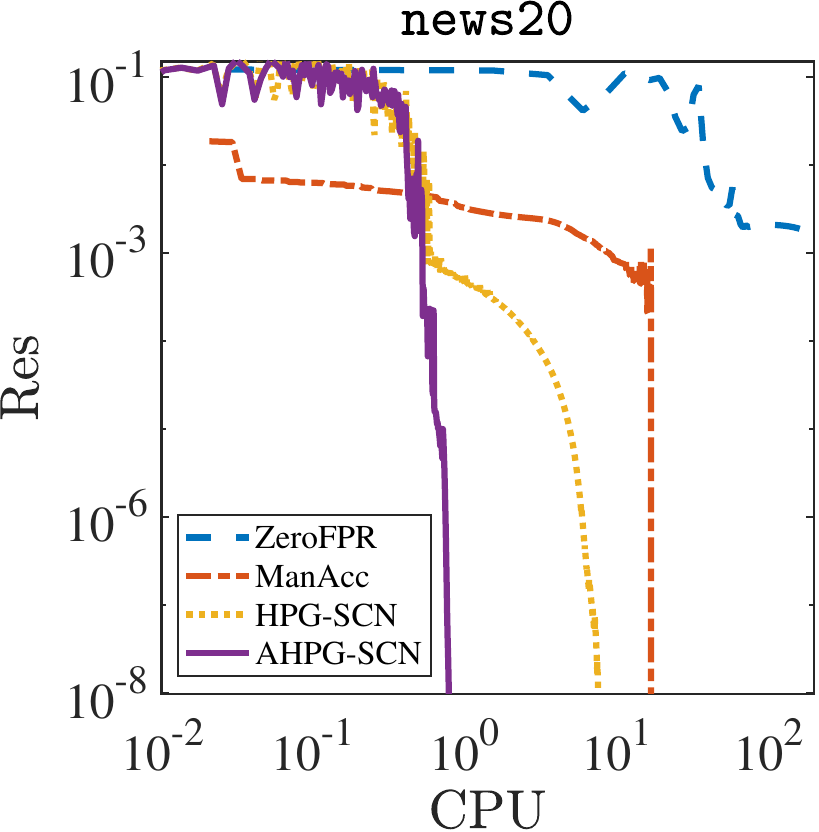}
%\bottomrule
\end{tabular}

\caption{Averaged stationarity residual versus cumulative CPU time for duke (left) and news20 (right).}
\label{fig:KKT_realdata}
\end{figure}
{Figure~\ref{fig:KKT_realdata} plots the mean stationarity residual against cumulative CPU time for $\mu=1$. On \texttt{duke}, HPG-SCN and AHPG-SCN both reduce the residual to the $10^{-8}$ range in about $0.1$ seconds. On \texttt{news20}, the mean residual of AHPG-SCN reaches the $10^{-8}$ range in under one second, followed by HPG-SCN. ManAcc takes longer, and ZeroFPR reaches the time limit in at least one trial.}
%Table~\ref{tab:real_data_binary} reports the performance of all methods in terms of {Obj}, {NNZ}, {TEAcc}, {KKT}, {CPU}, and {CPU Ratio} for $\mu \in \{0.5,1\}$.
%All methods produce nearly identical values of {Obj}, {NNZ}, and {TEAcc}, whereas AHPG-SCN achieves significantly smaller {KKT} residuals and markedly lower CPU time.
%{Among all the datasets, the CPU time of AHPG-SCN is approximately $5\%$ and $3\%$ of that required by PG and ManAcc, respectively.}
%Moreover, AHPG-SCN consistently outperforms HPG-SCN in terms of computational efficiency.
%This advantage is mainly due to the conservative estimate of $M(\x^k)$ in HPG-SCN, which leads to smaller steps in the cubic Newton update.
%Overall, these results demonstrate that HPG-SCN and AHPG-SCN attain high-accuracy solutions at a substantially lower computational cost than PG and ManAcc, highlighting the effectiveness of the proposed methods on both synthetic and real datasets.

{Table~\ref{tab:real_data_binary} reports the TEAcc, CPU times, and CPU ratios for the four methods with $\mu\in\{0.5,1\}$. All methods achieve comparable in-sample accuracies on \texttt{duke}. In particular, for $\mu=1$, AHPG-SCN achieves an accuracy of $0.99$ in only $0.079$ seconds, requiring less than $0.1\%$ of the CPU time of ZeroFPR. On \texttt{news20}, HPG-SCN, AHPG-SCN, and ManAcc achieve an accuracy of $0.89$, compared with $0.88$ for ZeroFPR. For $\mu=1$, AHPG-SCN requires only $0.554$ seconds, approximately $15\%$ and $4\%$ of the CPU times of HPG-SCN and ManAcc, respectively, and ZeroFPR reaches the $200$-second time limit in at least one trial.}

\begin{table}[htbp]
	\centering
	\small
	\setlength{\tabcolsep}{8pt}
	\renewcommand{\arraystretch}{1.12}
	\caption{{Mean classification accuracy, CPU time, and CPU ratio} for comparing algorithms. An asterisk ($^*$) indicates
		that the CPU-time limit of $200$ seconds was
		reached in at least one trial.}
	\label{tab:real_data_binary}
	\begin{tabular}{llcccc}
		\toprule
		\multirow{2}{*}{Algorithm}
		& \multirow{2}{*}{Measurement}
		& \multicolumn{2}{c}{\texttt{duke}}
		& \multicolumn{2}{c}{\texttt{news20}} \\
		\cmidrule(lr){3-4}
		\cmidrule(lr){5-6}
		& & $\mu=0.5$ & $\mu=1$
		& $\mu=0.5$ & $\mu=1$ \\
		\midrule
		\multirow{3}{*}{ZeroFPR}
		& TEAcc & 0.95 & 0.98 & 0.88 & 0.88 \\
		& CPU & 1.31e+02 & 1.25e+02
		& 2.00e+02$^*$ & 2.00e+02$^*$ \\
		& CPU Ratio & 0.07\% & 0.06\%
		& 0.34\% & 0.28\% \\
		\midrule
		\multirow{3}{*}{ManAcc}
		& TEAcc & 0.95 & 0.98 & 0.89 & 0.89 \\
		& CPU & 3.69e+00 & 5.84e+00
		& 1.40e+01 & 1.31e+01 \\
		& CPU Ratio & 2.34\% & 1.36\%
		& 4.84\% & 4.23\% \\
		\midrule
		\multirow{3}{*}{HPG-SCN}
		& TEAcc & 0.95 & 0.98 & 0.89 & 0.89 \\
		& CPU & 1.05e-01 & 1.08e-01
		& 6.80e+00 & 3.67e+00 \\
		& CPU Ratio & 82.29\% & 73.43\%
		& 9.97\% & 15.10\% \\
		\midrule
		\multirow{3}{*}{AHPG-SCN}
		& TEAcc & 0.95 & 0.99 & 0.89 & 0.89 \\
		& CPU & 8.64e-02 & 7.93e-02
		& 6.78e-01 & 5.54e-01 \\
		& CPU Ratio & 100.00\% & 100.00\%
		& 100.00\% & 100.00\% \\
		\bottomrule
	\end{tabular}
\end{table}

{Overall, these numerical results demonstrate the computational efficiency of HPG-SCN and AHPG-SCN across different sparse optimization problems.}

%\section{Conclusions}\label{sec:conclusion}
%
%
%We have developed a hybrid proximal-gradient and subspace cubic-Newton
%method for sparse composite optimization. Sufficient decrease and the proximal
%residual estimate give convergence to a first-order stationary point under the
%boundedness and KL assumptions. A threshold that decreases more slowly than
%the iterate error identifies the support of the limit by the selected coordinate
%set, without a nondegeneracy condition.
%
%Under the additional local restricted-Hessian and regularization assumptions,
%the method eventually accepts full cubic steps and admits a local outer-iteration
%bound for restricted approximate second-order stationarity. This guarantee is
%distinct from a full-space second-order condition for the nonsmooth objective,
%and the analysis does not bound the identification phase uniformly in the
%requested tolerances. The experiments illustrate the potential reduction in
%second-order cost when the selected subspace is small. Extending the framework
%to other second-order models requires verifying their descent, residual, and
%curvature conditions.
\section{Conclusions}\label{sec:conclusion}
HPG-SCN converges to a first-order stationary point under the
boundedness and KL assumptions. When a positive power-type KL exponent exists,
both the power and logarithmic thresholds identify the support of the limit
in finitely many iterations. The example in Section~\ref{sec:degenerate-example} illustrates this result in a degenerate setting where the proximal-gradient iterates themselves fail to identify the limiting support.
Under continuity of the restricted Hessian, the line search eventually
accepts full cubic steps. With bounded cubic regularization and
tolerance-dependent damping, we bound the number of outer iterations
from that point until the restricted gradient and curvature tolerances
are met. In the sparse recovery and logistic regression experiments, HPG-SCN and AHPG-SCN achieve lowest mean CPU time.

\vspace{0.2cm}

\noindent {\bf Data Availability} The datasets generated during and/or analyzed during the current study are available from
the corresponding author upon reasonable request.
\vspace{0.4cm}

\noindent {\Large \bf Declarations}
\vspace{0.2cm}

\noindent{\bf Conflict of interest} The authors declare that there is no conflict of interest related to this work.

\appendix

\section{Bounds for the Local Restricted Hessian Lipschitz Modulus}\label{app:LHLC}

{\begin{theorem}\label{LHbound:L1oL2}
	Let $F(\x)=\lambda\|\x\|_1/\|\x\|_2+
	\tfrac12\|A\x-\mathbf b\|_2^2$, with $\lambda>0$.
	For $\x\ne0$ and $\mathcal S=\supp(\x)$, the local constant relative
	to the fixed-support manifold satisfies
	\[
	\mathrm{lip}\,\nabla^2_{\mathcal S}F(\x)
	\le\frac{36\lambda\sqrt{\|\x\|_0}}{\|\x\|_2^3}.
	\]
\end{theorem}
\begin{proof}
	Let $s=|\mathcal S|$ and consider the restriction of $F$
	to the active coordinates, using the same notation for
	the restricted vector.
	In a sufficiently small neighborhood of $\x$, the signs
	of the nonzero components remain unchanged.
	Set $\mathbf a=\sign(\x)$, $r=\|\x\|_2$, and
	$u=\mathbf a^\top\x$.
	Only the regularizer contributes to Hessian variation because the
	quadratic loss has a constant Hessian.
	For any direction $\mathbf h\in\mathbb R^s$, let
	$\alpha=\mathbf a^\top\mathbf h$ and
	$\beta=\x^\top\mathbf h$.
	Differentiating the Hessian of $\lambda u/r$ yields
	\[
	\begin{aligned}
		D(\nabla^2 F)(\x)[\mathbf h]
		={}&\frac{3\lambda\beta}{r^5}
		(\mathbf a\x^\top+\x\mathbf a^\top+uI)
		-\frac{\lambda}{r^3}
		(\mathbf a\mathbf h^\top+\mathbf h\mathbf a^\top+\alpha I)\\
		&+\frac{3\lambda\alpha}{r^5}\x\x^\top
		+\frac{3\lambda u}{r^5}
		(\x\mathbf h^\top+\mathbf h\x^\top)
		-\frac{15\lambda u\beta}{r^7}\x\x^\top.
	\end{aligned}
	\]
	Using $\|\mathbf a\|_2=\sqrt s$, $|u|\le\sqrt s\,r$,
	$|\alpha|\le\sqrt s\|\mathbf h\|_2$, and
	$|\beta|\le r\|\mathbf h\|_2$, we obtain
	\[
	\|D(\nabla^2 F)(\x)[\mathbf h]\|
	\le(9+3+3+6+15)\frac{\lambda\sqrt s}{r^3}\|\mathbf h\|_2.
	\]
	Since the restricted function is $C^3$,
	its local Hessian Lipschitz constant equals
	the operator norm of $D(\nabla^2 F)(\x)$.
	The desired bound follows.
\end{proof}}

{\begin{theorem}\label{LHbound:Lp}
	Let
	\[
	F(\x)=\lambda\|\x\|_p^p
	+\frac12\|A\x-\mathbf b\|_2^2,
	\]
	where $\lambda>0$ and $0<p<1$.
	For $\x\ne0$ and $\mathcal S=\supp(\x)$, the local
	constant relative to the fixed-support manifold satisfies
	\[
	\mathrm{lip}\,\nabla^2_{\mathcal S}F(\x)
	\le
	\frac{\lambda p(1-p)(2-p)}
	{\left(\min_{i\in\mathcal S}|x_i|\right)^{3-p}}.
	\]
\end{theorem}
\begin{proof}
	Let $s=|\mathcal S|$ and consider the restriction of
	$F$ to the active coordinates, using the same notation
	for the restricted vector.
	Since the signs of the nonzero components remain
	unchanged in a sufficiently small neighborhood of
	$\x$, the restricted function is $C^3$, with Hessian
	\[
	\nabla^2 F(\x)
	=A_{\mathcal S}^{\top}A_{\mathcal S}
	+\lambda p(p-1)
	\operatorname{Diag}\left(|x_i|^{p-2}\right)_{i\in\mathcal S}.
	\]
	For any direction $\mathbf h\in\mathbb R^s$,
	the directional derivative of the Hessian is
	\[
	D(\nabla^2 F)(\x)[\mathbf h]
	=
	\lambda p(p-1)(p-2)
	\operatorname{Diag}\left(
	\operatorname{sign}(x_i)
	|x_i|^{p-3}h_i
	\right)_{i\in\mathcal S}.
	\]
	Using the spectral norm of a diagonal matrix,
	we obtain
	\[
	\begin{aligned}
		\|D(\nabla^2 F)(\x)[\mathbf h]\|
		=
		\lambda p(1-p)(2-p)
		\max_{i\in\mathcal S}
		\frac{|h_i|}{|x_i|^{3-p}}\le
		\frac{\lambda p(1-p)(2-p)}
		{\left(\min_{i\in\mathcal S}|x_i|\right)^{3-p}}
		\|\mathbf h\|_2.
	\end{aligned}
	\]
	Since the restricted function is $C^3$, its local
	Hessian Lipschitz constant equals the operator norm
	of $D(\nabla^2 F)(\x)$.
	The desired bound follows.
\end{proof}}

\begin{theorem}\label{LHbound:L1mL2}
Let
\[
F(\z,u)=\frac1m\sum_{i=1}^m
\log(1+e^{-b_i(\mathbf a_i^\top\z+u)})
+\lambda(\|\z\|_1-\mu\|\z\|_2),
\]
where $b_i\in\{\pm1\}$ and $\lambda,\mu>0$.
For $\x=(\z,u)\ne0$, put
$\mathcal S=\supp(\x)$ and {$\mathcal{S}_z=\supp(\z)$}, and define
{\[
M_{\mathrm{log}}=
\frac{1}{6\sqrt3m}\sum_{i=1}^m
(\|\mathbf a_{i,\mathcal{S}_z}\|_2^2+1)^{3/2}.
\]}
Then
\[
\mathrm{lip}\,\nabla^2_{\mathcal S}F(\x)
\le
\begin{cases}
M_{\mathrm{log}}+2\lambda\mu/(\sqrt{3}\|\z\|_2^2),&\z\ne0,\\
{M_{\mathrm{log}}\;(=1/(6\sqrt3))},&\z=0\quad(u\ne0).
\end{cases}
\]
\end{theorem}
\begin{proof}
Let $\mathbf v_i=(\mathbf a_i,1)_{\mathcal S}$ be the augmented
feature vector restricted to the selected coordinates. Thus the
intercept entry is included only when $u\ne0$, and
$\|\mathbf v_i\|_2^2\le\|\mathbf a_{i,\mathcal{S}_z}\|_2^2+1$.
Put $\theta_i=b_i(\mathbf a_i^\top\z+u)$,
$\sigma(t)=1/(1+e^{-t})$ and $w(t)=\sigma(t)(1-\sigma(t))$.
For a restricted direction $\mathbf h$,
\[
D\bigl(\nabla^2_{\mathcal S}\log(1+e^{-\theta_i})\bigr)[\mathbf h]
=w'(\theta_i)b_i(\mathbf v_i^\top\mathbf h)
\mathbf v_i\mathbf v_i^\top.
\]
Maximizing $|p(1-p)(1-2p)|$ over $p\in[0,1]$ gives
$\sup_t|w'(t)|=1/(6\sqrt3)$. Hence the logistic part has Hessian
derivative norm at most $M_{\mathrm{log}}\|\mathbf h\|_2$.

If $\z\ne0$, let {$r=\|\z_{\mathcal{S}_z}\|_2$},
and let $\mathbf{p}\neq\mathbf{0}$.
Writing {$\beta=\z_{\mathcal{S}_z}^\top\mathbf p$},
direct differentiation gives
{\[
\begin{aligned}
	D\bigl(\nabla^2\|\z_{\mathcal{S}_z}\|_2\bigr)[\mathbf p]
	={}&-\frac{\beta}{r^3}I
	+\frac{3\beta}{r^5}
	\z_{\mathcal{S}_z}\z_{\mathcal{S}_z}^\top-\frac{1}{r^3}
	(\z_{\mathcal{S}_z}\mathbf p^\top
	+\mathbf p\z_{\mathcal{S}_z}^\top) = -\frac{\|\mathbf p\|_2}{r^2}B(t),
\end{aligned}
\]}
where
{\[
\mathbf e=\frac{\z_{\mathcal{S}_z}}{r},
\quad
\mathbf v=\frac{\mathbf p}{\|\mathbf p\|_2},
\quad
t=\mathbf e^\top\mathbf v,
\quad
B(t)=tI-3t\mathbf e\mathbf e^\top
+\mathbf e\mathbf v^\top+\mathbf v\mathbf e^\top.
\]}
We estimate the spectral norm from the eigenvalues of $B(t)$.
When $|t|<1$, let
\[
\mathbf q=\frac{\mathbf v-t\mathbf e}{\sqrt{1-t^2}}.
\]
Then $\mathbf e$ and $\mathbf q$ are orthonormal, and
$\mathbf v=t\mathbf e+\sqrt{1-t^2}\mathbf q$.
Substituting this decomposition into $B(t)$ gives
\[
B(t)=tI-t\mathbf e\mathbf e^\top
+\sqrt{1-t^2}
(\mathbf e\mathbf q^\top+\mathbf q\mathbf e^\top).
\]
Since $\sqrt{4-3t^2}\ge |t|$, we obtain
\[
\|B(t)\|
=\frac{|t|+\sqrt{4-3t^2}}{2}
\leq\frac{2}{\sqrt3},
\]
with equality when $|t| = 1/\sqrt{3}$.
The cases $|t|=1$ follow directly, and the same bound holds in dimension one.
Consequently,
{\[
\left\|
D\bigl(\nabla^2\|\z_{\mathcal{S}_z}\|_2\bigr)
[\mathbf p]
\right\|
\le
\frac{2}{\sqrt3\,r^2}\|\mathbf p\|_2.
\]}
Since the {$L_1$} term is locally linear on a
fixed-support neighborhood, it contributes no
Hessian variation. Combining this estimate with
the bound for the logistic loss yields
\[
\operatorname{lip}_{\mathcal M_x}
\nabla_{\mathcal S}^2F(x)
\le
M_{\log}
+\frac{2\lambda\mu}{\sqrt{3}\|\z\|_2^2}.
\]

If $\z=0$ and $u\ne0$, then
$\mathcal S=\{n+1\}$. The regularizer vanishes
on this one-dimensional slice, so only the
logistic loss contributes to the Hessian
variation. Hence,
\[
\operatorname{lip}_{\mathcal M_x}
\nabla_{\mathcal S}^2F(x)
\le M_{\log}=\frac{1}{6\sqrt{3}}.
\]

%The $\ell_1$ term has zero
%restricted Hessian on a neighborhood with unchanged predictor signs.
%Adding the two bounds gives the first case.
%If $\z=0$ and $u\ne0$, then $\mathcal S=\{n+1\}$; the regularizer
%is identically zero on this one-dimensional slice. Only the logistic
%part contributes, which proves the second case.
\end{proof}

{\begin{proposition}\label{prop:toy_example}
Consider the function
\begin{equation*}
	F(\mathbf{x})
	:=
	\sum_{i=1}^n
	\left(\log(1+(x_i-b_i)^4)-x_i\right)
	+\|\mathbf{x}\|_1,
\end{equation*}
{Here $\mathbf b\in\{0,1\}^n$ has both zero and nonzero entries. Set
$g(\mathbf x)=\|\mathbf x\|_1$ and $h(\mathbf x)=F(\mathbf x)-g(\mathbf x)$.}
Then the following hold:
\begin{itemize}
	 \item[(i)] The point $\overline{\mathbf x}=\mathbf b$ is both
	the unique critical point and the unique global minimizer of $F$.
    The function $F$ is coercive, and the nondegeneracy condition fails at $\overline{\mathbf x}$.
	\item[(ii)] The gradient of $h$ is globally Lipschitz continuous with constant $L<4$.
	\item[(iii)] For any $\mathbf{0}\neq\mathbf{x}\in\mathbb{R}^n$ and $\emptyset\neq\mathcal{S}\subseteq\operatorname{supp}(\mathbf{x})$, it holds that
	$\operatorname{lip}_{\mathcal{M}_{\mathbf{x}}}\nabla^2_{\mathcal{S}}F(\mathbf{x})<9$.
	\item[(iv)] $F$ satisfies the KL property at $\overline{\mathbf{x}}$ with exact exponent $\theta=3/4$.
\end{itemize}
\end{proposition}
}

\begin{proof}
(i) Let $\phi(t)=\log(1+t^4)$. First note that $F(\mathbf{x})\geq 0$ and equality holds if and only if $\mathbf{x}=\mathbf{b}$.
Thus, $\overline{\mathbf{x}}=\mathbf{b}$ is the unique global minimizer, with $F(\overline{\mathbf{x}})=0$.
Coercivity follows from
$F(\mathbf x)\geq\log(1+\|\mathbf x-\mathbf b\|_\infty^4)$.
At $\overline{\mathbf x}=\mathbf b$, every component of
$-\nabla h(\overline{\mathbf x})$ equals $1$. For each $i$ with $b_i=0$,
the $i$th coordinate of $\operatorname{ri}\partial g(\overline{\mathbf x})$
ranges over $(-1,1)$. Thus the nondegeneracy condition fails.
To show that $\overline{\mathbf{x}}$ is also the unique critical point, observe that
\begin{equation*}
\partial_iF(\mathbf{x})=
\begin{cases}
	\{\phi'(x_i-b_i)\}, & x_i>0,\\
	[\phi'(-b_i)-2,\phi'(-b_i)], & x_i=0,\\
	\{\phi'(x_i-b_i)-2\}, & x_i<0.
\end{cases}
\end{equation*}
Since
\begin{equation*}
\phi'(t)=\frac{4t^3}{1+t^4},
\end{equation*}
we have $\phi'(t)=0$ if and only if $t=0$.
If $x_i>0$, the condition $0\in\partial_iF(\mathbf{x})$
implies $x_i=b_i=1$.
If $x_i<0$, then $\phi'(x_i-b_i)-2<0$, so this condition cannot hold.
If $x_i=0$, then $0\in\partial_iF(\mathbf{x})$ if and only if $b_i=0$.
Consequently, $0\in\partial F(\mathbf{x})$ if and only if
$\mathbf{x}=\mathbf{b}$, proving (i).

(ii) Next, direct computation gives
\begin{equation*}
\phi''(t)=\frac{4t^2(3-t^4)}{(1+t^4)^2},
\quad
\phi'''(t)=\frac{8t(t^8-12t^4+3)}{(1+t^4)^3}.
\end{equation*}
The nonzero critical points of $\phi''$ satisfy
$t^4=6\pm\sqrt{33}$, while those of $\phi'''$
satisfy $t^4\in\{1,15-4\sqrt{14},15+4\sqrt{14}\}$.
Evaluating the derivatives at these points yields
\begin{equation*}
\sup_{t\in\mathbb{R}}|\phi''(t)|<4,
\quad
\sup_{t\in\mathbb{R}}|\phi'''(t)|<9.
\end{equation*}
Since
\begin{equation*}
\nabla^2h(\mathbf{x})
=\operatorname{Diag}\big(\phi''(x_i-b_i)\big)_{i=1}^n,
\end{equation*}
we have
\begin{equation*}
\|\nabla^2h(\mathbf{x})\|
\leq \sup_{t\in\mathbb{R}}|\phi''(t)|<4.
\end{equation*}
Thus, $\nabla h$ is globally Lipschitz continuous
with constant $L<4$, proving (ii).

(iii) Fix $\mathbf{0}\neq\mathbf{x}\in\mathbb{R}^n$ and $\emptyset\neq\mathcal{S}\subseteq\operatorname{supp}(\mathbf{x})$.
Since $g$ is locally linear along $\mathcal{M}_{\mathbf{x}}$,
we have
\begin{equation*}
\nabla_{\mathcal{S}}^2F(\mathbf{z})
=\operatorname{Diag}\big(\phi''(z_i-b_i)\big)_{i\in\mathcal{S}}
\end{equation*}
for $\mathbf{z}\in\mathcal{M}_{\mathbf{x}}$
sufficiently close to $\mathbf{x}$.
By the mean value theorem, for any such
$\mathbf{z},\mathbf{w}$,
\begin{equation*}
\|\nabla_{\mathcal{S}}^2F(\mathbf{z})
-\nabla_{\mathcal{S}}^2F(\mathbf{w})\| \leq
\sup_{t\in\mathbb{R}}|\phi'''(t)|
\|\mathbf{z}-\mathbf{w}\|
\leq9\|\mathbf{z}-\mathbf{w}\|.
\end{equation*}
Hence,
$\operatorname{lip}_{\mathcal{M}_{\mathbf{x}}}
\nabla_{\mathcal{S}}^2F(\mathbf{x})\leq9$.

(iv) We next show that $F$ satisfies the KL inequality with exponent $3/4$ at $\overline{\mathbf{x}}$.
Choose a neighborhood $U$ of $\overline{\mathbf{x}}$
such that, for every $\mathbf{x}\in U$, $|x_i-b_i|<1$ and $F(\mathbf{x})<1$.

Consider first the case where $x_i\ge0$ for all $i$. Then
\begin{equation*}
F(\mathbf{x})
=\sum_{i=1}^n\phi(x_i-b_i)
\le\sum_{i=1}^n|x_i-b_i|^4
\le n\|\mathbf{x}-\mathbf{b}\|_\infty^4.
\end{equation*}
{For the derivative, we have}
\begin{equation*}
|\phi'(t)|
=\frac{4|t|^3}{1+t^4}
\ge 2|t|^3,
\qquad {|t|\le1.}
\end{equation*}
{Any coordinate with $x_i=0$ has $b_i=0$ in this neighborhood. If $\mathbf x\ne\mathbf b$, an index attaining $\|\mathbf x-\mathbf b\|_\infty$ therefore has $x_i>0$, and its partial subdifferential is $\{\phi'(x_i-b_i)\}$. Thus}
\begin{equation*}
\operatorname{dist}(0,\partial F(\mathbf{x}))
\ge 2\|\mathbf{x}-\mathbf{b}\|_\infty^3
\ge \frac{2}{n^{3/4}}F(\mathbf{x})^{3/4}.
\end{equation*}
Next, suppose that $x_i<0$ for some $i$.
Since $\mathbf{x}\in U$, we must have $b_i=0$.
Consequently,
\begin{equation*}
\partial_iF(\mathbf{x})
=\left\{
\frac{4x_i^3}{1+x_i^4}-2
\right\},
\end{equation*}
which implies
\begin{equation*}
\operatorname{dist}(0,\partial F(\mathbf{x}))
\ge2
\ge\frac{2}{n^{3/4}}F(\mathbf{x})^{3/4}.
\end{equation*}
Combining the two cases gives
\begin{equation*}
\operatorname{dist}(0,\partial F(\mathbf{x}))
\ge\frac{2}{n^{3/4}}
\big(F(\mathbf{x})-F(\overline{\mathbf{x}})\big)^{3/4}
\end{equation*}
for all $\mathbf{x}\in U$.
Hence, $F$ has the KL property at
$\overline{\mathbf{x}}$ with exponent $3/4$.

It remains to show that this exponent is exact.
Choose an index $j$ with $b_j=1$ and let
$\mathbf{x}(t)=\overline{\mathbf{x}}+t\mathbf{e}_j$
for $t>0$. Then
\begin{equation*}
F(\mathbf{x}(t))-F(\overline{\mathbf{x}})
=\log(1+t^4),\quad
\operatorname{dist}(0,\partial F(\mathbf{x}(t)))
=\frac{4t^3}{1+t^4}.
\end{equation*}
For any $\theta<3/4$, it follows that
\begin{equation*}
\frac{\operatorname{dist}
	(0,\partial F(\mathbf{x}(t)))}
{(F(\mathbf{x}(t))-F(\overline{\mathbf{x}}))^\theta}
=
\frac{4t^3}
{(1+t^4)[\log(1+t^4)]^\theta}\to0
\quad\text{as }t\downarrow0.
\end{equation*}
Thus, the KL inequality cannot hold with any
exponent $\theta<3/4$.
Therefore, the exact KL exponent of $F$ at
$\overline{\mathbf{x}}$ is $3/4$.
\end{proof}

\section[Sufficient Conditions for Boundedness]{{Sufficient Conditions for Boundedness}}\label{sec:boundness_sequence}
\begin{lemma}
\label{cor:L1_L2_ratio_boundedness}
Consider the $L_1/L_2$ regularized least-squares problem  (\ref{Problem}) with $h({\h x})=\frac12\|A\mathbf{x}-\mathbf{b}\|_2^2$.
Suppose $\ker(A)\neq\{\mathbf{0}\}$ and define
\begin{equation*}
\rho_A
:=
\min_{\mathbf{d}\in\ker(A)\setminus\{\mathbf{0}\}}
\frac{\|\mathbf{d}\|_1}{\|\mathbf{d}\|}.
\end{equation*}
If the following condition holds
$
F(\mathbf{x}^0) <
\frac{1}{2}\inf_{\mathbf{x}\in\mathbb{R}^n}\|A\mathbf{x}-\mathbf{b}\|^2+\lambda\rho_A,
$
then the sequence $\{\mathbf{x}^k\}$ generated by
Algorithm~\ref{alg:apgcn} is bounded.
\end{lemma}

\begin{proof}
Let $ h^*:=\inf_{\mathbf{x}\in\mathbb{R}^n}\frac{1}{2}\|A\mathbf{x}-\mathbf{b}\|^2.$
Suppose by contradiction that $\{\mathbf{x}^k\}$ is unbounded. Then there exists a subsequence, still denoted by $\{\mathbf{x}^k\}$ for simplicity, such that
\[
\lim_{k\to\infty}\|\mathbf{x}^k\|=+\infty.
\]
Since $\{F(\mathbf{x}^k)\}$ is non-increasing and
\begin{equation*}
\frac{1}{2}\|A\mathbf{x}^k-\mathbf{b}\|^2
\leq
F(\mathbf{x}^k)
\leq
F(\mathbf{x}^0),
\end{equation*}
the sequence $\{A\mathbf{x}^k\}$ is bounded.
By taking a further subsequence if necessary, we have
\[
\lim_{k\to\infty}\mathbf{x}^k/\|\mathbf{x}^k\|_2=\overline{\mathbf{d}}
\quad\text{with}\quad \|\overline{\mathbf{d}}\|=1.
\]
Since $\{A\mathbf{x}^k\}$ is bounded and
$\lim_{k\to\infty}\|\mathbf{x}^k\|=+\infty$, it follows that
$A\overline{\mathbf{d}}=\mathbf{0}$. Hence,
$\overline{\mathbf{d}}\in\ker(A)$ and $ \|\overline{\mathbf{d}}\|_1\geq\rho_A.$
Therefore,
\begin{equation*}
F(\mathbf{x}^0) \geq
\lim_{k\to\infty}F(\mathbf{x}^k)
\geq h^*+ \lambda\|\overline{\mathbf{d}}\|_1 \geq h^*+\lambda\rho_A,
\end{equation*}
which contradicts the assumption. Hence, $\{\mathbf{x}^k\}$ is bounded.
\end{proof}
%\begin{remark}
%(1) The quantity $\rho_A$ was introduced in \cite{zeng2021analysis} to study the existence of optimal solutions to linearly constrained $L_1/L_2$ minimization problems. Related null-space quantities were also employed in \cite{tao2023study} to establish the existence of global solutions for $L_1/L_2$ minimization models. However, neither work relates $\rho_A$ to the boundedness of the sequence $\{\mathbf{x}^k\}$.
%
%(2) The quantity $\rho_A$ is also highly related to the spherical section property \cite{vavasis2009spherical,zhang2013theory}. Recall that $\ker(A)$ is said to have the $s$-spherical section
%property if $\rho_A \geqslant \sqrt{m/s}$.
%In this case, the more explicit condition
%\begin{equation*}
%F(\mathbf{x}^0) < h^*+\lambda\sqrt{\frac{m}{s}}
%\end{equation*}
%is sufficient for the boundedness of $\{\mathbf{x}^k\}$.
%Moreover, it was shown in \cite[Theorem 3.1]{zhang2013theory} that if $A$ is a random matrix with independent and identically distributed (i.i.d.) standard Gaussian entries, $\ker(A)$ has the $s$-spherical section property with $s=c_1(\log(n/m)+1)$ with probability at least $1-e^{-c_0(n-m)}$, where $c_0$ and $c_1$ are positive constants independent of $m$ and $n$.
%\end{remark}
\begin{remark}
The quantity $\rho_A$ is related to null-space quantities used in
existence analyses
\cite{zeng2021analysis,tao2023study}.
{If $\ker(A)$ has the $s$-spherical section property, i.e.,
\[
\frac{\|\mathbf d\|_1^2}{\|\mathbf d\|_2^2}\ge\frac ms
\quad\text{for every }\mathbf d\in\ker(A)\setminus\{0\}{,}
\]
then $\rho_A\ge\sqrt{m/s}$.} Hence the simpler condition
$
F(\mathbf{x}^0)<h^*+\lambda\sqrt{\frac{m}{s}}
$
also guarantees boundedness. For matrices with independent standard Gaussian entries and $m<n$, the null space
satisfies this property with high probability for
$s=c_0(\log(n/m)+1)$; see \cite[Theorem~3.1]{zhang2013theory}.
\end{remark}
\begin{lemma}
\label{cor:L1_muL2_logistic_boundedness}
Consider the $L_1-\mu L_2$ regularized logistic regression problem (\ref{eq:intro_sparse_logistic}) with $0<\mu\leq1$.
Suppose
\begin{equation}
\left\{
(\mathbf{d},v)\in\mathbb{R}^n\times\mathbb{R}:
b_i(\mathbf{a}_i^\top\mathbf{d}+v)\geq0, i\in[m],
\|\mathbf{d}\|_1-\mu\|\mathbf{d}\|=0
\right\}
=
\{(\mathbf{0},0)\}.
\label{eq:L1_muL2_logistic_condition}
\end{equation}
Then the sequence $\{(\mathbf{z}^k,u^k)\}$ generated by
Algorithm~\ref{alg:apgcn} is bounded.
Moreover, let
$\mathcal{I}_+:=\{i\in[m]:b_i=1\}$ and $\mathcal{I}_-:=\{i\in[m]:b_i=-1\}$.

\begin{itemize}
\item[(i)] If $0<\mu<1$, condition
\eqref{eq:L1_muL2_logistic_condition} is equivalent to
$\mathcal{I}_+\neq\emptyset$ and
$\mathcal{I}_-\neq\emptyset$.
\item[(ii)] If $\mu=1$, condition
\eqref{eq:L1_muL2_logistic_condition} is equivalent to
$\mathcal{I}_+\neq\emptyset$ and
$\mathcal{I}_-\neq\emptyset$, together with
\begin{equation*}
	\min_{i\in\mathcal{I}_+}a_{ij} < \max_{i\in\mathcal{I}_-}a_{ij},
	\quad
	\min_{i\in\mathcal{I}_-}a_{ij} <\max_{i\in\mathcal{I}_+}a_{ij},
	\quad \forall j\in[n],
	\label{eq:L1_muL2_logistic_interval_condition}
\end{equation*}
where $a_{ij}$ denotes the $(i,j)$-th component of $A$.
\end{itemize}
\end{lemma}

\begin{proof}
Suppose by contradiction that $\{(\mathbf{z}^k,u^k)\}$ is unbounded.
Then there exists a subsequence, still denoted by
$\{(\mathbf{z}^k,u^k)\}$ for simplicity, such that $\lim_{k\to\infty}$ $\|(\mathbf{z}^k,u^k)\|=+\infty.$
Since $0<\mu\leq1$, we have
$ \|\mathbf{z}\|_1-\mu\|\mathbf{z}\|\geq(1-\mu)\|\mathbf{z}\|\geq0.$
Hence, by the monotonicity of
$\{F(\mathbf{z}^k,u^k)\}$, we obtain $
h(\mathbf{z}^k,u^k) \leq F(\mathbf{z}^k,u^k) \leq F(\mathbf{z}^0,u^0) $ and
$ \lambda ( \|\mathbf{z}^k\|_1-\mu\|\mathbf{z}^k\|)\leq F(\mathbf{z}^0,u^0).$
Let $t_k:=\|(\mathbf{z}^k,u^k)\|,(\mathbf{d}^k,v^k):=(\mathbf{z}^k,u^k)/t_k.$ By taking a further subsequence if necessary, we have
$ (\mathbf{d}^k,v^k) \to (\overline{\mathbf{d}},\overline{v})$ with $
\|(\overline{\mathbf{d}},\overline{v})\|=1.$
We first show that
\begin{equation*}
b_i(\mathbf{a}_i^\top\overline{\mathbf{d}}+\overline{v})
\geq0, \quad i=1,\ldots,m.
\end{equation*}
Suppose otherwise that there exists some $i$ such that
$ b_i(\mathbf{a}_i^\top\overline{\mathbf{d}}+\overline{v})<0.$
Since
\begin{equation*}
\lim_{k\to\infty}\frac{
	b_i(\mathbf{a}_i^\top\mathbf{z}^k+u^k)}
{t_k}
=
\lim_{k\to\infty}b_i(\mathbf{a}_i^\top\mathbf{d}^k+v^k)
=
b_i(\mathbf{a}_i^\top\overline{\mathbf{d}}+\overline{v})<0,
\end{equation*}
we have
$\lim_{k\to\infty}b_i(\mathbf{a}_i^\top\mathbf{z}^k+u^k)=-\infty$.
Consequently,
$\lim_{k\to\infty}\log( 1+e^{-b_i(\mathbf{a}_i^\top\mathbf{z}^k+u^k)})$ $=+\infty,
$
which contradicts the boundedness of
$\{h(\mathbf{z}^k,u^k)\}$.

Moreover, note that
$ \lambda ( \|\mathbf{z}^k\|_1-\mu\|\mathbf{z}^k\|_2)
= \lambda t_k( \|\mathbf{d}^k\|_1-\mu\|\mathbf{d}^k\|_2).$
Since the left-hand side is bounded and $\lim_{k\to\infty}t_k=+\infty$, we have
\begin{equation*}
\lim_{k\to\infty}\|\mathbf{d}^k\|_1-\mu\|\mathbf{d}^k\| =\|\overline{\mathbf{d}}\|_1 - \mu\|\overline{\mathbf{d}}\| =0.
\end{equation*}
Therefore, by assumption it follows that
$(\overline{\mathbf{d}},\overline{v})=(\mathbf{0},0)$,
which contradicts
$\|(\overline{\mathbf{d}},\overline{v})\|_2=1$.
Hence, $\{(\mathbf{z}^k,u^k)\}$ is bounded.

We next establish the stated equivalences.

(i) Suppose $0<\mu<1$. Since $\|\mathbf{d}\|_1-\mu\|\mathbf{d}\|_2
\geq(1-\mu)\|\mathbf{d}\|_2$, condition \eqref{eq:L1_muL2_logistic_condition} is equivalent to
\begin{equation*}
\left\{ v\in\mathbb{R}: b_i v\geq0,\ i\in[m]\right\}=\{0\},
\end{equation*}
which is equivalent to $\mathcal{I}_+\neq\emptyset$ and
$\mathcal{I}_-\neq\emptyset$.

(ii) Suppose $\mu = 1$. In this case $\|\mathbf{d}\|_1 - \|\mathbf{d}\| = 0$ holds if and only if $\|\mathbf{d}\|_0 \leq 1$. If $\mathbf{d} = \mathbf{0}$, {the argument in (i) shows that $v=0$ is the only feasible intercept direction if and only if} $\mathcal{I}_+\neq\emptyset$ and
$\mathcal{I}_-\neq\emptyset$. It remains to consider the case where $\|\mathbf{d}\|_0 = 1$. {Write $\mathbf{d}=t\mathbf{e}_j$ for some $j\in[n]$ and $t\neq0$.} Then condition \eqref{eq:L1_muL2_logistic_condition} holds if and only if for any $j\in[n]$, $t\neq0$ and $v\in\mathbb{R}$, there exists $i\in[m]$ such that $b_i(ta_{ij} + v) <0$, i.e.
$
\min_{i\in\mathcal{I}_+}a_{ij} < \max_{i\in\mathcal{I}_-}a_{ij},
\;
\min_{i\in\mathcal{I}_-}a_{ij} <\max_{i\in\mathcal{I}_+}a_{ij},
\; \forall j\in[n],
$
which proves the result.
\end{proof}

\begin{remark}
Lemma \ref{cor:L1_muL2_logistic_boundedness} gives a geometric condition
that depends on the value of the regularization parameter.
First, when $0<\mu<1$, the condition is equivalent to requiring
that both positive and negative classes are present in the data set.
Second, when $\mu=1$, the condition requires that the two classes
cannot be weakly separated based on a single
feature~\cite[Sec.~2.5.1]{BoydVandenberghe2004}.
\end{remark}


\begin{thebibliography}{99}

\bibitem[Attouch and Bolte(2009)]{attouch2009convergence}
Attouch, H., Bolte, J.: On the convergence of the proximal algorithm for nonsmooth functions involving analytic features.
Mathematical Programming \textbf{116}(1), 5--16 (2009).
\url{https://doi.org/10.1007/s10107-007-0133-5}

\bibitem[Attouch et al.(2010)]{attouch2010proximal}
Attouch, H., Bolte, J., Redont, P., Soubeyran, A.: Proximal alternating minimization and projection methods for nonconvex problems: An approach based on the Kurdyka--{\L ojasiewicz} inequality.
Mathematics of Operations Research \textbf{35}(2), 438--457 (2010).
\url{https://doi.org/10.1287/moor.1100.0449}

\bibitem[Bareilles, Iutzeler, and Malick(2023)]{bareilles2023newton}
Bareilles, G., Iutzeler, F., Malick, J.: Newton acceleration on manifolds identified by proximal gradient methods.
Mathematical Programming \textbf{200}(1), 37--70 (2023).
\url{https://doi.org/10.1007/s10107-022-01873-w}

\bibitem[Beck(2017)]{beck2017first}
Beck, A.: First-order methods in optimization.
SIAM (2017).
\url{https://doi.org/10.1137/1.9781611974997}

\bibitem[Birgin et al.(2017)]{birgin2017worst}
Birgin, E.G., Gardenghi, J.L., Mart\'inez, J.M., Santos, S.A., Toint, P.L.: Worst-case evaluation complexity for unconstrained nonlinear optimization using high-order regularized models.
Mathematical Programming \textbf{163}(1--2), 359--368 (2017).
\url{https://doi.org/10.1007/s10107-016-1065-8}

\bibitem[Blumensath and Davies(2009)]{blumensath2009iterative}
Blumensath, T., Davies, M.E.: Iterative Hard Thresholding for Compressed Sensing.
Applied and Computational Harmonic Analysis \textbf{27}(3), 265--274 (2009).
\url{https://doi.org/10.1016/j.acha.2009.04.002}

% \bibitem[Bolte et al.(2014)]{bolte2014proximal}
% J. Bolte, S. Sabach, and M. Teboulle.
% \textit{Proximal alternating linearized minimization for nonconvex and nonsmooth problems}.
% Mathematical Programming, 146.1--2 (2014), pp.~459--494.

\bibitem[R. I. Bo\c{t}, E. R. Csetnek, D.-K. Nguyen(2019)]{BCN19}
Bo\c{t}, R.I., Csetnek, E.R., Nguyen, D.-K.: A proximal minimization algorithm for structured nonconvex and nonsmooth problems.
SIAM Journal on Optimization \textbf{29}, 1300--1328 (2019).
\url{https://doi.org/10.1137/18M1190689}

\bibitem[Bo{\c{t}}, Li, and Tao(2025)]{boct2025full}
Bo{\c{t}}, R.I., Li, G.Y., Tao, M.: Full splitting algorithms for fractional programs with structured numerators and denominators.
SIAM Journal on Optimization \textbf{35}(4), 2623--2653 (2025).
\url{https://doi.org/10.1137/23M1626591}

\bibitem[S.~Boyd and L.~Vandenberghe(2004)]{BoydVandenberghe2004}
Boyd, S., Vandenberghe, L.: Convex Optimization.
Cambridge University Press, Cambridge (2004).
\url{https://doi.org/10.1017/CBO9780511804441}

\bibitem[Carmon and Duchi(2019)]{carmon2019gradient}
Carmon, Y., Duchi, J.: Gradient descent finds the cubic-regularized nonconvex Newton step.
SIAM Journal on Optimization \textbf{29}(3), 2146--2178 (2019).
\url{https://doi.org/10.1137/17M1113898}

\bibitem[Cartis, Gould, and Toint(2011a)]{cartis2011adaptive}
Cartis, C., Gould, N.I.M., Toint, P.L.: Adaptive cubic regularisation methods for unconstrained optimization. Part I: motivation, convergence and numerical results.
Mathematical Programming \textbf{127}(2), 245--295 (2011).
\url{https://doi.org/10.1007/s10107-009-0286-5}

\bibitem[Cartis, Gould, and Toint(2011b)]{cartis2011adaptive2}
Cartis, C., Gould, N.I.M., Toint, P.L.: Adaptive cubic regularisation methods for unconstrained optimization. Part II: worst-case function- and derivative-evaluation complexity.
Mathematical Programming \textbf{130}(2), 295--319 (2011).
\url{https://doi.org/10.1007/s10107-009-0337-y}

\bibitem[Chen, Curtis, and Robinson(2017)]{chen2017reduced}
Chen, T.Y., Curtis, F.E., Robinson, D.P.: A reduced-space algorithm for minimizing $\ell_1$-regularized convex functions.
SIAM Journal on Optimization \textbf{27}(3), 1583--1610 (2017).
\url{https://doi.org/10.1137/16M1062259}

%\bibitem[Chen, Jiang, Lin, and Zhang(2022)]{chen2022accelerating}
%X. Chen, B. Jiang, T. Y. Lin, and S. Z. Zhang.
%\textit{Accelerating adaptive cubic regularization of Newton's method via random sampling}.
%Journal of Machine Learning Research, 23.90 (2022), pp.~1--38.

\bibitem[Curtis, Dai, and Robinson(2022)]{curtis2022subspace}
Curtis, F.E., Dai, Y., Robinson, D.P.: A subspace acceleration method for minimization involving a group sparsity-inducing regularizer.
SIAM Journal on Optimization \textbf{32}(2), 545--572 (2022).
\url{https://doi.org/10.1137/21M1411111}

\bibitem[Gao, Yue, and Ng(2022)]{gao2022approximate}
Gao, Y., Yue, M.-C., Ng, M.K.: Approximate secular equations for the cubic regularization subproblem.
In: Advances in Neural Information Processing Systems \textbf{35} (2022).
\url{https://doi.org/10.52202/068431-1036}

\bibitem[Gong et al.(2013)]{gong2013general}
Gong, P., Zhang, C., Lu, Z., Huang, J., Ye, J.: A General Iterative Shrinkage and Thresholding Algorithm for Non-convex Regularized Optimization Problems.
In: Proceedings of the 30th International Conference on Machine Learning, Proceedings of Machine Learning Research \textbf{28}(2), 37--45 (2013).

\bibitem[Hare and Lewis(2004)]{hare2004identifying}
Hare, W.L., Lewis, A.S.: Identifying active constraints via partial smoothness and prox-regularity.
Journal of Convex Analysis \textbf{11}(2), 251--266 (2004).
\url{https://doi.org/10.68381/jca11017}

\bibitem[Huang et al.(2010)]{huang2010predicting}
Huang, L., Jia, J., Yu, B., Chun, B.-G., Maniatis, P., Naik, M.: Predicting execution time of computer programs using sparse polynomial regression.
Advances in Neural Information Processing Systems \textbf{23}, 883--891 (2010).
\url{https://papers.nips.cc/paper_files/paper/2010/hash/995665640dc319973d3173a74a03860c-Abstract.html}


\bibitem[Jia et al.(2022)]{jia2022solving}
Jia, X., Liang, X., Shen, C., Zhang, L.-H.: Solving the cubic regularization model by a nested restarting Lanczos method.
SIAM Journal on Matrix Analysis and Applications \textbf{43}(2), 812--839 (2022).
\url{https://doi.org/10.1137/21M1436324}

\bibitem[Jiang, Yue, and Zhou(2021)]{jiang2021accelerated}
Jiang, R.J., Yue, M.-C., Zhou, Z.S.: An accelerated first-order method with complexity analysis for solving cubic regularization subproblems.
Computational Optimization and Applications \textbf{79}(2), 471--506 (2021).
\url{https://doi.org/10.1007/s10589-021-00274-7}

\bibitem[Jiang, Zhou, and Zhou(2022)]{jiang2022cubic}
Jiang, R.J., Zhou, Z.S., Zhou, Z.R.: Cubic regularization methods with second-order complexity guarantee based on a new subproblem reformulation.
Journal of the Operations Research Society of China \textbf{10}(3), 471--506 (2022).
\url{https://doi.org/10.1007/s40305-022-00398-5}

\bibitem[Lee(2023)]{lee2023accelerating}
Lee, C.-P.: Accelerating inexact successive quadratic approximation for regularized optimization through manifold identification.
Mathematical Programming \textbf{201}, 599--633 (2023).
\url{https://doi.org/10.1007/s10107-022-01916-2}


\bibitem[Lewis(2002)]{lewis2002active}
Lewis, A.S.: Active sets, nonsmoothness, and sensitivity.
SIAM Journal on Optimization \textbf{13}(3), 702--725 (2002).
\url{https://doi.org/10.1137/S1052623401387623}

\bibitem[Lewis and Wright(2016)]{lewis2016proximal}
Lewis, A.S., Wright, S.J.: A proximal method for composite minimization.
Mathematical Programming \textbf{158}, 501--546 (2016).
\url{https://doi.org/10.1007/s10107-015-0943-9}

\bibitem[Lewis and Wylie(2021)]{lewis2021active}
Lewis, A.S., Wylie, C.: Active-set Newton methods and partial smoothness.
Mathematics of Operations Research \textbf{46}(2), 712--725 (2021).
\url{https://doi.org/10.1287/moor.2020.1075}

\bibitem[Lewis and Zhang(2013)]{lewis2013partial}
Lewis, A.S., Zhang, S.: Partial smoothness, tilt stability, and generalized Hessians.
SIAM Journal on Optimization \textbf{23}(1), 74--94 (2013).
\url{https://doi.org/10.1137/110852103}


\bibitem[G. Li and T.~K. Pong(2015)]{LiPong15}
Li, G., Pong, T.K.: Global convergence of splitting methods for nonconvex composite optimization.
SIAM Journal on Optimization \textbf{25}, 2434--2460 (2015).
\url{https://doi.org/10.1137/140998135}

\bibitem[Liang, Fadili, and Peyr\'e(2017)]{liang2017activity}
Liang, J., Fadili, J., Peyr\'e, G.: Activity identification and local linear convergence of forward--backward-type methods.
SIAM Journal on Optimization \textbf{27}(1), 408--437 (2017).
\url{https://doi.org/10.1137/16M106340X}

\bibitem[Lieder(2020)]{lieder2020solving}
Lieder, F.: Solving large-scale cubic regularization by a generalized eigenvalue problem.
SIAM Journal on Optimization \textbf{30}(4), 3345--3358 (2020).
\url{https://doi.org/10.1137/19M1291388}

\bibitem[Lin et al.(2024)]{lin2024highly}
Lin, M., Yuan, Y., Sun, D., Toh, K.-C.: A Highly Efficient Algorithm for Solving Exclusive Lasso Problems.
Optimization Methods and Software \textbf{39}(3), 489--518 (2024).
\url{https://doi.org/10.1080/10556788.2023.2253356}

\bibitem[Liu, Deng, and He(2024)]{liu2024jaccard}
Liu, R., Deng, G., He, H.: Generalized Jaccard feature screening for ultra-high dimensional survival data.
AIMS Mathematics \textbf{9}(10), 27607--27626 (2024).
\url{https://doi.org/10.3934/math.20241341}

\bibitem[Lou and Yan(2018)]{lou2018fast}
Lou, Y., Yan, M.: Fast $L_1$--$L_2$ Minimization via a Proximal Operator.
Journal of Scientific Computing \textbf{74}(2), 767--785 (2018).
\url{https://doi.org/10.1007/s10915-017-0463-2}

\bibitem[Nesterov and Polyak(2006)]{nesterov2006cubic}
Nesterov, Y., Polyak, B.: Cubic regularization of Newton method and its global performance.
Mathematical Programming \textbf{108}(1), 177--205 (2006).
\url{https://doi.org/10.1007/s10107-006-0706-8}

\bibitem[Nutini, Schmidt, and Hare(2019)]{nutini2019active}
Nutini, J., Schmidt, M., Hare, W.: Active-set complexity of proximal gradient: How long does it take to find the sparsity pattern?
Optimization Letters \textbf{13}, 645--655 (2019).
\url{https://doi.org/10.1007/s11590-018-1325-z}

\bibitem[Qin and Lou(2019)]{qin20191}
Qin, J., Lou, Y.F.: $\ell_{1-2}$ regularized logistic regression.
In: Proceedings of the 53rd Asilomar Conference on Signals, Systems, and Computers, pp.~779--783 (2019).
\url{https://doi.org/10.1109/IEEECONF44664.2019.9048830}

\bibitem[Rockafellar and Wets(1998)]{rockafellar1998variational}
Rockafellar, R.T., Wets, R.J.B.: Variational Analysis.
Springer (1998).


\bibitem[Simon et al.(2013)]{simon2013sparse}
Simon, N., Friedman, J., Hastie, T., Tibshirani, R.: A sparse-group lasso.
Journal of Computational and Graphical Statistics \textbf{22}(2), 231--245 (2013).
\url{https://doi.org/10.1080/10618600.2012.681250}


\bibitem[Tao(2022)]{Tao20}
Tao, M.: Minimization of ${L}_1$ over ${L}_2$ for sparse signal recovery with convergence guarantee.
SIAM Journal on Scientific Computing \textbf{44}(2), A770--A797 (2022).
\url{https://doi.org/10.1137/20M136801X}

\bibitem[Tao and Zhang(2023)]{tao2023study}
Tao, M., Zhang, X.-P.: Study on $L_1$ over $L_2$ minimization for nonnegative signal recovery.
Journal of Scientific Computing \textbf{95}(3), Article~94 (2023).
\url{https://doi.org/10.1007/s10915-023-02225-2}

\bibitem[Tao and Zhang(2026)]{tao2026complexity}
Tao, M., Zhang, X.-P.: A complexity analysis framework for active manifold identification with applications to $L_0$ and $L_p$ regularization models.
Applied Mathematics \& Optimization \textbf{93}, Article~50 (2026).
\url{https://doi.org/10.1007/s00245-026-10410-6}

\bibitem[Tao, Zhang, and Xia(2024)]{tao2024partlysmooth}
Tao, M., Zhang, X.P., Xia, Z.H.: On partial smoothness, activity identification and faster algorithms of $L_1$ over $L_2$ minimization.
IEEE Transactions on Signal Processing \textbf{72}, 2874--2889 (2024).
\url{https://doi.org/10.1109/TSP.2024.3404250}

\bibitem[Themelis, Stella, and Patrinos(2018)]{themelis2018forward}
Themelis, A., Stella, L., Patrinos, P.: Forward-backward envelope for the sum of two nonconvex functions: Further properties and nonmonotone linesearch algorithms.
SIAM Journal on Optimization \textbf{28}(3), 2274--2303 (2018).
\url{https://doi.org/10.1137/16M1080240}

\bibitem[Tibshirani(1996)]{tibshirani1996regression}
Tibshirani, R.: Regression Shrinkage and Selection via the Lasso.
Journal of the Royal Statistical Society: Series B \textbf{58}(1), 267--288 (1996).
\url{https://doi.org/10.1111/j.2517-6161.1996.tb02080.x}

\bibitem[Wang, Yang, and Zhu(2024)]{wang2024alternating}
Wang, H., Yang, X.Y., Zhu, Y.C.: Alternating iteratively reweighted $\ell_1$ and subspace Newton algorithms for nonconvex sparse optimization.
arXiv preprint arXiv:2407.17216 (2024).

\bibitem[Wu, Pan, and Yang(2023)]{wu2023regularized}
Wu, Y.Q., Pan, S.H., Yang, X.Q.: A Regularized Newton Method for $\ell_q$-Norm Composite Optimization Problems.
SIAM Journal on Optimization \textbf{33}(3), 1676--1706 (2023).
\url{https://doi.org/10.1137/22M1482822}

\bibitem[Xu et al.(2012)]{xu2012half}
Xu, Z., Chang, X., Xu, F., Zhang, H.: $L_{1/2}$ Regularization: A Thresholding Representation Theory and a Fast Solver.
IEEE Transactions on Neural Networks and Learning Systems \textbf{23}(7), 1013--1027 (2012).
\url{https://doi.org/10.1109/TNNLS.2012.2197412}

\bibitem[Zeng, Yu, and Pong(2021)]{zeng2021analysis}
Zeng, L., Yu, P., Pong, T.K.: Analysis and algorithms for some compressed sensing models based on $L_1/L_2$ minimization.
SIAM Journal on Optimization \textbf{31}(2), 1576--1603 (2021).
\url{https://doi.org/10.1137/20M1355380}

\bibitem[Zhang(2013)]{zhang2013theory}
Zhang, Y.: Theory of compressive sensing via $\ell_1$-minimization: A non-RIP analysis and extensions.
Journal of the Operations Research Society of China \textbf{1}(1), 79--105 (2013).
\url{https://doi.org/10.1007/s40305-013-0010-2}

\bibitem[Zhang and Xin(2017)]{zhang2017minimization}
Zhang, S., Xin, J.: Minimization of Transformed $L_1$ Penalty: Closed Form Representation and Iterative Thresholding Algorithms.
Communications in Mathematical Sciences \textbf{15}(2), 511--537 (2017).
\url{https://doi.org/10.4310/CMS.2017.v15.n2.a9}

\bibitem[Zhao, Doikov, and Lucchi(2025)]{zhao2024cubic}
Zhao, J., Doikov, N., Lucchi, A.: Cubic regularized subspace Newton for non-convex optimization.
Proceedings of the 28th International Conference on Artificial Intelligence and Statistics, Proceedings of Machine Learning Research \textbf{258}, 811--819 (2025).


\bibitem[Zou and Hastie(2005)]{zou2005regularization}
Zou, H., Hastie, T.: Regularization and Variable Selection via the Elastic Net.
Journal of the Royal Statistical Society: Series B \textbf{67}(2), 301--320 (2005).
\url{https://doi.org/10.1111/j.1467-9868.2005.00503.x}

\end{thebibliography}
\end{document}